%% file: main.tex
\documentclass[a4paper,11pt,oneside]{article}

\input{header.tex}
\usepackage{siunitx}

\newcommand\plot[1]{\let\frame\relax
	\frame{\includegraphics[clip,trim=0 220 0 60,width=8cm]{#1}}}
\newcommand{\dx}{\text{ d}x}

\newcommand{\nn}{\nonumber}

\usepackage{authblk}

\usepackage[top=2cm, bottom=2cm, left=2cm, right=2cm]{geometry}

\theoremstyle{plain} 
\newtheorem{theorem}{Theorem}[section]

\newtheorem{Remark}[theorem]{Remark}
\newtheorem{Definition}[theorem]{Definition}

\newtheorem{Problem}[theorem]{Problem}

\theoremstyle{definition} %

\theoremstyle{remark} %

\usepackage{newunicodechar}
\newunicodechar{ℝ}{\mathbb{R}}
\newunicodechar{Ψ}{\Psi}
\newunicodechar{α}{\alpha}
\newunicodechar{∇}{\nabla}
\newunicodechar{σ}{\sigma}
\newunicodechar{ν}{\nu}
\newunicodechar{ρ}{\rho}
\newunicodechar{ϑ}{\vartheta}
\newunicodechar{Γ}{\Gamma}

\newcommand{\Om}{\Omega}

\usepackage{xcolor}

\newcommand{\black}[1]{{\color{black}{#1}}}
\newcommand{\new}[1]{{\color{black}{#1}}}
\newcommand{\T}{{\vartheta}}

\begin{document}
	
	\title{
		Modeling and multigoal-oriented a posteriori error control 
		for heated material processing 
		using a generalized Boussinesq model
	}
	
	\author[1,2]{S. Beuchler}
	\author[1,2]{B. Endtmayer}
	\author[1]{J. Lankeit}
	\author[1,2]{T. Wick}

	\affil[1]{Leibniz Universit\"at Hannover, Institut f\"ur Angewandte
		Mathematik, Welfengarten 1, 30167 Hannover, Germany}
	
	\affil[2]{Cluster of Excellence PhoenixD (Photonics, Optics, and
		Engineering -- Innovation Across Disciplines), Leibniz Universit\"at Hannover, Germany}

	\date{}
	
	\maketitle

	\begin{abstract}
		 In this work, we develop a posteriori error control for a generalized 
		Boussinesq model in which thermal 
		conductivity and viscosity are temperature-dependent. Therein, 
		the stationary Navier-Stokes equations 
		are coupled with a stationary heat equation. 
		The coupled problem is modeled and solved
		in a monolithic fashion. The focus is on multigoal-oriented error estimation 
		with the dual-weighted residual method in which an adjoint problem 
		is utilized to obtain sensitivity measures with respect to several 
		goal functionals. The error localization is achieved with the help
		of a partition-of-unity in a weak formulation, which is specifically 
		convenient for coupled problems as we have at hand. The error indicators 
		are used to employ adaptive algorithms, which are substantiated with 
		several numerical tests such as one benchmark and two further experiments that 
		are motivated from laser material processing.
		Therein, error reductions and effectivity indices are consulted to 
		establish the robustness and efficiency of our framework.
		\textbf{Keywords: } Boussinesq; finite elements; 
		multigoal error control; partition-of-unity dual-weighted residuals; Y-beam splitter\\
		\textbf{MSC2020: 76M10, 76D05, 65N30, 65N50} 
	\end{abstract}

\section{Introduction}
\label{sec_intro}
This work considers a nonlinear coupled fluid flow heat system. 
Fluid flow is described by the 
incompressible Navier-Stokes equations \cite{Temam2001,Ga11} (for important
numerical developments, we refer to 
\cite{GiRa1986,GloTa89,GLOWINSKI20033,GlPe87,BristeauGlowinskiPeriaux1987,Ra00,Tu99,TuRiHrGl05,HeRa90}) 
and the heat distribution 
by an advection-diffusion equation. The resulting PDE (partial differential equation)
system is known as Boussinesq model \cite{DrRe81}. This model has been 
widely applied in various fields such as climate modeling \cite{Etling2008} 
or earth mantle convection problems \cite{Kronbichler2012}. 
Furthermore, the Boussinesq equation can serve as a sub-model within 
laser material processing \cite{OTTO201035} in wave guide modeling (e.g., \cite{Pae17,Chenetal18}) 
where heated material starts to flow 
due to local heat sources. A mathematical analysis of the stationary 
model that serves as our point of departure was done in \cite{LoBo96}.

The objective of this work is to design a robust and efficient framework
using adaptive finite elements for the numerical discretization of the 
Boussinesq system proposed in \cite{LoBo96}. 
Specifically, we derive multigoal a posteriori error estimates with respect 
to one or several quantities of interest \cite{Ha08,EnLaWi18,EnLaNeWiWo20}. This is \black{intriguing} since we deal 
with a coupled system of partial differential equations in which various 
parts of the solution might be of interest simultaneously. 
We notice that related results of coupling the stationary Navier-Stokes equations 
to the heat equation are published in some conference proceedings and the PhD thesis 
of the second author \cite{Endt21,BeuEndtWi21}. Moreover, there is only one other study 
\cite{AhEndtSteiWi22} in which
this multigoal-framework was applied so far to a nonlinear coupled system.  

In more detail, we formulate an optimization problem in which the discretization 
error measured in the goal functional is minimized with respect to a constraint. This 
constraint is nothing else than the PDE problem itself. For a very detailed 
description we refer the reader to the introduction of \cite{AhEndtSteiWi22}.
The resulting optimality system consists of the primal problem (the PDE, here 
the Boussinesq model) and a linear adjoint problem \cite{BeRa01,BaRa03}.
These results allow to design error 
identities and estimators for model errors (\cite{OdPr02,BraackErn02}), 
discretization and linearization
errors \cite{RanVi2013,EnLaWi18,EndtLaWi20_smart}. 
In this work, we consider discretization and linearization errors.
In order to use the error estimators for local mesh adaptivity we localize them to single 
mesh elements using a partition-of-unity localization \cite{RiWi15_dwr}.
This allows us to employ the algorithms from \cite{EndtLaWi20} and to apply them 
in this work to the Boussinseq system. For verification, we use one benchmark and 
we design two novel prototype experiments. Therein, we study error reductions and 
effectivity indices.

The outline of this paper is as follows: In Section \ref{sec_NSE_heat}, 
we explain our problem statement, derive the weak form and briefly 
explain the finite element discretization. Next, in Section \ref{sec_goal},
goal oriented adaptivity including multigoal estimates are addressed.
Then, in Section \ref{sec_tests}, we conduct three numerical tests in order 
to substantiate our algorithmic developments. Our work is summarized 
in Section \ref{sec_conclusions}.

\section{Boussinesq model: coupling Navier-Stokes to the heat equation}
\label{sec_NSE_heat}
Let $d=2$ (we notice that $d=3$ is possible as well) be the problem dimension and let $\Omega\subset\mathbb{R}^d$ 
be a bounded domain with boundary 
$\partial\Omega$.
For flow boundary conditions, $\partial\Omega$ is decomposed into non-overlapping parts
$\Gamma_{vD}$ and $\Gamma_{vN}:= \partial\Omega\setminus\Gamma_{vD}$,
where $\Gamma_{vD}$ indicate homogeneous 
or inhomogenous Dirichlet conditions, respectively, 
and $\Gamma_{vN}$ homogeneous 
or inhomogenous Neumann conditions, respectively. 
For the temperature equation, 
$\partial\Omega$ is decomposed into non-overlapping parts
$\Gamma_{\T D}$ and $\Gamma_{\T N}:= \partial\Omega\setminus\Gamma_{\T D}$.
We note
that for each numerical example, we specify the boundaries
separately.
Moreover, we denote the $L^2$ scalar product in $L^2(\Om;ℝ)$, $L^2(\Om,ℝ^d)$ or $L^2(\Om;ℝ^{d\times d})$ by $(\cdot,\cdot)$ 
and the $L^2$ scalar product over a boundary by $\langle \cdot, \cdot\rangle$.

\subsection{Model in strong form}
\label{sec_model_strong_form}
In this section, the strong form of the governing model is formulated.

\paragraph{Parameters and constitutive laws}
Let $\alpha \in \mathbb{R}$ be related to 
the coefficient of volume expansion, $g:\Omega\to\mathbb{R}^d$ be 
the external forces (for instance gravity) of the flow, 
$f:\Omega\to[0,\infty)$ 
be a heat source, and $k\colon ℝ\to (0,\infty)$ be
the thermal conductivity. 
Moreover, let the Cauchy stress tensor be given by
\begin{equation}\label{def:sigma}
\sigma:= \sigma(v,p,\T) = -pI + \rho\nu(\T) (\nabla v+\nabla v^T),
\end{equation}
where $I$ is the identity matrix and $\rho >0$ is the density.
The temperature-dependent kinematic viscosity is given by
\begin{equation}
\label{def:ny}
\nu(\T):=  \nu_0 e^{\frac{E_A}{R\T}},
\end{equation}
where $E_A>0$ and $\nu_0>0$ are material constants and $R>0$
is the universal gas constant. Specifically, \eqref{def:ny}
is the Arrhenius equation; see \cite{Arrhenius+1889+226+248,Arrhenius+1889+96+116} for chemical reactions and \cite{de1913relation,raman1923theory,andrade1934lviii,ward1937viscosity,haj2014contribution} for viscosity.

\paragraph{Strong form}
With these definitions at hand, our problem statement reads:
Find vector-valued velocities $v:\Omega\to\mathbb{R}^d$,
a scalar-valued pressure $p:\Omega\to\mathbb{R}$, and a 
scalar-valued temperature $\T:\Omega\to(0,\infty)$ such that
\begin{align}
(\rho v \cdot \nabla) v - \nabla\cdot \sigma  -\alpha \T g  =& 0   \qquad\text{in } \Omega,\nonumber\\ 
\nabla \cdot v =&  0   \qquad\text{in } \Omega,  \label{Eqn: Model Problem}\\
- \nabla \cdot (k(\T) \nabla \T) + v \cdot \nabla \T =& f   \qquad\text{in } \Omega. \nonumber
\end{align}
In the manner of the Boussinseq approximation, \cite{approximation,boussinesq_review}, possible variations of the density due to temperature differences are neglected except for their most significant effect in the form of buoyancy forces ($-α\T g$), so that in \eqref{Eqn: Model Problem} $ρ$ is constant.

\paragraph{Boundary conditions}
Furthermore, we have fluid flow boundary conditions 	
\begin{equation}\label{bc:fluid}
\begin{aligned}
v=& v_D  \qquad \text{on } \Gamma_{vD}, \\ 
\rho\nu(\T) \frac{\partial v}{\partial {n}} - p \cdot n =& \black{\sigma_N}
\qquad\text{on } \Gamma_{vN},
\end{aligned}
\end{equation}
where $v_D$ is some given Dirichlet data and $\black{\sigma_N}$ is the Neumann data.
It is assumed that $\Gamma_{vD}$ has a non-zero ($(d-1)$-dimensional) measure. 
Moreover, $n$ denotes the outer normal vector.
Next, we have the
temperature boundary conditions
\begin{align*}
\T = & \T_D    \qquad\text{on } \Gamma_{\T D}, \nonumber \\
k(\T)\frac{\partial \T}{\partial n}  = & 0
\qquad\text{on }\Gamma_{\T N}, \nonumber
\end{align*}
with given Dirichlet data $\T_D> 0$. 
Here, we assume that $\Gamma_{\T D}$ has a non-zero measure. 
The different boundary parts and their values are specified for each 
example.

\subsection{Weak form}
\label{sec_weak_form}
In this subsection, we present the weak formulation via variational-monolithic coupling by standard arguments. First we construct the function spaces:
\begin{align*}
V^v &:= \{ v\in [H^1(\Omega)]^d | \; v=0 \text{ on } \Gamma_{vD}\},\\
V^p &:= L^2(\Omega)/\mathbb{R},\\
V^{\T} &:= \{ \T\in H^1(\Omega) | \; \T=0 \text{ on } \Gamma_{\T D}\}.
\end{align*}
Let us introduce the space $X:= V^v \times V^p \times V^{\T}$. 
In the following, we give a weak formulation using 
the notation of a semi-linear form $A(U)(\Psi)$ which is nonlinear in its 
first argument (i.e., the trial function $U$) and linear 
with respect to the second argument (i.e., the test function $Ψ$). 
This gives
\begin{Problem}
	Let $\{\pi v_D, 0, \pi \T_D\}$ be an extension of nonhomogeneous Dirichlet 
	data. Furthermore, let the semi-linear form 
	$A(U)(\Psi)$ be given by
	\begin{align}\label{def:A}
	A(U)(\Psi) := &
	(\rho(v \cdot \nabla) v,\psi^v) + (\sigma, \nabla \psi^v) - \langle \black{\sigma_N}, \psi^v  \rangle_{\Gamma_{vN}}  - (\alpha \T g,\psi^v)\nn\\
	& + (\nabla \cdot v, \psi^p)\nn\\
	& + (k(\T) \nabla \T,∇ \psi^{\T}) 
	+ (v \cdot \nabla \T, \psi^{\T})
	-(f,\psi^{\T}),\qquad U\in X^D,Ψ\in X,
	\end{align}
	where $\sigma$ is as in \eqref{def:sigma}, $\nu$ is taken from \eqref{def:ny} and $k\colon ℝ\to (0,\infty)$ is a continuous positive function.\\
	Find $U=(v,p,\T)\in X^D:=\{\pi v_D, 0, \pi \T_D\} + X$ such that
	\begin{equation}
	\label{eq_1}
	A(U)(\Psi) = 0 \quad\forall \Psi:= (\psi^v,\psi^p,\psi^{\T})\in X.
	\end{equation}
	
\end{Problem}

We recall Theorem~2.1 of \cite{LoBo96}, merely adjusting the notation to the problem description given above:
\begin{theorem}
	Let $\Om\subset ℝ^d$, $d\in\{2,3\}$, be a bounded domain with Lipschitz 
	boundary, $ν,k\in C^0(ℝ)$ positive functions, 
	assume $\alpha\in\mathbb{R}$, $g\in [L^2(\Omega)]^d$, 
	$\T_D\in H^{\frac12}(\partial\Om)\cap L^\infty(\partial\Om)$. 
	Let $f=0$, $v_D=0$, 
	$Γ_{vD}=\partial\Om$, 
	$Γ_{vN}= \emptyset$, 
	$Γ_{\T D}=\partial\Om$.
	and $Γ_{\T N}=\emptyset$. 
	Then the problem has a weak solution.
\end{theorem}
The extension of the problem of \cite{LoBo96} by a non-zero external source $f$ is straight-forward. For the change of the fluid boundary conditions to the mixed conditions in \eqref{bc:fluid} we refer to \cite[Sec. 3]{guerra}, where solvability of a Navier-Stokes system with these boundary conditions, but without any influence of the temperature, was proven.

\subsection{Discretization and numerical solution}
\label{sec_disc}
The problem in equation \eqref{eq_1} is discretized with a Galerkin 
finite element scheme \cite{Cia87} \black{using quadrilaterals 
	with hanging nodes for local mesh refinement. The choice of employing 
	quadrilaterals is motivated by our finite element library deal.II \cite{dealII91,deal2020}.} 
To this end, we introduce
finite dimensional conforming subspaces $X_h\subset X$,
where $X_h=V_h^v \times V_h^p \times V_h^{\T}$.
Furthermore, let $\{\pi v_{h,D}, 0, \pi \T_{h,D}\}$
be an extension of the discretized boundary data.
Then, the problem statement reads:
Find $U_h=(v_h,p_h,\T_h)\in X^D_h =\{\pi v_{h,D}, 0, \pi \T_{h,D}\} + X_h$
such that 
\begin{equation}
\label{eq_2}
A(U_h)(\Psi_h) = 0 \quad\forall \Psi_h\in X_h.
\end{equation}
This finite-dimensional nonlinear system is solved with Newton's method:
Given an initial guess $U_h^0 \in \{\pi v_{h,D}, 0, \pi \T_{h,D}\} + X_h$,
find $\delta U_h\in X_h$ for $j=1,2,3,\ldots$ such that
\begin{align}
A'(U_h^j)(\delta U_h, \Psi_h) &= 
- A(U_h^j)(\Psi_h) \quad\forall \Psi_h\in X_h.\\
U_h^{j+1} &= U_h^{j} + \omega \delta U_h,
\end{align}
where $\omega\in (0,1]$ is a line-search parameter for globalization.
Inside Newton's method, the arising systems of linear equations 
are solved with a sparse direct solver (UMFPACK \cite{DaDu97}).

\section{Goal-oriented error control}
\label{sec_goal}
In goal-oriented error estimation the aim is to estimate the 
error in a certain quantity of interest $J:X^D \mapsto \mathbb{R}$. 
Examples for such quantities of interest could be a point evaluation, 
an integral evaluation of any solution component 
or some other possibly nonlinear quantity $J$.
In the following, first the abstract primal problem 
from before is stated, and subsequently the associated 
adjoint problem is given. Both are employed to derive
an error identity.

\subsection{Primal problem}
The primal problem is given by: Find $U=(v,p,\T)\in X^D$ such that 
\begin{equation}
\label{eq_primal_nochmal}
A(U)(\Psi) = 0 \quad\forall \Psi\in X.
\end{equation}
The discrete version of this problem reads as 
discussed above: Find $U_h=(v_h,p_h,\T_h)\in X^D_h$, such that 
\begin{equation*}
A(U_h)(\Psi_h) = 0 \quad\forall \Psi_h\in X_h.
\end{equation*}
Our aim is to obtain $J(U)$, however all we can compute is $J(U_h)$. To estimate the error we use the adjoint problem for $J$ as proposed in \cite{BeckerRannacher1995,BeRa01}. 
\black{This approach is known 
	as dual-weighted residual method (DWR), which is inspired by optimal 
	control, and therefore both are conceptionally similar. In the DWR
	method we aim to minimize the approximation error subject to 
	a PDE constraint, here $A(U)(\Psi) = 0$. The approximation error 
	may consist of the discretization error only, but can also include 
	iteration errors \cite{MeiRaVih109,RanVi2013,dolejvsi2021goal} or model errors \cite{OdPr02,BraackErn02}.
	This minimization problem is given by \black{\cite[Section 2.2]{BeRa01}}
	\[
	\min J(U) \text{ s.t. } A(U)(\Psi) = 0,
	\]
	which can be solved by formulating the Lagrangian $L(U,Z)$ with 
	the adjoint variable $Z\in X$. The resulting optimality system 
	is obtained by differentiation with respect to $U$ and $Z$,
	which is conceptionally similar to numerical optimization such 
	as optimal control or topology optimization \cite{Troe09,HiPiUlUl09,All15,Allaire2006}.
}

\subsection{Adjoint problem}
The adjoint problem is given by: Find $Z=(z^v,z^p,z^\T) \in X$ such that 
\begin{equation}
\label{Eq: adjoint}
A'(U)(\Psi,Z) = J'(U)(\Psi) \quad\forall \Psi\in X,
\end{equation}
where $A'$ and $J'$ are the Fr\'echet derivatives with respect to $U$. However, also the adjoint problem has to be discretized. The discretized adjoint problem 
reads:  Find $Z_h=(z_h^v,z_h^p,z_h^\T) \in X_h$ such that 
\begin{equation}
\label{Eq: adjointh}
A'(U_h)(\Psi_h,Z_h) = J'(U_h)(\Psi_h) \quad\forall \Psi_h\in X_h.
\end{equation}

\subsection{Error representation}
Using the solutions of the primal and adjoint problem, we obtain the following theorem:

\begin{theorem}\label{Theorem: Error Representation} 
	Let $A$ be as in \eqref{def:A} and $J \in \mathcal{C}^3(X^D,\mathbb{R})$. 
	If $U$ solves (\ref{eq_primal_nochmal})
	and $Z$ solves (\ref{Eq: adjoint}) for $U$, 
	then for every $\tilde{U}\in X^D$ and $\tilde{Z} \in X$, the error $J(U)-J(\tilde{U})$ can be written as 
	\begin{align} \label{Error Representation}
	\begin{split}
	J(U)-J(\tilde{U})&= \frac{1}{2}\rho(\tilde{U})(Z-\tilde{Z})+\frac{1}{2}\rho^*(\tilde{U},\tilde{Z})(U-\tilde{U}) 
	-\rho (\tilde{U})(\tilde{Z}) + \mathcal{R}^{(3)}(\tilde{U},\tilde{Z},e_u,e_z),
	\end{split}
	\end{align}
	where the primal and adjoint residuals are given by
	\begin{align*}
	\rho(\tilde{U})(\cdot) &:= -A(\tilde{U})(\cdot),\\
	\rho^*(\tilde{U},\tilde{Z})(\cdot) &:= J'(\tilde{U})-A'(\tilde{U})(\cdot,\tilde{Z}),
	\end{align*}
	respectively, 
	and the remainder term
	\begin{equation}
	\begin{split}	\label{Error Estimator: Remainderterm}
	\mathcal{R}^{(3)}(\tilde{U},\tilde{Z},e_u,e_z):=&\frac{1}{2}\int_{0}^{1}\big[J'''(\tilde{U}+se_u)(e_u,e_u,e_u) \\
	-&A'''(\tilde{U}+se_u)(e_u,e_u,e_u,\tilde{Z}+se_z)
	-3A''(\tilde{U}+se_u)(e_u,e_u,e_u)\big]s(s-1)\,ds,
	\end{split} 
	\end{equation}
	with $e_u=U-\tilde{U}$ and $e_z =Z-\tilde{Z}$. 	
\end{theorem}
This error representation allows us to represent the error in a different way. However  (\ref{Error Representation}) still depends on $U$ and $Z$, which are both unknown.
\begin{proof}
	For information on the proof, we refer to 
	\cite{BeRa01,RanVi2013,EndtLaWi20}. Note that
	we use positivity of  
	$\T$ in order 
	to avoid singularities in $\rho, \rho^*$ and $\mathcal{R}^{(3)}$.
\end{proof}

\begin{Remark}
	Since this error respresentation holds for all 
	$\tilde{U}$ and $\tilde{Z}$, it also holds 
	for $\tilde{U}=U_h$ and $\tilde{Z}=Z_h$, 
	provided that $U_h\in X^D$ and $Z_h\in X$. 
	We note that $X_h\subset X$, but for non-trivial 
	boundary data, $X_h^D\not\subset X^D$.
\end{Remark}

\subsection{Error estimators} 
If we replace $U$ and $Z$ in  (\ref{Error Representation}) by approximations, we obtain an error estimator instead of an error representation.
This can be realized by higher order interpolation or enriched approximation. Both methods are described in more details in \cite{BeRa01} and a mixed method is presented in \cite{EndtLaWi20_smart,BeuEndtWi21}. In this work, we will use and describe enriched approximation in more detail. We consider $X_h^{(2)}$ and $X_h^{0,(2)}$ to be enriched spaces, i.e $X_h \subset X_h^{(2)} \subset X$ and $X_h^0 \subset X_h^{0,(2)} \subset X^0$. Examples of such enriched spaces can be generated by refining the mesh or using other finite elements.
This leads us to the enriched model problem: Find $U_h^{(2)}=(v_h^{(2)},p_h^{(2)},\T_h^{(2)})\in X_h^{(2)}$, 
such that 
\begin{equation} \label{Eq: primal enriched}
A(U_h^{(2)})(\Psi_h^{(2)}) = 0 \quad\forall \Psi_h^{(2)}\in X_h^{0,(2)}.
\end{equation}
The enriched adjoint problem reads:  Find $Z_h^{(2)}=(z_h^{v,(2)},z_h^{p,(2)},z_h^{\T,(2)}) \in X_h^{0,(2)}$ such that 
\begin{equation} \label{Eq: adjoint enriched}
A'(U_h^{(2)})(\Psi_h^{(2)},Z_h^{(2)}) = J'(U_h^{(2)})(\Psi_h^{(2)}) \quad\forall \Psi_h^{(2)}\in X_h^{0,(2)}.
\end{equation}
As above, $U$ and $Z$ in the right hand side of  the error representation (\ref{Error Representation}) are replaced  by $U_h^{(2)}$ and $Z_h^{(2)}$ and $\tilde{U}$ and $\tilde{Z}$ by $U_h$ and $Z_h$, respectively.  We obtain the  error estimation formula
\begin{equation} \label{Error estimator full}
J(U)-J(U_h)\approx \underbrace{\frac{1}{2}\rho(U_h)(Z_h^{(2)}-Z_h)+\frac{1}{2}\rho^*(U_h,Z_h)(U_h^{(2)}-U_h)}_{\eta_h} 
\underbrace{-\rho (U_h)(Z_h)}_{\eta_k} + \underbrace{\mathcal{R}^{(3)}(U_h,Z_h,e_u^{(2)},e_z^{(2)})}_{\eta_\mathcal{R}},
\end{equation}
where $e_u^{(2)}:=U_h^{(2)}-U_h$ and $Z_h^{(2)}-Z_h$.
Here $U_h^{(2)}$ and $Z_h^{(2)}$ are the solutions of (\ref{Eq: primal enriched}) and (\ref{Eq: adjoint enriched}) respectively. 
\begin{Remark}
	For the error estimator (\ref{Error estimator full}), $U_h$ and $Z_h$ 
	need not be 
	the exact solutions of the discrete problem but can be some 
	approximations of it as well.
\end{Remark}

\paragraph*{The first part of the error estimator $\eta_h$}$~$\newline
The part $\eta_h$ represents the discretization error as proposed in 
\cite{RanVi2013,EndtLaWi20,Endt21}. \new{Here $$\eta_{h,p}:=\rho(U_h)(Z_h^{(2)}-Z_h)$$ is the primal part 
	of the error estimator and $$\eta_{h,a}:=\rho^*(U_h,Z_h)(U_h^{(2)}-U_h)$$ 
	is the adjoint part.} In the literature the adjoint part of the error estimator 
is often replaced with the primal error part $\rho$.
In \cite{BeRa01,bruchhauser2020dual}, it is proven that the adjoint 
part can be expressed as the primal part and higher order terms depending 
on the problem and the goal functionals. For moderate nonlinear 
problems, this approximation often works very well \cite{BeRa01,BrRi06}.
In this work, both parts are considered.
The localization is done by using the partition of unity technique 
proposed in \cite{RiWi15_dwr}. Alternatives 
are the filtering approach \cite{BraackErn02} (which works 
as well on the variational level) or integration 
by parts \cite{BeRa01}. 
However, specifically for coupled problems 
(as in the current work) the latter is error prone
and computationally expensive
since the strong form operators must be evaluated.

\paragraph*{The second part of the error estimator $\eta_k$}$~$\newline
This part mimics the iteration error as in \cite{RanVi2013,EndtLaWi20}. 
It can be used as stopping criterion of
the solver on each level. If $U_h$ is the exact 
solution of (\ref{eq_2}) then this part vanishes.

\paragraph*{The third part of the error estimator $\eta_\mathcal{R}$}$~$\newline
This part is usually of higher order \cite{BeRa01,RanVi2013}. Mostly, 
$\eta_\mathcal{R}$ is neglected in the evaluation of the error estimator.
Recent studies and investigations of $\eta_\mathcal{R}$
with the help of numerical examples were undertaken in \cite{EndtLaWi20}.

\paragraph*{The practical error estimator}$~$\newline
After the previous assumptions and explanations, the 
practical error estimator is given by $\eta:=\eta_h+\eta_k$.
As proven in \cite{EndtLaWi20} (see also \cite{Endt21}), 
this error estimator is 
efficient and reliable if a certain saturation 
assumption is fulfilled. Furthermore, interpolation techniques
for a new class of algorithms were established in \cite{EndtLaWi20_smart}.

\subsection{Finite element discretization and polynomial spaces}
Having the primal and adjoint problems at hand, we employ the following 
finite elements in our algorithms and numerical experiments.
We use 
\begin{itemize}
	\item continous piecewise bi-quadratic functions $Q_2^c$ for the velocity $v$ and
	\item continuous piecewise bilinear functions $Q_1^c$ for the pressure $p$ and temperature $\T$.
\end{itemize}
\black{The adatively refined mesh will lead to hanging nodes in the mesh \cite{CaOd84}. These nodes are constrained such that generated functions in the finite element space are continous. For more information about this topic we refer to \cite{rheinboldt1980data,bangerth2009data}.}
For the enriched space we have
\begin{itemize}
	\item continous piecewise bi-quartic functions $Q_4^c$ for the velocity $v$ and
	\item continuous piecewise quadratic functions $Q_2^c$ for the pressure $p$ and temperature $\T$.
\end{itemize}
A comparison of different finite elements and uniform mesh refinement for the enriched space can be found in \cite{EnLaThieWienumath}.

\subsection{Multiple goal functionals and algorithms}
In many applications, such as multiphysics problems, or coupled problems in general (as in this work), more than one goal functional is of interest. 
Let us assume we are interested in $N$ goal functionals $J_1, \ldots, J_N$.
Later in Section \ref{sec_tests}, we have up to $N=7$.
A straightforward application of the previous concepts would be to 
compute an error estimator for each functional and to combine them afterwards. 
However, for this approach we have to solve the adjoint problem $N$ times;
see \cite{HaHou03,Ha08}. This would lead to a non-acceptable computational cost.
To overcome this problem, 
several techniques have been proposed such as
a combined functional by solving an additional dual-dual problem for 
the sign computation in the combined functional \cite{HaHou03,Ha08}, 
using generalized Green's functions \cite{EsHoLaMulti2005},
by a linear combination of the functionals \cite{PARDO20101953,AlvParBar2013},
by reformulation to a minimization problem where the quantites of interest serve as constraints \cite{KerPruChaLaf2017,kergrene2018goal},
and combined functionals with hierarchical higher order approximations 
for the sign computation \cite{EndtWi17,EnLaWi18,Endt21}.
In this work we follow the approach proposed in 
\cite{EndtWi17,EnLaWi18}. 
We combine the functionals to one by using
\begin{equation}
J_c:= \sum_{i=1}^{N} \omega_i J_i,
\end{equation}
with 
$$
\omega_i:= 
\begin{cases}\text{sign}(J_i(U_h^{(2)})-J_i(U_h))/|J_i(U_h)| 
\quad  for |J_i(U_h)|\geq 10^{-12},\\
10^3 \text{sign}(J_i(U_h^{(2)})-J_i(U_h)).
\end{cases} 
$$ 
The overall algorithmic realization using $J_c$ has been described 
in detail in \cite{EndtLaWi20} and we simply apply exactly these 
schemes to the Boussinesq model in the current paper.

\section{Numerical experiments}
\label{sec_tests}
In this section, we investigate three numerical examples. The programming 
code is based on the open-source finite element library 
deal.II \cite{dealII91,deal2020}. All geometry data 
and material parameters are given in SI units.
To measure the quality of our error estimators we employ 
so-called effectivity indices.
\black{
	\begin{Definition}
		For the functional $J$ the efficitivty index $I_{eff}$, the primal effectivity index $I_{eff,p}$ and the adjoint effectivity index $I_{eff,a}$  are  defined as $$I_{eff}:=\frac{\eta_h}{J(U)-J(U_h)}, \qquad I_{eff,p}:=\frac{\eta_{h,p}}{J(U)-J(U_h)}\quad \text{and } \quad I_{eff,a}:=\frac{\eta_{h,a}}{J(U)-J(U_h)}.$$
	\end{Definition}
	For the real error $J(U)-J(U_h)$, we compute a
	reference solution on sufficiently refined meshes as it is often done.
}

\subsection{A flow benchmark}
In this first example we apply our method to a problem featuring a flow around a cylinder as in \cite{TurSchabenchmark1996}. 


\subsubsection{Configuration, geometry, parameters, boundary conditions}
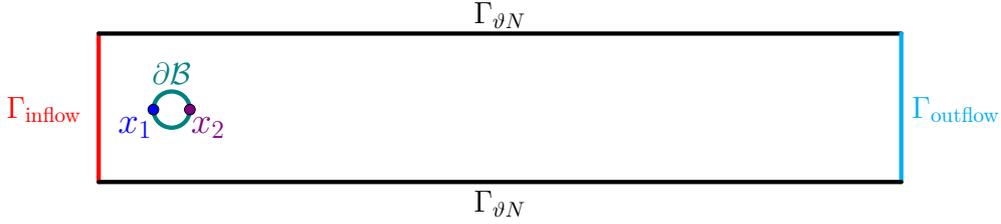
\begin{figure}[H]
	\scalebox{0.8}{
		\begin{tikzpicture}[line cap=round,line join=round,>=triangle 45,x=6cm,y=6cm]
		\clip(-0.4,-0.2) rectangle (2.6,0.6);
		\draw [line width=2pt,color=red] (0,0)-- (0,0.41);
		\draw [line width=2pt,color=black] (0,0.41)-- (2.2,0.41);
		\draw [line width=2pt,color=cyan] (2.2,0.41)-- (2.2,0);
		\draw [line width=2pt,color=black] (2.2,0)-- (0,0);
		\draw [line width=2pt,color=teal] (0.2,0.2) circle (0.05);
		\begin{scriptsize}
		\draw[color=red] (-0.15,0.2) node {\Large $\Gamma_{\text{inflow}}$};
		\draw[color=black] (1.1,0.46) node {\Large$\Gamma_{\T N}$};
		\draw[color=cyan] (2.35,0.2) node {\Large$\Gamma_{\text{outflow}}$};
		\draw[color=black] (1.1,-0.06) node {\Large$\Gamma_{\T N}$};
		\draw[color=teal] (0.2,0.3) node {\Large$\partial \mathcal{B}$};
		\draw[color=blue] (0.1,0.15) node {{\LARGE$x_1$}};
		\draw [fill=blue] (0.15,0.2) circle (2.5pt);
		\draw[color=violet] (0.3,0.15) node {{\LARGE$x_2$}};
		\draw [fill=violet] (0.25,0.2) circle (2.5pt);
		\end{scriptsize}
		\end{tikzpicture}
	}	
	\caption{The domain $\Omega$ with boundary conditions. \label{Fig:OmegaBoundary}}
\end{figure}

The domain is given by $\Omega:=(0,2.2)\times (0,0.41)\setminus \mathcal{B}$ where $\mathcal{B}:=\{x \in \black{\mathbb{R}^2}:|x-(0.2,0.2)|<0.05\}$. The domain as well as the boundary conditions are depicted in Figure~\ref{Fig:OmegaBoundary}.

In the strong model given in 
Section \ref{sec_model_strong_form}, we set $g=(0,0)^T, k(\T)=1$ and\footnote{The values are computed by using the data for 293.15K and 353.15K in \url{https://www.lss.ovgu.de/lss_media/Downloads/Lehre/Strömungsmechanik/Arbeitsheft/IV.pdf}.} 
$E_A = 1.49\times 10^{4}$, 
$\nu_0=2.22\times 10^{-6}$, 
and $R=8.31$. 
Next, the thermal expansion coefficient is $\alpha=6.88\times 10^{-5}$ 
and the density is $\rho=998.21$. \black{Constant conductivity was chosen for 
	simplicity, in order to focus on prototype situations within the development 
	of the multigoal framework for the Boussinesq model.}

Furthermore we have no-slip boundary conditions on \black{$\Gamma_{\T N}$} and $\partial \mathcal{B}$, do-nothing  conditions on $\Gamma_{\text{outflow}}$ and an inflow on $\Gamma_{\text{inflow}}$, i.e
\begin{align*}
v=&0 \qquad \qquad\text{ on } \black{\Gamma_{\T N}} \cup \partial \mathcal{B},\\
v=& v_{in} \qquad \quad\text{ on } \Gamma_{\text{inflow}},\\
\sigma\cdot n= &0 \qquad\qquad \text{ on } \Gamma_{\text{outflow}},\\
\frac{\partial\T}{\partial n} = &0 \qquad\qquad \text{ on }\partial \Omega \setminus (\partial \mathcal{B} \cup \Gamma_{\text{inflow}}),\\
\T= & \T_{\text{inflow}} \qquad \text{ on }\Gamma_{\text{inflow}},\\
\T= &\T_{\partial \mathcal{B}} \qquad\quad \text{ on }\partial \mathcal{B},\\
\end{align*}
where $v_{in}(x,y):= 4 v_m \frac{y(H-y)}{H^2}$ with $v_m=0.3$ and $H=0.41$.
In the following, 
two possible configurations are considered:
\begin{itemize}
	\item "cold to warm": $\T_{\text{inflow}}=278.15$, $\T_{\partial \mathcal{B}}=353.15$,
	\item "warm to cold": $\T_{\text{inflow}}=353.15$, $\T_{\partial \mathcal{B}}=278.15$.
\end{itemize}

\subsubsection{Goal functionals}
In this example, a pressure difference serves as goal functional:
$$p_{\text{diff}}(U):= p(x_1)-p(x_2),$$
where $x_1:=(0.15,0.2)$ and $x_2:=(0.25,0.2)$ \new{as in the original benchmark problem \cite{TurSchabenchmark1996}}.
The reference values are 
\begin{align*}
p_{\text{diff}}(U)=& 114.68898895581040\quad \text{for "cold to warm"} \qquad \text{ and} \\
p_{\text{diff}}(U)=& 101.97737719601436\quad \text{for "warm to cold".}
\end{align*}
We remark that due to the pointwise evaluation $p_{\text{diff}}(U)$ is not well-defined on the solution space $X^D$ introduced in Section \ref{sec_weak_form} \black{and the adjoint equation features Dirac delta distributions at $x_1$ and $x_2$ on the right hand side.} Therefore higher regularity of 
solutions to \eqref{eq_1} (respectively \eqref{Eqn: Model Problem}) 
has to be assumed here. 
For corresponding conditions on data and domain cf. also 
\cite[Theorem 2.3]{LoBo96}.

\subsubsection{Discussion of our findings}
\begin{figure}[H]
	\includegraphics[scale=0.28]{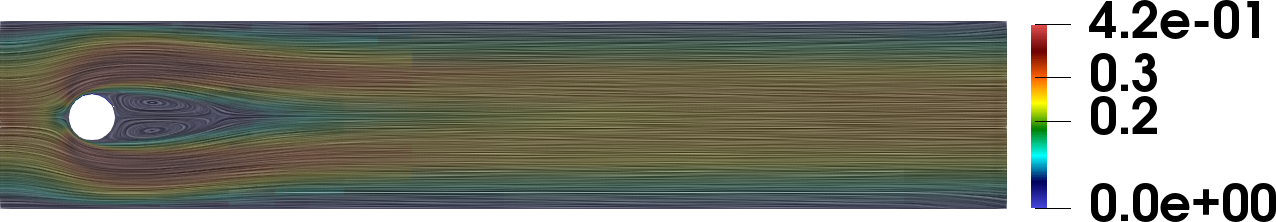}\\[0.2cm]
	\includegraphics[scale=0.28]{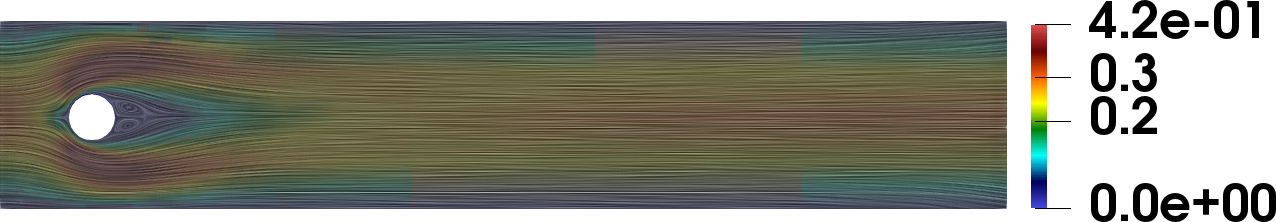}
	\caption{The magnitude of the velocity and the streamlines for "cold to warm" (above) and "warm to cold" (below)\label{Figure: Flow around cylinder V}}
\end{figure}

\begin{figure}[H]
	\includegraphics[scale=0.25]{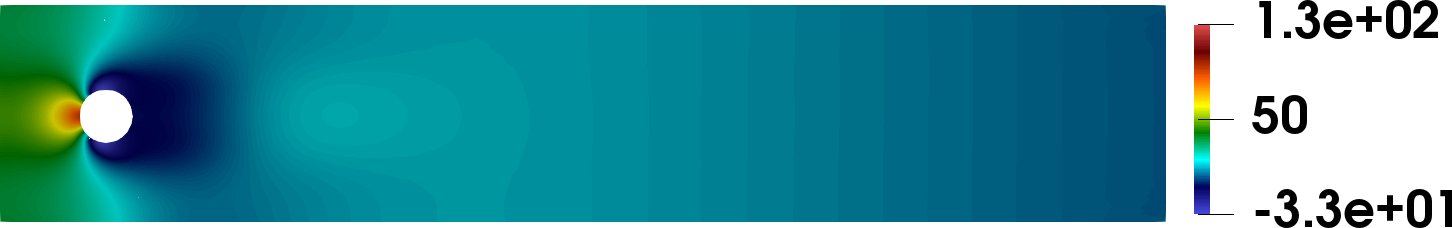}\\[0.2cm]
	\includegraphics[scale=0.25]{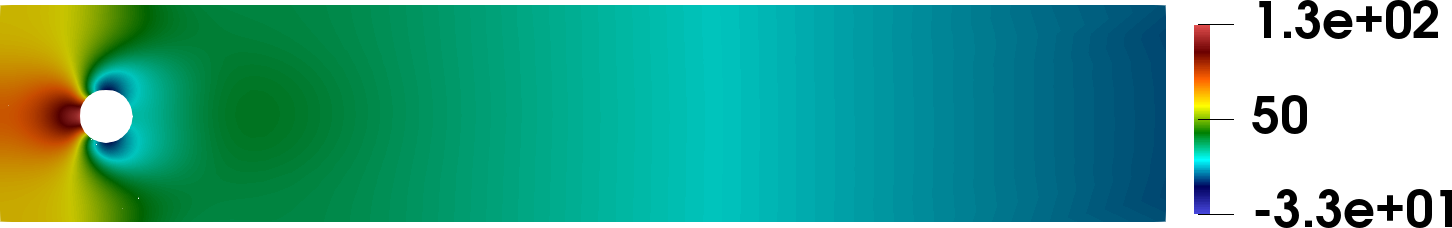}
	\caption{The pressure $p$ for the configuration "cold to warm" (above) and "warm to cold" (below)\label{Figure: Flow around cylinder P}}
\end{figure}

\begin{figure}[H]
	\includegraphics[scale=0.25]{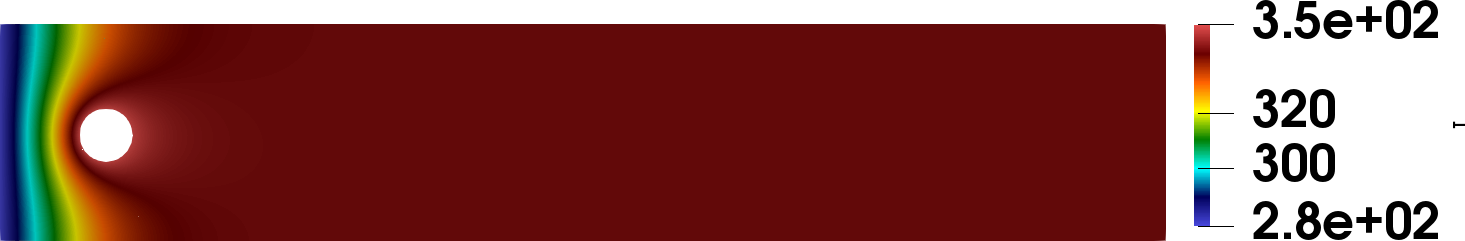}\\[0.2cm]
	\includegraphics[scale=0.25]{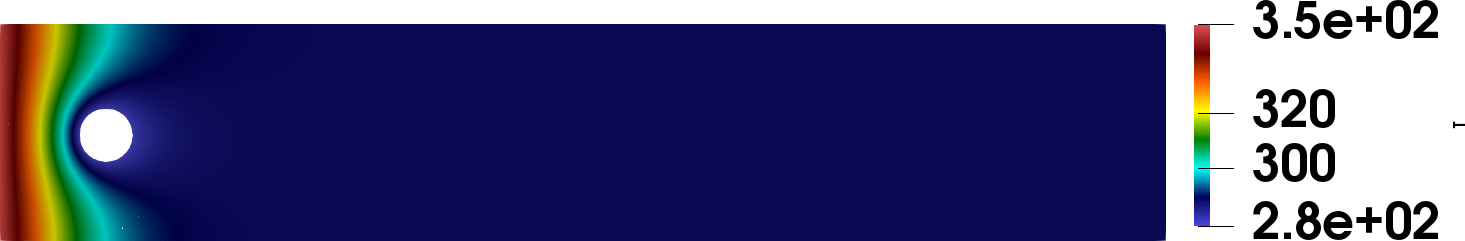}
	\caption{The temperature $\T$ for the configuration "cold to warm" (above) and "warm to cold" (below)\label{Figure: Flow around cylinder T}}
\end{figure}

\begin{figure}[H]
	\includegraphics[scale=0.28]{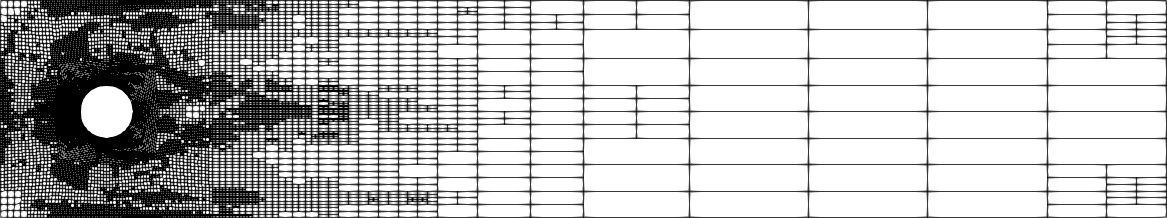}\\[0.2cm]
	\includegraphics[scale=0.28]{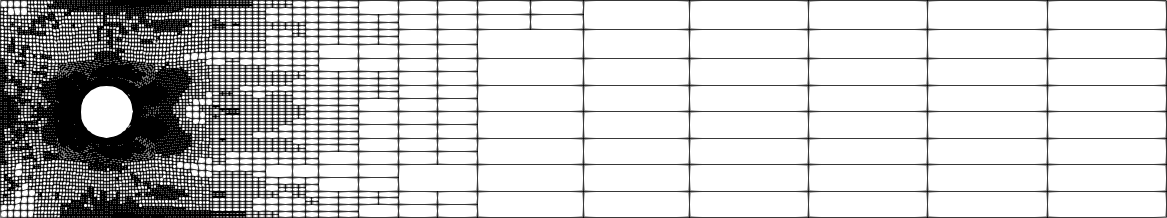}
	\caption{The mesh for the configuration "cold to warm" (above) and "warm to cold" (below)\label{Figure: Flow around cylinder Mesh}}
\end{figure}

The numerical results are displayed in Figures~\ref{Figure: P_diff Error Cold to warm} to \ref{Figure: Flow around cylinder Mesh}. More precisely, 
Figure \ref{Figure: P_diff Error Cold to warm} and Figure~\ref{Figure: P_diff Error Warm_to_cold} show the error and error estimator for the configurations "cold to warm" and "warm to cold", respectively.
The magnitude of the velocity (including streamlines), the pressure and the temperature are visualized in the next three figures for both configurations.
Finally, the meshes are displayed in  Figure~\ref{Figure: Flow around cylinder Mesh}.

We observe that in the configuration "cold to warm" the vortices are much bigger than in the configuration "warm to cold". This is a result of the temperature dependent viscosity, which is smaller around the cold cylinder. 
This also leads to higher convection terms in the configuration "cold to warm". 
In Figure~\ref{Figure: Flow around cylinder Mesh}, we notice that this effect also has a big influence in the adaptive refinement. 
In both configurations the temperature is almost constant 
on the other side of the cylinder, see  Figure~\ref{Figure: Flow around cylinder T}.
Furthermore, it is close to the temperature at $\partial \mathcal{B}$.
From Figure~\ref{Figure: Flow around cylinder P}, we deduce that the high dependency of the viscosity on the temperature has a big impact on the pressure $p$. 
The error reduces approximately with the rate $\mathcal{O}(\text{DOFs}^{-1})$ for both configurations, cf.  Figure~\ref{Figure: P_diff Error Cold to warm} and Figure~\ref{Figure: P_diff Error Warm_to_cold}, respectively.
Surprisingly, the error estimator shows a more uniform behavior in the convergence than the error itself.
One can observe a strong refinement of the mesh around the cylinder. 
On the other hand, the mesh has almost no refinement on the right side,
namely on the outflow boundary.

\begin{figure}[H]
	\begin{minipage}{0.49\linewidth} 
		\ifMAKEPICS
		\begin{gnuplot}[terminal=epslatex,terminaloptions=color]
			set output "Figures/Examplecoldtowarm.tex"
			set title "cold to warm"
			set key top left
			set key opaque
			set datafile separator "|"
			set logscale x
			set logscale y
			set grid ytics lc rgb "#bbbbbb" lw 1 lt 0
			set grid xtics lc rgb "#bbbbbb" lw 1 lt 0
			set xlabel '\text{DOFs}'
			set format '
			plot \
			'< sqlite3 Data/Turek_benchmark/p_Diff/Cold_to_Warm.db "SELECT DISTINCT DOFS_primal, abs(Exact_Error) from data "' u 1:2 w  lp lt 1 lw 3 title ' \footnotesize Error in $p_{\text{diff}}$', \
			'< sqlite3 Data/Turek_benchmark/p_Diff/Cold_to_Warm.db "SELECT DISTINCT DOFS_primal, abs(Estimated_Error) from data "' u 1:2 w  lp  lw 3 title ' \footnotesize Error Estimator of $p_{\text{diff}}$', \
			1000/x  title '$\mathcal{O}(\text{DOFs}^-1)$' lw  10										
			#					 '< sqlite3 Data/Multigoalp4/Higher_Order/dataHigherOrderJE.db "SELECT DISTINCT DOFs, abs(Exact_Error) from data "' u 1:2 w  lp lw 3 title ' \footnotesize Error in $J_\mathfrak{E}$', \
		\end{gnuplot}
		\fi
		{		\scalebox{0.55}{\input{Figures/Examplecoldtowarm.tex}}}
		\caption{The error and error estimator for $p_{\text{diff}}$ for the configuration "cold to warm".\label{Figure: P_diff Error Cold to warm}}
	\end{minipage}%
	\hfill
	\begin{minipage}{0.49\linewidth}
		\ifMAKEPICS
		\begin{gnuplot}[terminal=epslatex,terminaloptions=color]
			set output "Figures/Examplewarmtocold.tex"
			set title "warm to cold"
			set key top left
			set key opaque
			set datafile separator "|"
			set logscale x
			set logscale y
			set grid ytics lc rgb "#bbbbbb" lw 1 lt 0
			set grid xtics lc rgb "#bbbbbb" lw 1 lt 0
			set xlabel '\text{DOFs}'
			set format '
			plot \
			'< sqlite3 Data/Turek_benchmark/p_Diff/Warm_to_Cold.db "SELECT DISTINCT DOFS_primal, abs(Exact_Error) from data "' u 1:2 w  lp lt 1 lw 3 title ' \footnotesize Error in $p_{\text{diff}}$', \
			'< sqlite3 Data/Turek_benchmark/p_Diff/Warm_to_Cold.db "SELECT DISTINCT DOFS_primal, abs(Estimated_Error) from data "' u 1:2 w  lp  lw 3 title ' \footnotesize Error Estimator of $p_{\text{diff}}$', \
			1000/x  title '$\mathcal{O}(\text{DOFs}^-1)$' lw  10										
			#					 '< sqlite3 Data/Multigoalp4/Higher_Order/dataHigherOrderJE.db "SELECT DISTINCT DOFs, abs(Exact_Error) from data "' u 1:2 w  lp lw 3 title ' \footnotesize Error in $J_\mathfrak{E}$', \
		\end{gnuplot}
		\fi
		{		\scalebox{0.55}{\input{Figures/Examplewarmtocold.tex}}}
		\caption{The error and error estimator for $p_{\text{diff}}$ for the configuration "warm to cold".\label{Figure: P_diff Error Warm_to_cold}}
	\end{minipage}
\end{figure}

\subsection{Laser point source}
\label{sec_test_laser_point_source}
In this second example, we consider the flow in a square without 
inflow and outflow. The temperature on the 
boundary is constant. A heat point source\footnote{The motivation of this example 
	is due to laser wave guide modeling in which a laser causes temperature changes and 
	for which material starts to flow.} enters as right hand side 
into the flow equations and generates a flow field.

\subsubsection{Configuration, geometry, parameters, boundary conditions}	
We consider $\Omega := (0,1)^2$. The right hand side $f$ of \eqref{Eqn: Model Problem} is
\begin{equation*}
f(x) := 10 \psi_{x_0}(x) \qquad \text{where } \quad \psi_{x_0}(x):=10^4\sqrt{2\pi}e^{-{10^4|x-x_0|^2}},
\end{equation*}
with $x_0:=(0.75,0.75)$. This models a laser pointing at $x_0$. The boundary conditions are
\begin{equation*}
v=0 \quad \text{on } \partial \Omega \qquad \text{and} \qquad \T=293.15 \quad \text{on }\partial \Omega. 
\end{equation*}
The gravity is given by $g=(0,-9.81)^T$ and the thermal 
expansion coefficient is $\alpha=6.88\times 10^{-5}$.

\subsubsection{Goal functionals}
Our quantities of interest are the mean value of the velocity and the mean value of the temperature
\begin{equation*}
\bar{|v|}(U):=\frac{1}{|\Omega|}\int_{\Omega}|v| \dx \qquad \text{and} \qquad  \bar{\T}(U):=\frac{1}{|\Omega|}\int_{\Omega}\T \dx.
\end{equation*}
We use the adaptive strategy for $\bar{|v|}$, $\bar{\T}$, 
and for both of these at once in the combined 
functional $J_\mathfrak{E}$.

\subsubsection{Discussion of our findings}		

\begin{figure}[H]
	\centering
	\includegraphics[scale=0.210]{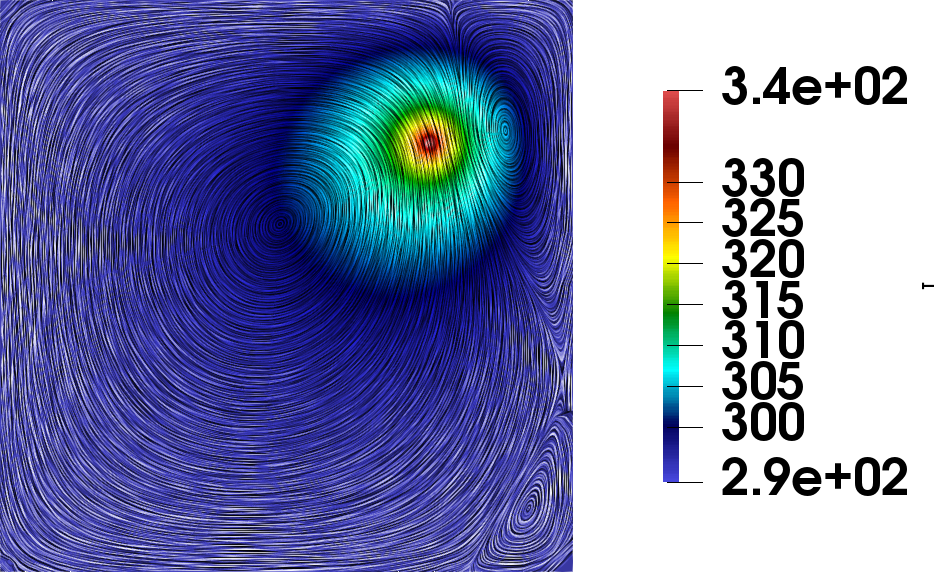}\hfil
	\includegraphics[scale=0.270]{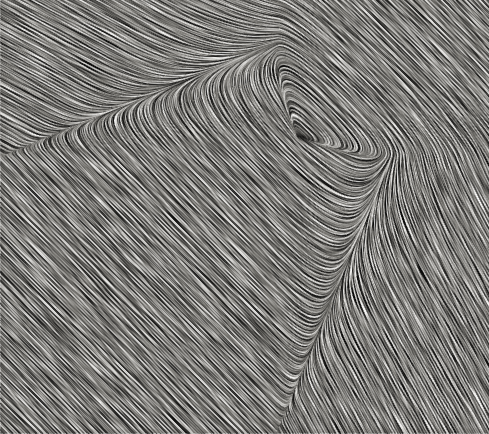}\hfil
	\caption{
		Section \ref{sec_test_laser_point_source}. 
		Velocity streamlines and temperature (left) and zoom into the small vortex at left bottom corner (right).
	}
	\label{fig_pick_4_out_of_6_A}
\end{figure}

\begin{figure}[H]
	
	\includegraphics[scale=0.12]{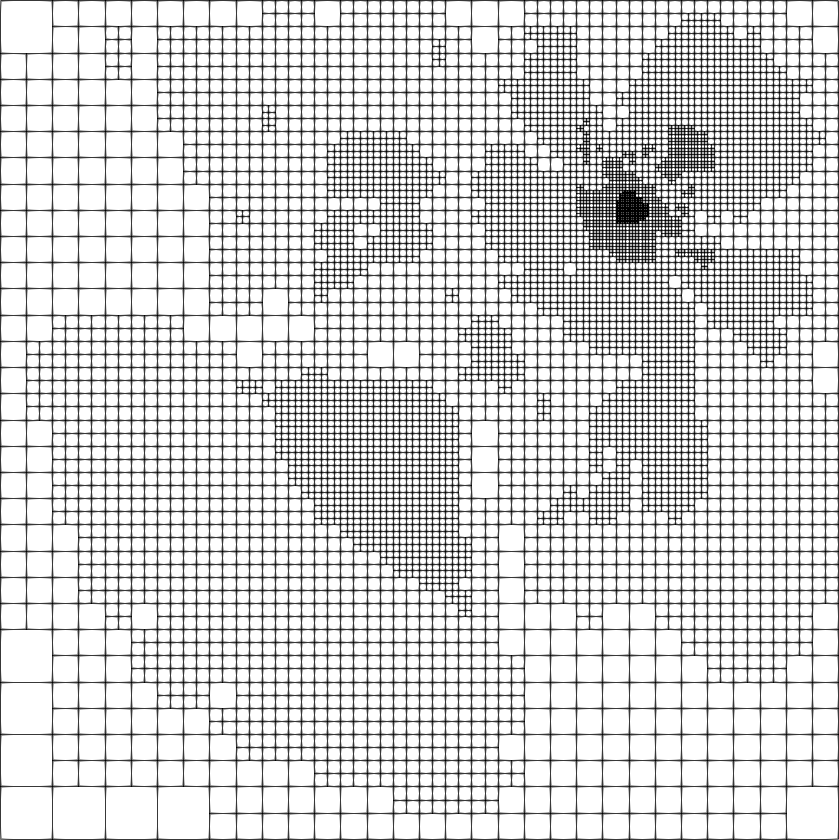} \hfil
	\includegraphics[scale=0.12]{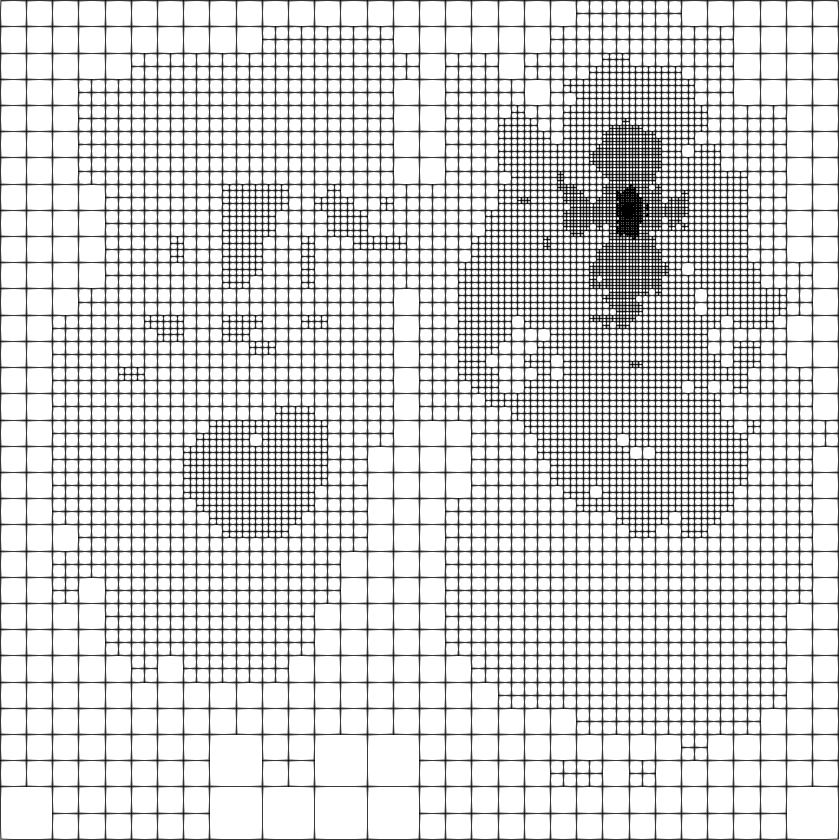}\hfil
	\includegraphics[scale=0.12]{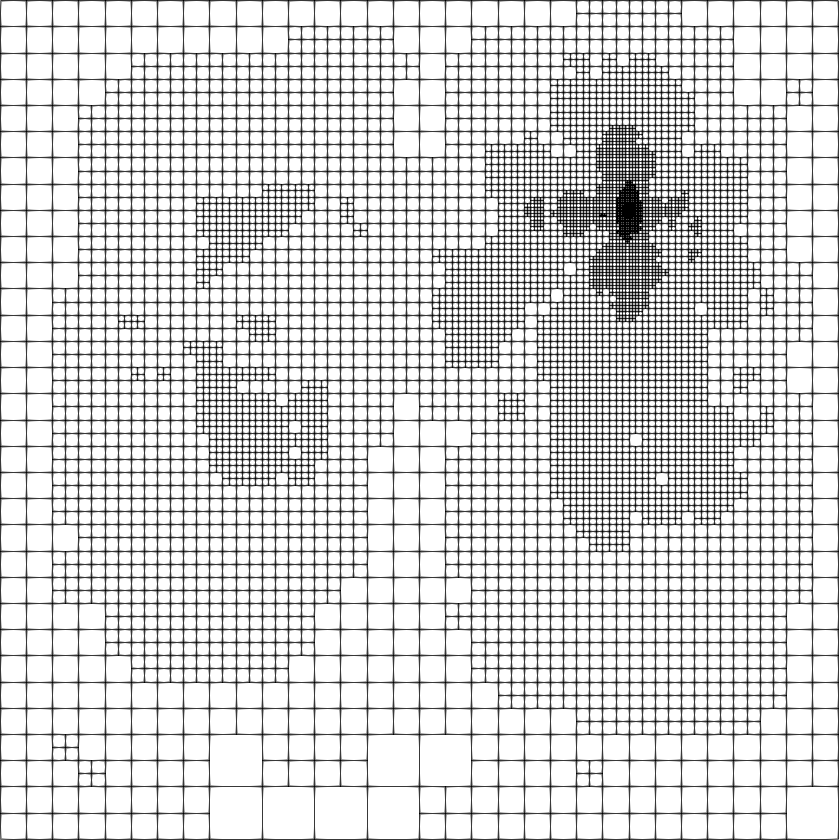}
	\caption{Section \ref{sec_test_laser_point_source}. The mesh after 21 refinements for $\bar{\T}$ (left), $J_\mathfrak{E}$ (center) and $\bar{|v|}$(right). \label{Figure: Laser on Square Combined Mesh}}
\end{figure}

\begin{figure}[H]
	
	\centering
	\includegraphics[scale=0.12]{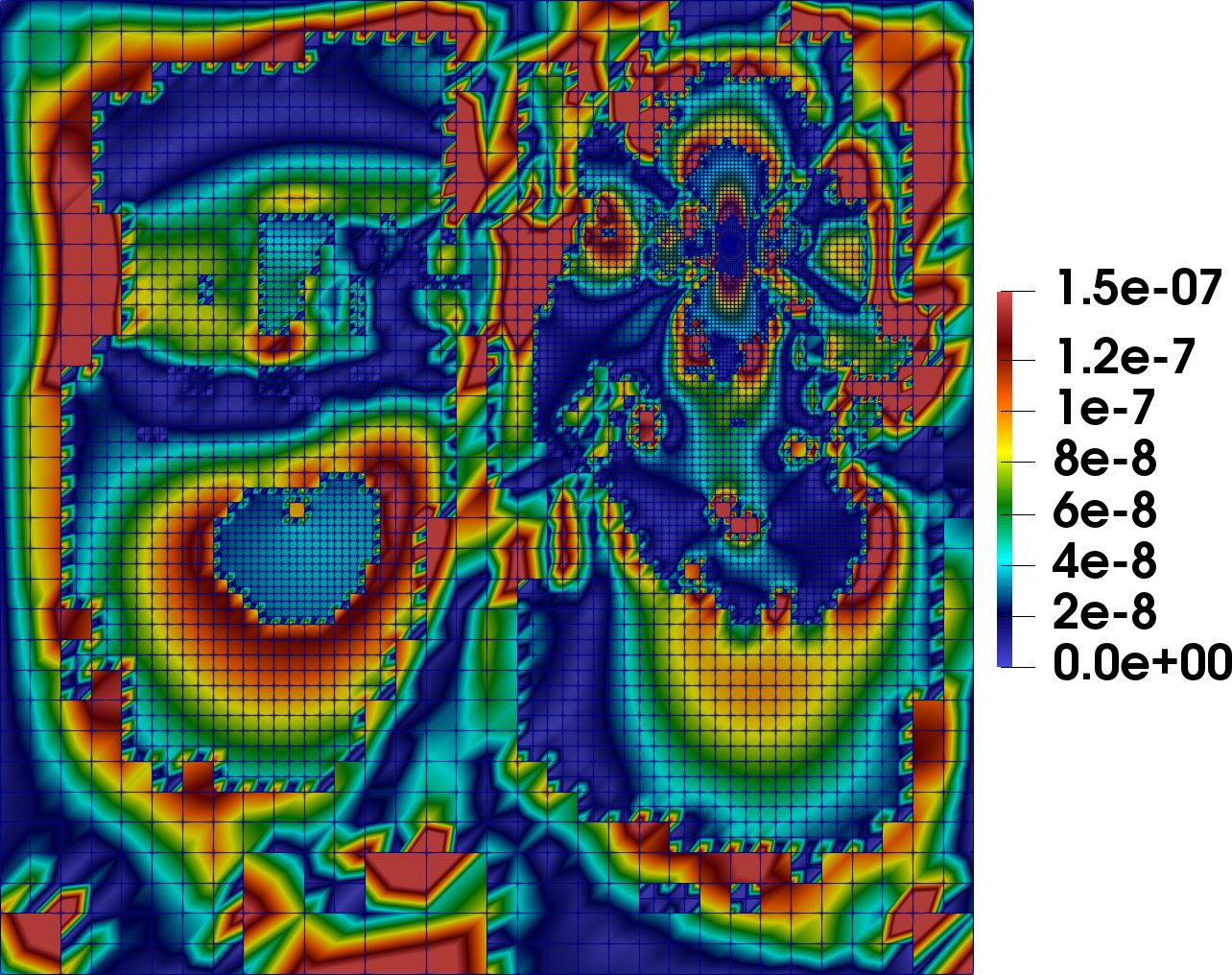}
	\hfil
	\includegraphics[scale=0.12]{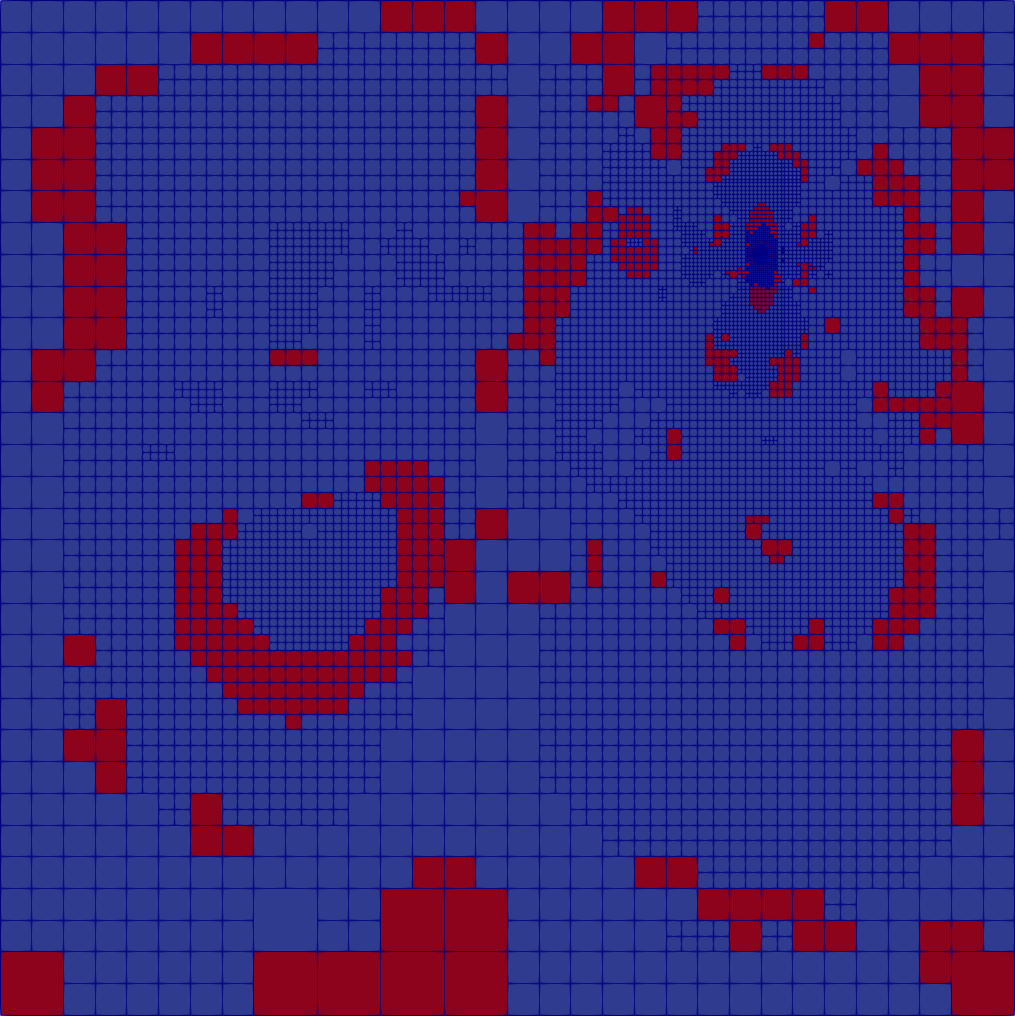}
	\caption{
		Section \ref{sec_test_laser_point_source}. 
		Error estimator (left) and marked elements for refinement (right) for the mesh $J_\mathfrak{E}$ displayed at the center of 
		Figure \ref{Figure: Laser on Square Combined Mesh}.
		\label{fig_pick_4_out_of_6_B}
	}
\end{figure}

The numerical results are presented in Figures \ref{Figure: Laser on Square Meantemperature} 
to \ref{fig_pick_4_out_of_6_B}. More precisely,
the first six figures display error and error estimator in one picture and \black{effectivity} index in a second picture for the temperature $\bar{\T}$, the absolute value of the velocity $\bar{|v|}$  and the 
combination of both, respectively.
The streamlines of the flow and the magnitude of the temperature 
$\T$ are displayed in Figure~\ref{fig_pick_4_out_of_6_A} left.
The refined meshes after 21 refinement steps are 
depicted in Figure \ref{Figure: Laser on Square Combined Mesh}.
In Figure \ref{fig_pick_4_out_of_6_B}, the value 
of the error estimator and the marked elements are displayed.
Good effectivity indices for $\bar{|v|}$, $\bar{\T}$ and for the combination $J_\mathfrak{E}$ are observed on the refined meshes.
Figure \ref{Figure: Laser on Square Combined} provides us information on the error of the individual functionals, the combined error and the error estimator. 
The error for $\bar{|v|}$ dominates the error of $\bar{\T}$. Therefore, $J_\mathfrak{E}$ has a similar behavior as $\bar{|v|}$. 
This is an explanation for the similar behaviour of the refined meshes. In all three cases, there is a strong refinement around the point source. 
This clover is a typical refinement structure around point sources. 
There is a big vortex in the center, two smaller vortices on the right side, 
c.f. Figure~\ref{fig_pick_4_out_of_6_A}. 
Moreover, there are smaller vortices in the left vertices of the square. This is similar to the driven cavity problem.
A more detailed picture of this part of the domain is displayed in 
Figure \ref{fig_pick_4_out_of_6_A} (right).

\begin{figure}[H]
	\begin{minipage}{0.49\linewidth}
		\ifMAKEPICS
		\begin{gnuplot}[terminal=epslatex,terminaloptions=color]
			set output "Figures/LaseronSquareMeanError.tex"
			set title "Laser on Square"
			set key top right
			set key opaque
			set datafile separator "|"
			set logscale x
			set logscale y
			set grid ytics lc rgb "#bbbbbb" lw 1 lt 0
			set grid xtics lc rgb "#bbbbbb" lw 1 lt 0
			set xlabel '\text{DOFs}'
			set format '
			plot \
			'< sqlite3 Data/LaseronSquare/New/Rhodaten/MeanTemp/data.db "SELECT DISTINCT DOFS_primal, abs(Exact_Error) from data WHERE DOFS_primal<200000"' u 1:2 w  lp lt 1 lw 3 title ' \footnotesize Error in $\bar{\T}$', \
			'< sqlite3 Data/LaseronSquare/New/Rhodaten/MeanTemp/data.db "SELECT DISTINCT DOFS_primal, abs(Estimated_Error) from data WHERE DOFS_primal<200000"' u 1:2 w  lp  lw 3 title ' \footnotesize Error Estimator of $\bar{\T}$', \
			1000/x**1.5  title '$\mathcal{O}(\text{DOFs}^-\frac{3}{2})$' lw  10									
			#					 '< sqlite3 Data/Multigoalp4/Higher_Order/dataHigherOrderJE.db "SELECT DISTINCT DOFs, abs(Exact_Error) from data "' u 1:2 w  lp lw 3 title ' \footnotesize Error in $J_\mathfrak{E}$', \
		\end{gnuplot}
		\fi
		{		\scalebox{0.55}{\input{Figures/LaseronSquareMeanError.tex}}}
		\caption{The error and error estimator for $\bar{\T}$.\label{Figure: Laser on Square Meantemperature}}
	\end{minipage}%
	\hfill
	\begin{minipage}{0.49\linewidth}
		
		\ifMAKEPICS
		\begin{gnuplot}[terminal=epslatex,terminaloptions=color]
			set output "Figures/LaseronSquareMeanIeff.tex"
			set title "Laser on Square"
			set key top left
			set key opaque
			set datafile separator "|"
			set logscale x
			set logscale y
			set grid ytics lc rgb "#bbbbbb" lw 1 lt 0
			set grid xtics lc rgb "#bbbbbb" lw 1 lt 0
			set xlabel '\text{DOFs}'
			set format '
			plot \
			'< sqlite3 Data/LaseronSquare/New/Rhodaten/MeanTemp/data.db "SELECT DISTINCT DOFS_primal, Ieff from data WHERE DOFS_primal>1000"' u 1:2 w  lp lt 1 lw 4 title ' \footnotesize $I_{eff}$ $\bar{\T}$', \
			'< sqlite3 Data/LaseronSquare/New/Rhodaten/MeanTemp/data.db "SELECT DISTINCT DOFS_primal, Ieff_primal from data WHERE DOFS_primal>1000"' u 1:2 w  lp lt 2 lw 2 title ' \footnotesize $I_{eff,p}$ $\bar{\T}$', \
			'< sqlite3 Data/LaseronSquare/New/Rhodaten/MeanTemp/data.db "SELECT DISTINCT DOFS_primal, Ieff_adjoint from data WHERE DOFS_primal>1000"' u 1:2 w  lp lt 3 lw 2 title ' \footnotesize $I_{eff,a}$ $\bar{\T}$', \
			1 lw 4										
			#					 '< sqlite3 Data/Multigoalp4/Higher_Order/dataHigherOrderJE.db "SELECT DISTINCT DOFs, abs(Exact_Error) from data "' u 1:2 w  lp lw 3 title ' \footnotesize Error in $J_\mathfrak{E}$', \
		\end{gnuplot}
		\fi
		{		\scalebox{0.55}{\input{Figures/LaseronSquareMeanIeff.tex}}}
		
		\caption{The \black{effectivity} indices for $\bar{\T}$. \label{Figure: Laser on Square Meantemperature Ieff}}
	\end{minipage}
\end{figure}

\begin{figure}[H]
	\begin{minipage}{0.49\linewidth}
		\ifMAKEPICS
		\begin{gnuplot}[terminal=epslatex,terminaloptions=color]
			set output "Figures/LaseronSquareMeanVelociityError.tex"
			set title "Laser on Square"
			set key top right
			set key opaque
			set datafile separator "|"
			set logscale x
			set logscale y
			set grid ytics lc rgb "#bbbbbb" lw 1 lt 0
			set grid xtics lc rgb "#bbbbbb" lw 1 lt 0
			set xlabel '\text{DOFs}'
			set format '
			plot \
			'< sqlite3 Data/LaseronSquare/New/Rhodaten/MeanVelo/data.db "SELECT DISTINCT DOFS_primal, abs(Exact_Error) from data WHERE DOFS_primal<200000"' u 1:2 w  lp lt 1 lw 3 title ' \footnotesize Error in $\bar{|v|}$', \
			'< sqlite3 Data/LaseronSquare/New/Rhodaten/MeanVelo/data.db "SELECT DISTINCT DOFS_primal, abs(Estimated_Error) from data WHERE DOFS_primal<200000"' u 1:2 w  lp  lw 3 title ' \footnotesize Error Estimator of $\bar{|v|}$', \
			1.0/x**1.5  title '$\mathcal{O}(\text{DOFs}^-\frac{3}{2})$' lw  10														
			#					 '< sqlite3 Data/Multigoalp4/Higher_Order/dataHigherOrderJE.db "SELECT DISTINCT DOFs, abs(Exact_Error) from data "' u 1:2 w  lp lw 3 title ' \footnotesize Error in $J_\mathfrak{E}$', \
		\end{gnuplot}
		\fi
		{		\scalebox{0.55}{\input{Figures/LaseronSquareMeanVelociityError.tex}}}
		\caption{The error and error estimator for $\bar{|v|}$.\label{Figure: Laser on Square Meanvelocity}}
	\end{minipage}%
	\hfill
	\begin{minipage}{0.49\linewidth}
		
		\ifMAKEPICS
		\begin{gnuplot}[terminal=epslatex,terminaloptions=color]
			set output "Figures/LaseronSquareMeanVelociityIeff.tex"
			set title "Laser on Square"
			set key top left
			set key opaque
			set datafile separator "|"
			set logscale x
			set logscale y
			set grid ytics lc rgb "#bbbbbb" lw 1 lt 0
			set grid xtics lc rgb "#bbbbbb" lw 1 lt 0
			set xlabel '\text{DOFs}'
			set format '
			plot \
			'< sqlite3 Data/LaseronSquare/New/Rhodaten/MeanVelo/data.db "SELECT DISTINCT DOFS_primal, Ieff from data WHERE DOFS_primal>1000"' u 1:2 w  lp lt 1 lw 4 title ' \footnotesize $I_{eff}$ $\bar{|v|}$', \
			'< sqlite3 Data/LaseronSquare/New/Rhodaten/MeanVelo/data.db "SELECT DISTINCT DOFS_primal, Ieff_primal from data WHERE DOFS_primal>1000"' u 1:2 w  lp lt 2 lw 2 title ' \footnotesize $I_{eff,p}$ $\bar{|v|}$', \
			'< sqlite3 Data/LaseronSquare/New/Rhodaten/MeanVelo/data.db "SELECT DISTINCT DOFS_primal, Ieff_adjoint from data WHERE DOFS_primal>1000"' u 1:2 w  lp lt 3 lw 2 title ' \footnotesize $I_{eff,a}$ $\bar{|v|}$', \
			1 lw 4	
			#					 '< sqlite3 Data/Multigoalp4/Higher_Order/dataHigherOrderJE.db "SELECT DISTINCT DOFs, abs(Exact_Error) from data "' u 1:2 w  lp lw 3 title ' \footnotesize Error in $J_\mathfrak{E}$', \
		\end{gnuplot}
		\fi
		{		\scalebox{0.55}{\input{Figures/LaseronSquareMeanVelociityIeff.tex}}}
		
		\caption{The \black{effectivity} indices for $\bar{|v|}$. \label{Figure: Laser on Square Meanvelocity Ieff}}
	\end{minipage}
\end{figure}

\begin{figure}[H]
	\begin{minipage}[t]{0.49\linewidth}
		\ifMAKEPICS
		\begin{gnuplot}[terminal=epslatex,terminaloptions=color]
			set output "Figures/LaseronSquareCombinedError.tex"
			set title "Laser on Square"
			set key top right
			set key opaque
			set datafile separator "|"
			set logscale x
			set logscale y
			set grid ytics lc rgb "#bbbbbb" lw 1 lt 0
			set grid xtics lc rgb "#bbbbbb" lw 1 lt 0
			set xlabel '\text{DOFs}'
			set format '
			plot \
			'< sqlite3 Data/LaseronSquare/New/Rhodaten/Combined/data.db	 "SELECT DISTINCT DOFS_primal, abs(Exact_Error) from data WHERE DOFS_primal<200000"' u 1:2 w  lp lt 1 lw 3 title ' \footnotesize Error in $J_\mathfrak{E}$', \
			'< sqlite3 Data/LaseronSquare/New/Rhodaten/Combined/data.db	 "SELECT DISTINCT DOFS_primal, abs(Estimated_Error) from data WHERE DOFS_primal<200000"' u 1:2 w  lp  lw 3 title ' \footnotesize Error Estimator of $J_\mathfrak{E}$', \
			'< sqlite3 Data/LaseronSquare/New/Rhodaten/Combined/data.db	 "SELECT DISTINCT DOFS_primal, abs(relativeError0) from data WHERE DOFS_primal<200000"' u 1:2 w  lp  lw 3 title ' \footnotesize Error in $\bar{|v|}$', \
			'< sqlite3 Data/LaseronSquare/New/Rhodaten/Combined/data.db	 "SELECT DISTINCT DOFS_primal, abs(relativeError1) from data WHERE DOFS_primal<200000"' u 1:2 w  lp  lw 3 title ' \footnotesize Error in $\bar{T}$', \
			1000000/x**1.5  title '$\mathcal{O}(\text{DOFs}^-\frac{3}{2})$' lw  10														
			#					 '< sqlite3 Data/Multigoalp4/Higher_Order/dataHigherOrderJE.db "SELECT DISTINCT DOFs, abs(Exact_Error) from data "' u 1:2 w  lp lw 3 title ' \footnotesize Error in $J_\mathfrak{E}$', \
		\end{gnuplot}
		\fi
		{		\scalebox{0.55}{\input{Figures/LaseronSquareCombinedError.tex}}}
		\caption{The error and error estimator for $J_\mathfrak{E}$ for $\bar{|v|}$ and $\bar{T}$.\label{Figure: Laser on Square Combined}}
	\end{minipage}%
	\hfill
	\begin{minipage}[t]{0.49\linewidth}	
		\ifMAKEPICS
		\begin{gnuplot}[terminal=epslatex,terminaloptions=color]
			set output "Figures/LaseronSquareCombinedIeff.tex"
			set title "Laser on Square"
			set key top right
			set key opaque
			set datafile separator "|"
			set logscale x
			set logscale y
			set grid ytics lc rgb "#bbbbbb" lw 1 lt 0
			set grid xtics lc rgb "#bbbbbb" lw 1 lt 0
			set xlabel '\text{DOFs}'
			set format '
			plot \
			'< sqlite3 Data/LaseronSquare/New/Rhodaten/Combined/data.db	 "SELECT DISTINCT DOFS_primal, Ieff from data WHERE DOFS_primal<200000"' u 1:2 w  lp lt 1 lw 4 title ' \footnotesize  $I_{eff}$ $J_\mathfrak{E}$', \
			'< sqlite3 Data/LaseronSquare/New/Rhodaten/Combined/data.db	 "SELECT DISTINCT DOFS_primal, Ieff_primal from data WHERE DOFS_primal<200000"' u 1:2 w  lp lt 2 lw 2 title ' \footnotesize  $I_{eff,p}$ $J_\mathfrak{E}$', \
			'< sqlite3 Data/LaseronSquare/New/Rhodaten/Combined/data.db	 "SELECT DISTINCT DOFS_primal, Ieff_adjoint from data WHERE DOFS_primal<200000"' u 1:2 w  lp lt 3 lw 2 title ' \footnotesize  $I_{eff,a}$ $J_\mathfrak{E}$', \
			1 lw 4						
			#					 '< sqlite3 Data/Multigoalp4/Higher_Order/dataHigherOrderJE.db "SELECT DISTINCT DOFs, abs(Exact_Error) from data "' u 1:2 w  lp lw 3 title ' \footnotesize Error in $J_\mathfrak{E}$', \
		\end{gnuplot}
		\fi
		{		\scalebox{0.55}{\input{Figures/LaseronSquareCombinedIeff.tex}}}
		
		\caption{The \black{effectivity} indices  for $J_\mathfrak{E}$ for $\bar{|v|}$ and $\bar{T}$.\label{Figure: Laser on Square Combined Ieff}}
	\end{minipage}
	%
	%
\end{figure}
%
%

\subsection{Y-beam splitter}
\label{sec_Y_beam}
In this third example, we consider 
a Y-beam splitter with a laser that generates a flow field due to gravity. 
This application is motivated from 
PhoenixD\footnote{\url{https://www.phoenixd.uni-hannover.de/en/}}
and is an important configuration in wave guide modeling \cite{Peetal21}.

\subsubsection{Configuration, geometry, parameters, boundary conditions}
\begin{figure}[t]
	\begin{minipage}[t]{0.49\linewidth}
		\includegraphics[scale=0.34]{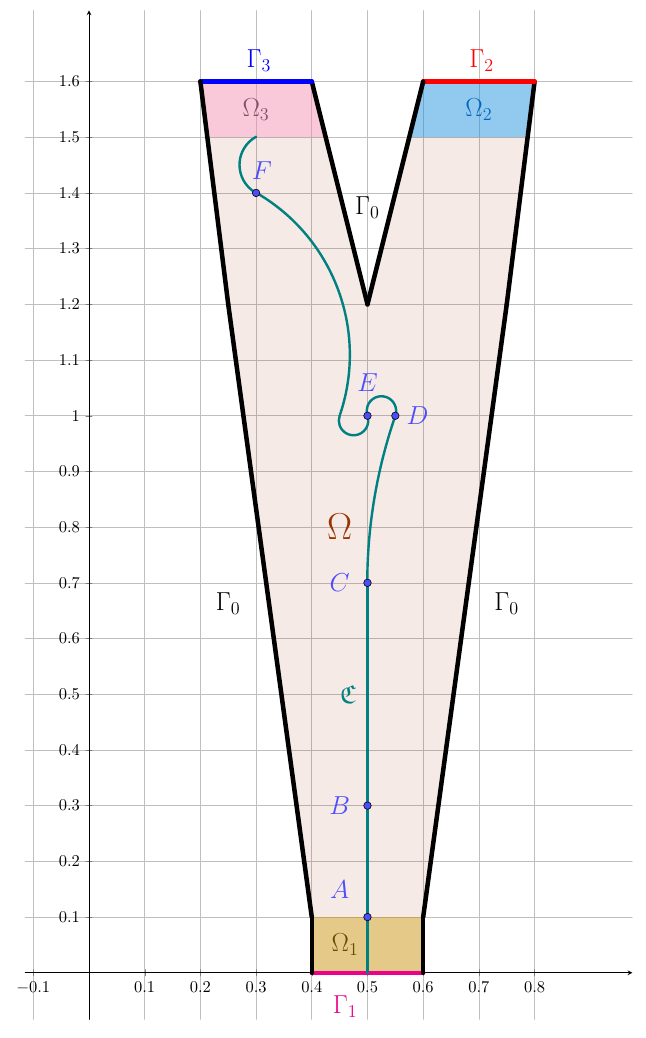}
		\caption{The domain $\Omega$, the boundary parts, the domains for our quantites of interest, the curve $\mathfrak{C}$ and the points for the different configurations.} 
		\label{YBeamsplitterData}
	\end{minipage}%
	\hfill
	\begin{minipage}[t]{0.49\linewidth}
		\centering
		\includegraphics[scale=0.29]{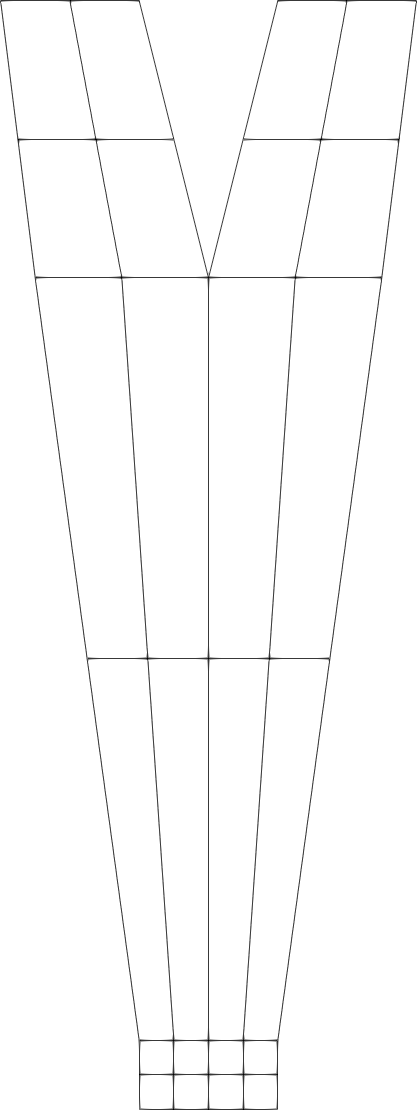}	
		\caption{The initial mesh for all configurations.}
	\end{minipage}%
\end{figure}

The domain $\Omega$ \black{and its subdomains $\Omega_1, \Omega_2, \Omega_3$ are} depicted in Figure~\ref{YBeamsplitterData}. 
Furthermore, we have the fluid flow boundary conditions	
\begin{align*}
v=&0   \qquad\text{on } \Gamma_{0}  = \Gamma_{vD}, \\
\nu (\T) \frac{\partial v}{\partial n} - p \cdot n=& 0  
\qquad\text{on } \Gamma_{1} \cup \Gamma_{2} \cup \Gamma_{3} = \Gamma_{vN}, \nonumber\\
\end{align*}
the temperature boundary conditions
\begin{align*}
\T = & 293.15   \qquad\text{on } \Gamma_{1} = \Gamma_{\T D}, \nonumber \\
\frac{\partial \T}{\partial n}  = & 0  \qquad \qquad\text{on } 
\partial \Omega \setminus \Gamma_{1} = \Gamma_{\T N}. \nonumber
\end{align*}
The right hand side $f$ is chosen as
\begin{itemize}
	\item Configuration 1: $f(x):= \psi_{A}(x)$, $A=(0.5,0.1)$,
	\item Configuration 2: $f(x):= \frac{1}{2}\psi_{B}(x)+\frac{1}{2}\psi_{C}(x)$, $B=(0.5,0.3)$, $C=(0.5,0.7)$,
	\new{
		\item Configuration 3: $f(x):= \psi_{E}(x)$, $E=(0.5,1)$,
		\item Configuration 4: $f(x):= \psi_{D}(x)$, $D=(0.55,1)$,
	}
	\item Configuration 5: $f(x):= \psi_{F}(x)$, $F=(0.3,1.4)$,
	\item Configuration 6: $f(x):= \psi_{\mathfrak{C}}(x)$.
\end{itemize}
The function $$\psi_{x_0}(x):=10^4\sqrt{2\pi}e^{-{10^4|x-x_0|^2}},$$ resembles the laser centred at $x_0\in\{A, B, C, D, E, F\}$, and the function $$\psi_{\mathfrak{C}}(x):=\frac{500}{|\mathfrak{C}|}\sqrt{2\pi}e^{{{-10^4}dist(x,\mathfrak{C})^2}},$$ resembles the laser along the curve $\mathfrak{C}$ where $dist(x,\mathfrak{C}):=\inf_{x^*\in \mathfrak{C}} |x-x^*|$ and $|\mathfrak{C}|$ is the length of the curve $\mathfrak{C}$.
The curve $\mathfrak{C}$ is given by circular arcs, with continuous tangent through the points $((0.5,0),A,B,C,D,E,(0.45,1),F,(0.3,1.5))$, where we start with a line between  $(0.5,0)$ and $A$.

\subsubsection{Goal functionals}
\label{sec_Y_beam_goal}
We consider seven goal functionals:
\begin{figure}[H]
	\begin{itemize}	
		\begin{minipage}[t]{0.49\linewidth}
			\item[]  $J_1(U):= \int_{\Gamma_{1}}{v\cdot n}\dx $,
			\item[]  $J_3(U):= \int_{\Gamma_{3}}{v\cdot n}\dx $,
			\item[]  $J_5(U):=  \frac{1}{|\Omega_{2}|}\int_{\Omega_{2}} \T \dx $,
			\item[]  $J_7(U):= (J_5(U)-J_6(U))^2$.
		\end{minipage}%
		\hfill
		\begin{minipage}[t]{0.49\linewidth}
			\item[]  $J_2(U):= \int_{\Gamma_{2}}{v\cdot n}\dx $,
			\item[]  $J_4(U):= \frac{1}{|\Omega_{1}|} \int_{\Omega_{1}} \T \dx $,
			\item[]  $J_6(U):= \frac{1}{|\Omega_{3}|}\int_{\Omega_{3}} \T  \dx $,
		\end{minipage}%
	\end{itemize}
\end{figure}
\black{\begin{Remark}
		Due to symmetry $J_7$ vanishes for Configuration 1-3.
\end{Remark}}

\subsubsection{Discussion of our findings}

\begin{figure}[H]
	\begin{minipage}{1.0 \linewidth}
		{
			\captionsetup{labelformat=empty}
			\begin{minipage}{0.16\linewidth}
				\centering
				\includegraphics[scale=0.104]{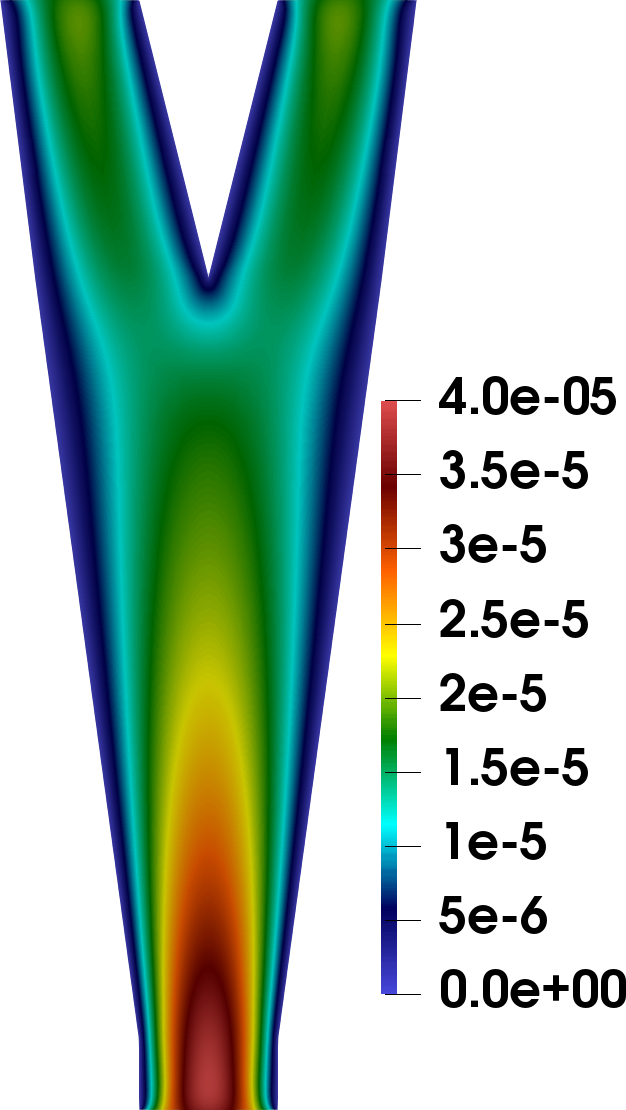}
				\includegraphics[scale=0.104]{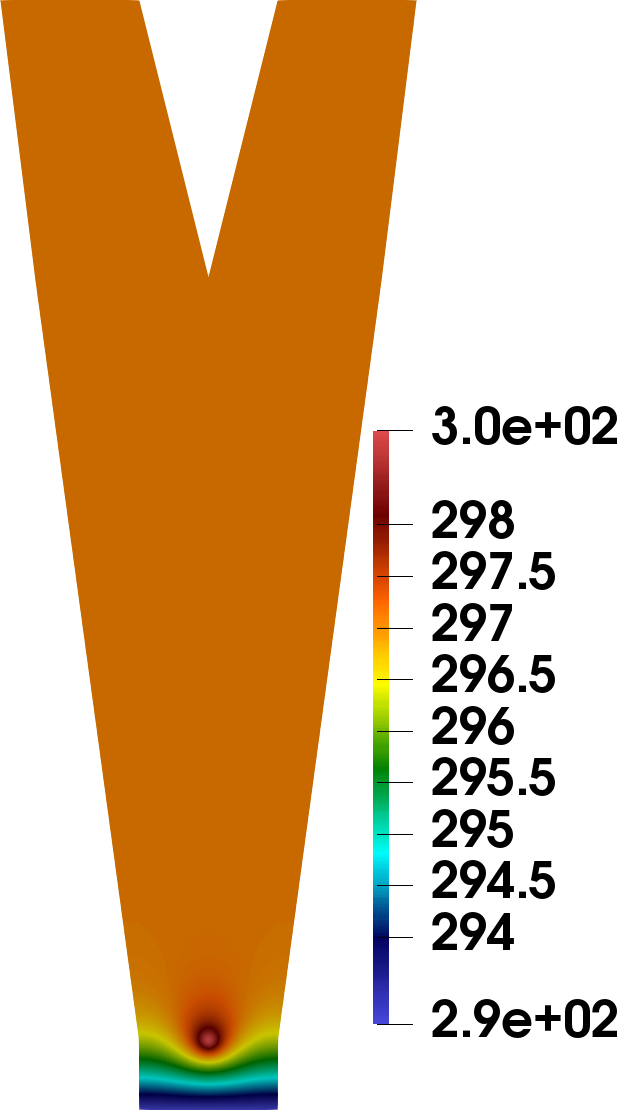}
			\end{minipage}%
			\hfill
			\begin{minipage}{0.16\linewidth}	
				\centering
				\includegraphics[scale=0.14412]{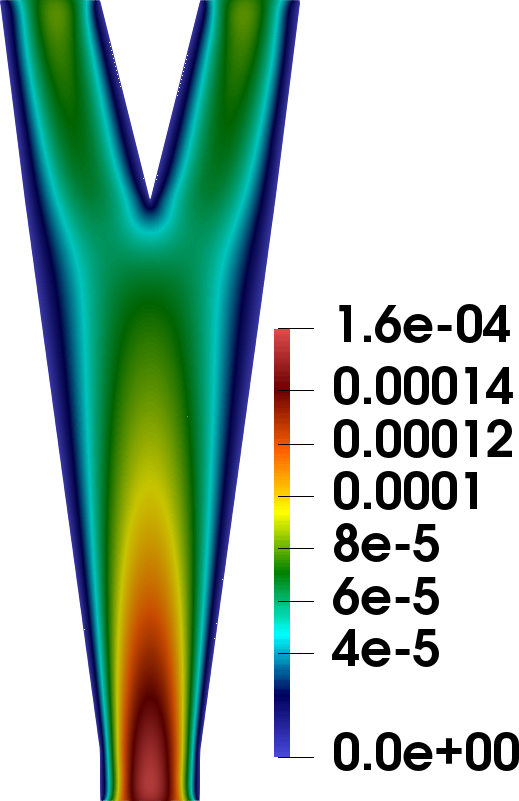}
				\includegraphics[scale=0.14412]{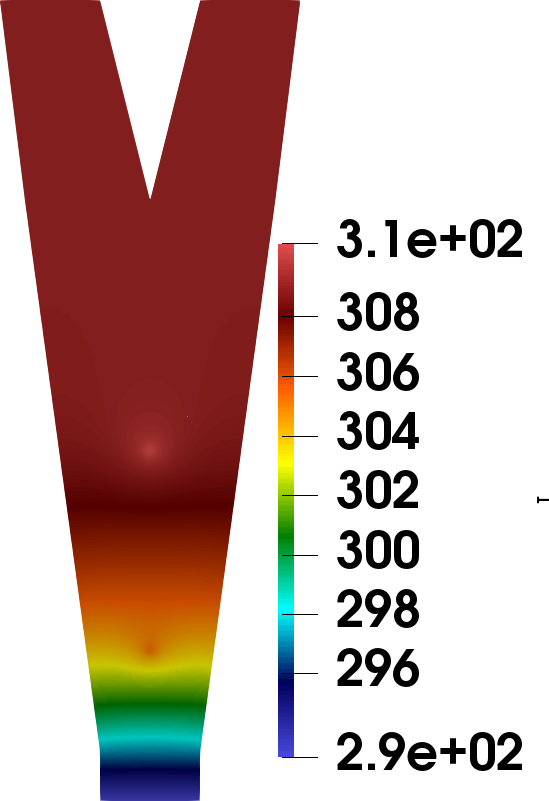}
			\end{minipage}%
			\hfill
			\begin{minipage}{0.16\linewidth}	
				\centering
				\includegraphics[scale=0.104]{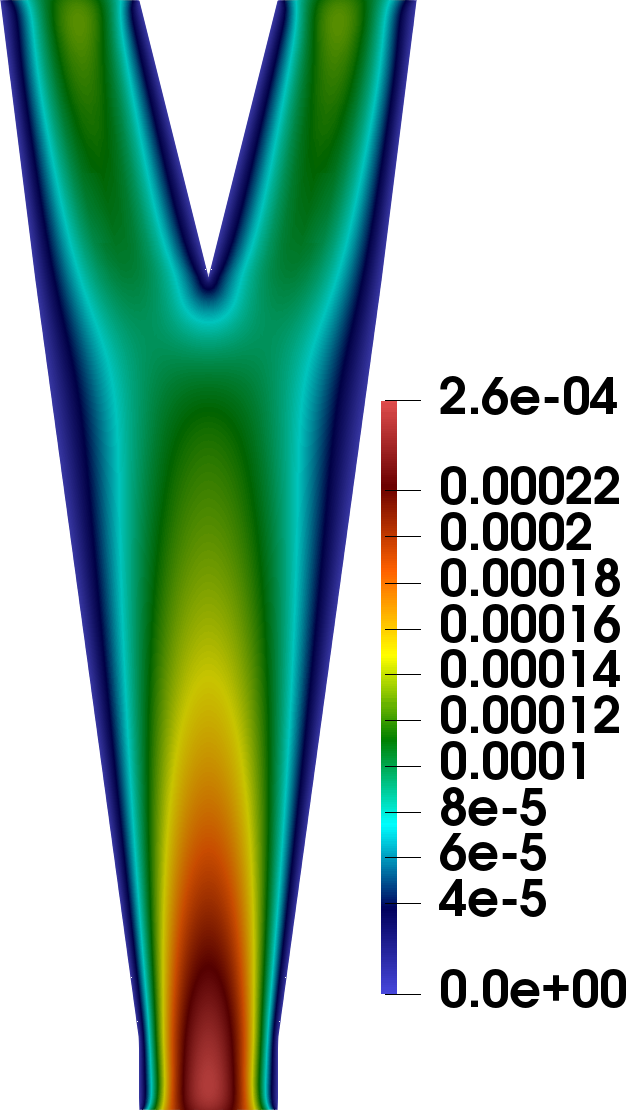}
				\includegraphics[scale=0.104]{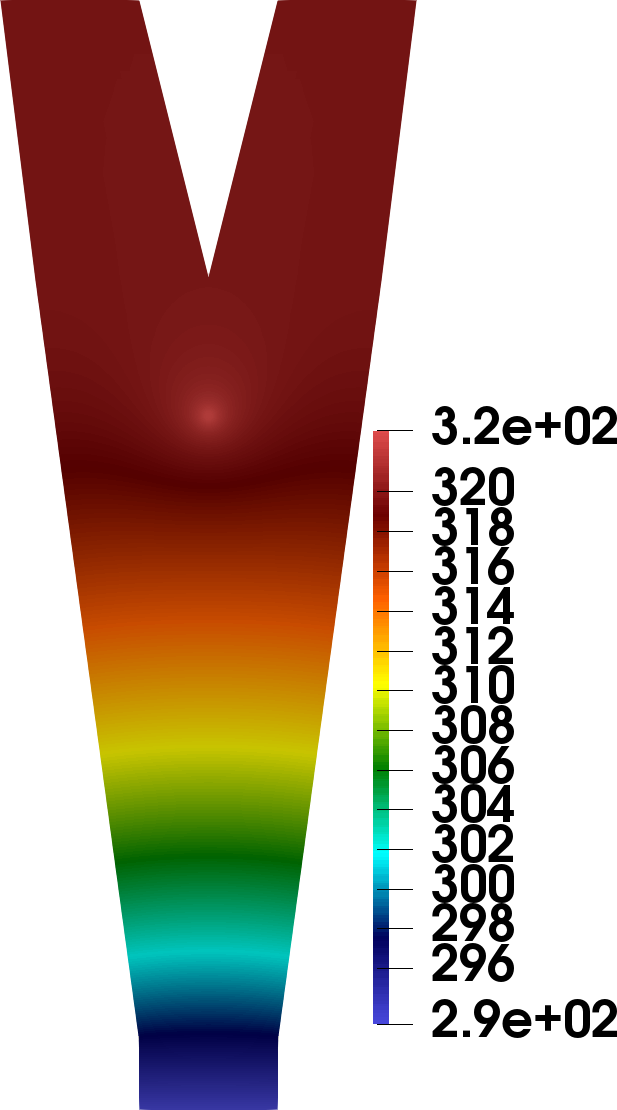}
			\end{minipage}%
			\hfill
			\begin{minipage}{0.16\linewidth}
				\centering
				\includegraphics[scale=0.104]{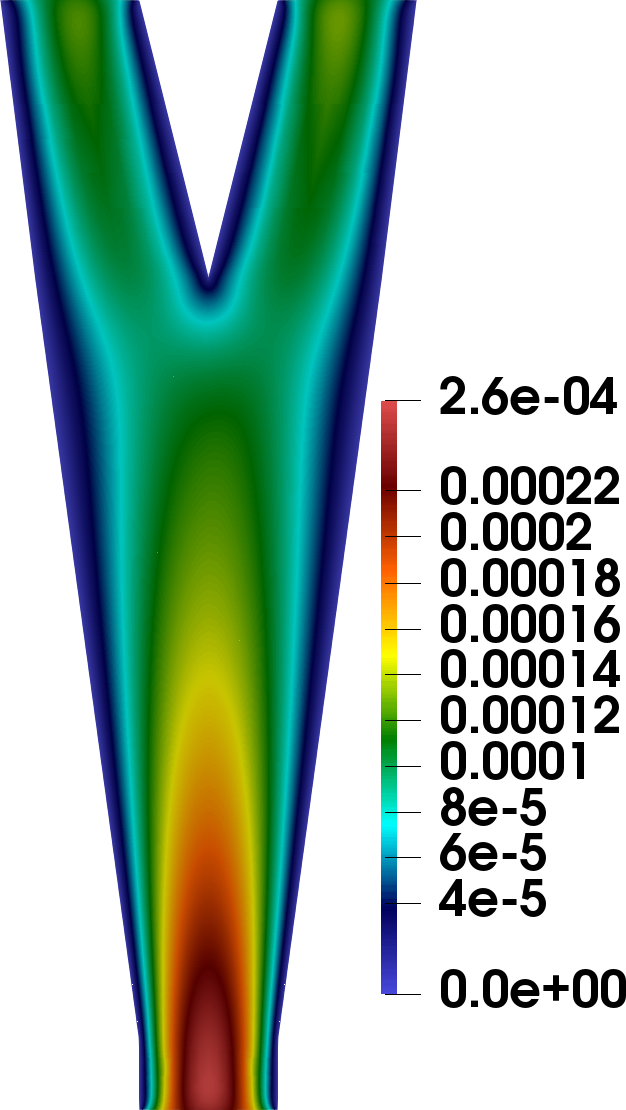}
				\includegraphics[scale=0.104]{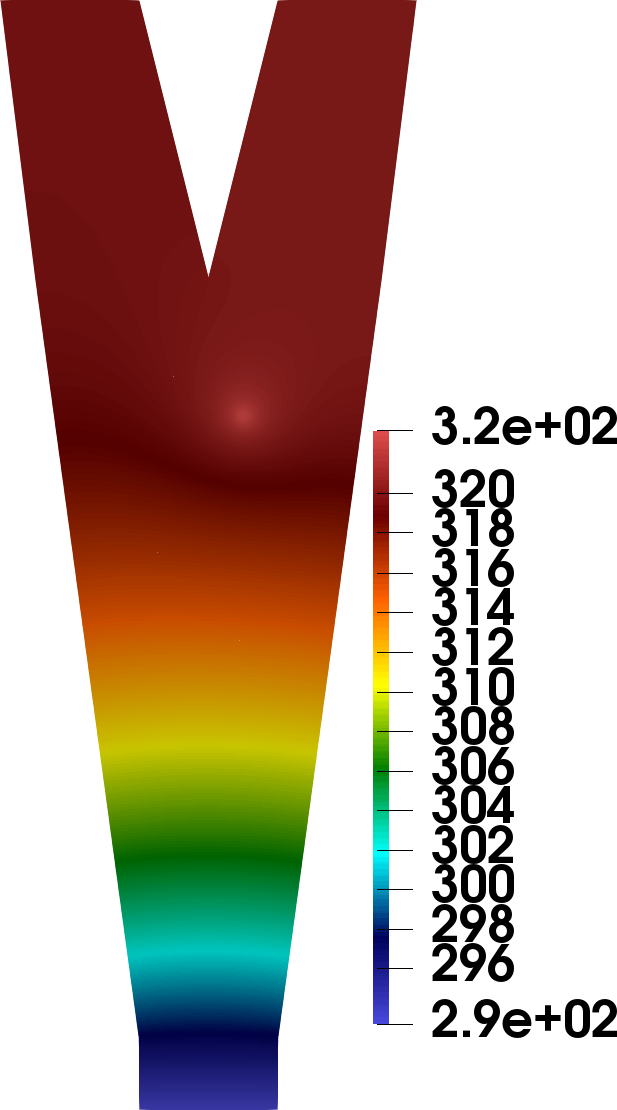}
			\end{minipage}%
			\hfill
			\begin{minipage}{0.16\linewidth}	
				\centering
				\includegraphics[scale=0.104]{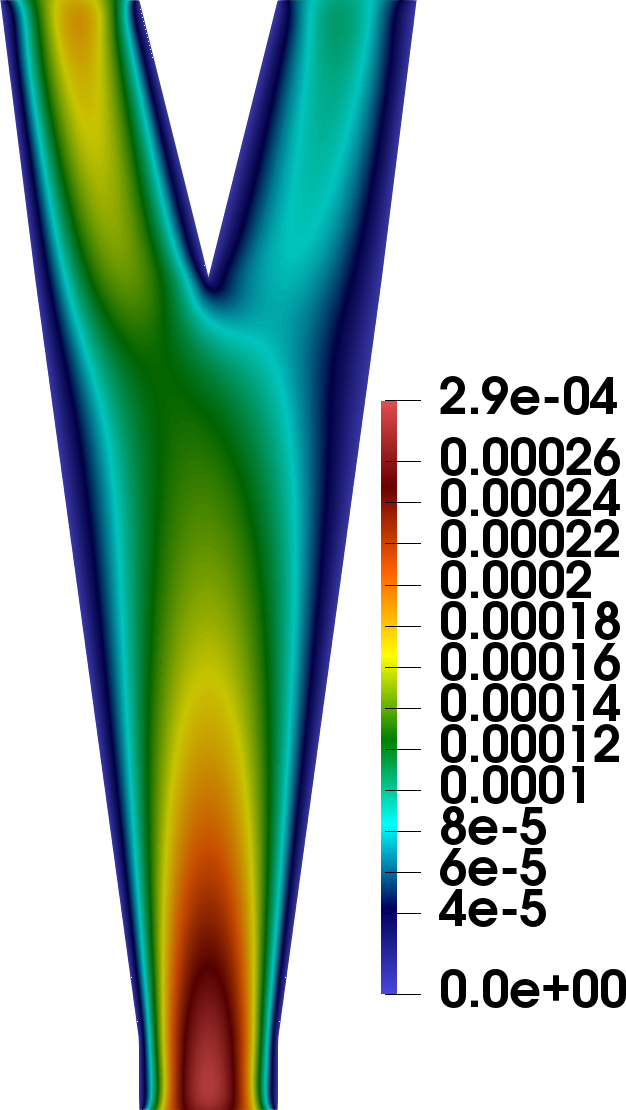}
				\includegraphics[scale=0.104]{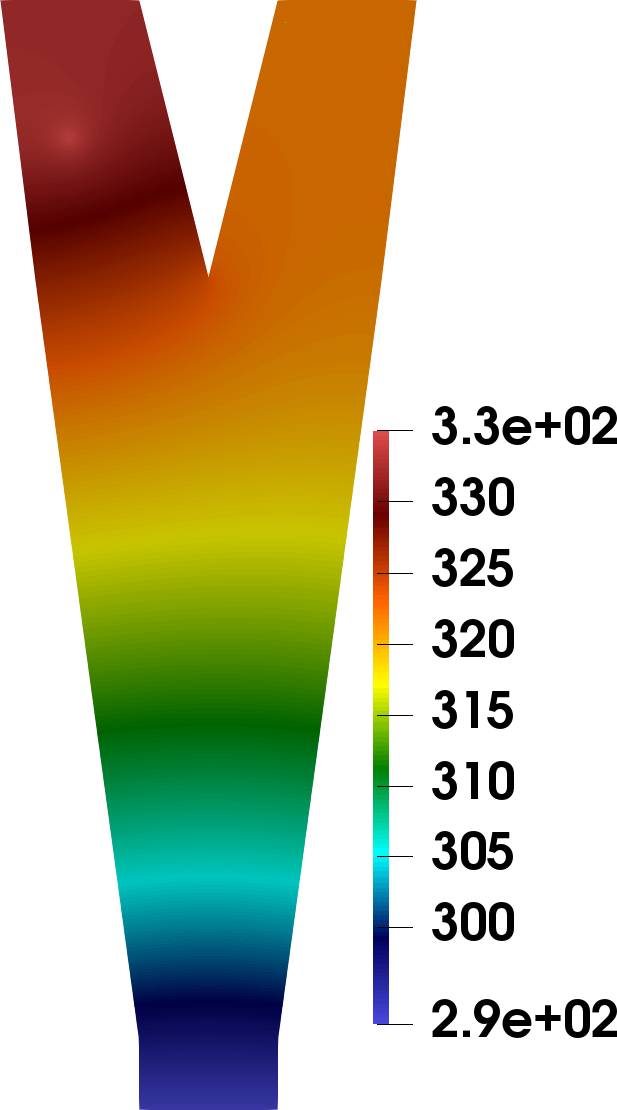}
			\end{minipage}%
			\hfill
			\begin{minipage}{0.16\linewidth}	
				\centering
				\includegraphics[scale=0.13376]{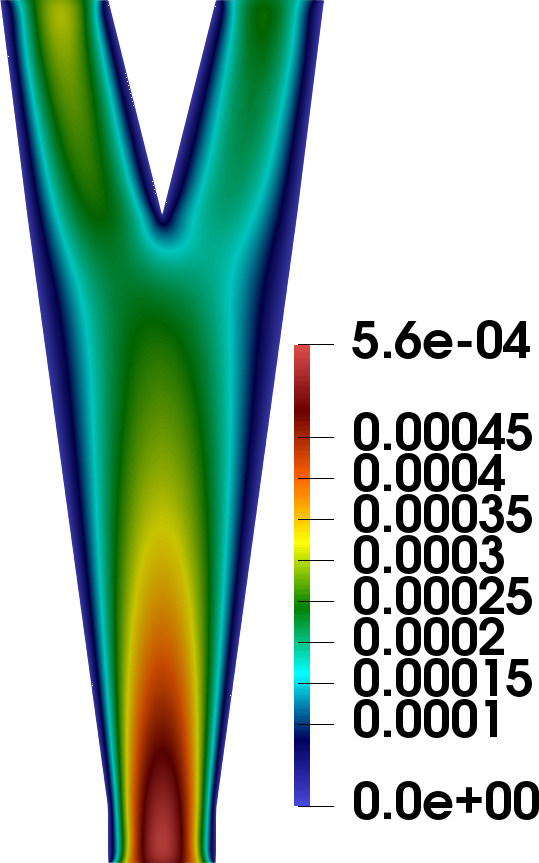}
				\includegraphics[scale=0.13376]{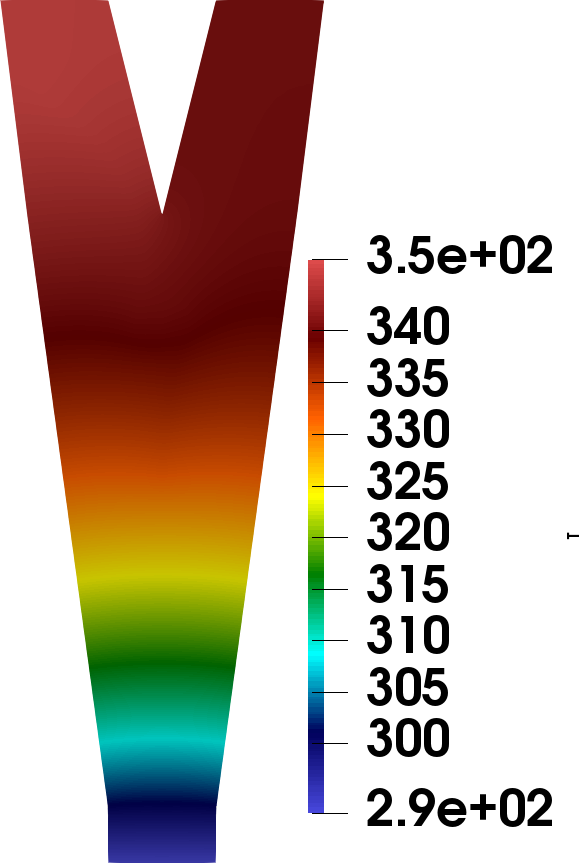}
		\end{minipage}}
		\caption{The magnitute of the velocity (above) and the temperature (below) 
			for Configuration 1-6 from left to right.
			\label{fig_test_3_a}}
		
	\end{minipage}

\end{figure}

\begin{figure}[H]
	{
		\captionsetup{labelformat=empty}
		\begin{minipage}{0.16\linewidth}
			\centering
			\includegraphics[scale=0.15]{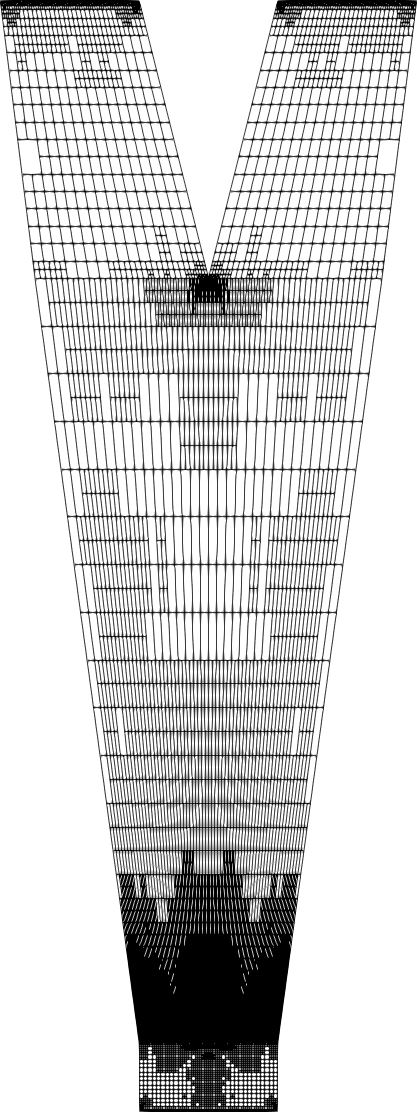}
		\end{minipage}%
		\hfill
		\begin{minipage}{0.16\linewidth}
			\centering
			\includegraphics[scale=0.193055555555]{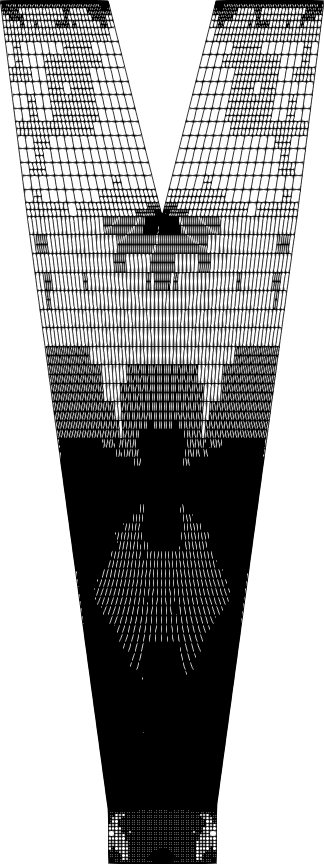}
		\end{minipage}%
		\hfill
		\begin{minipage}{0.16\linewidth}	
			\centering
			\includegraphics[scale=0.15]{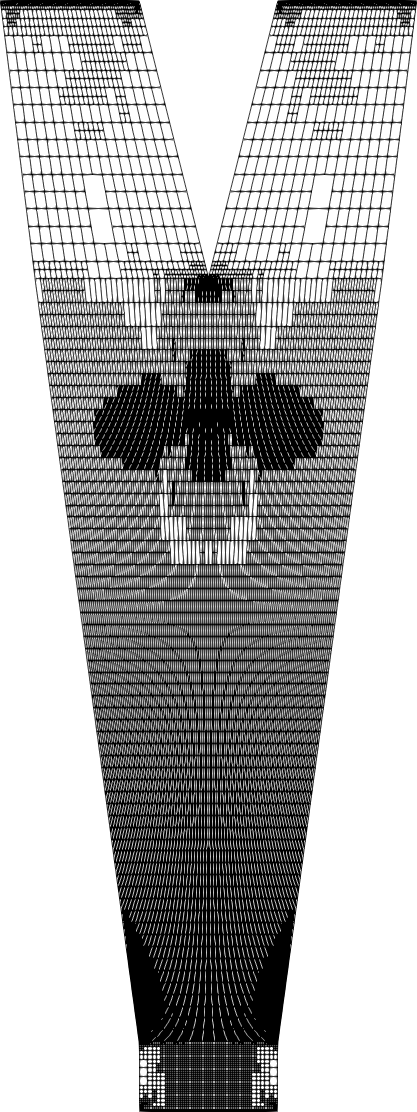}
		\end{minipage}%
		\hfill
		\begin{minipage}{0.16\linewidth}
			\centering
			\includegraphics[scale=0.15]{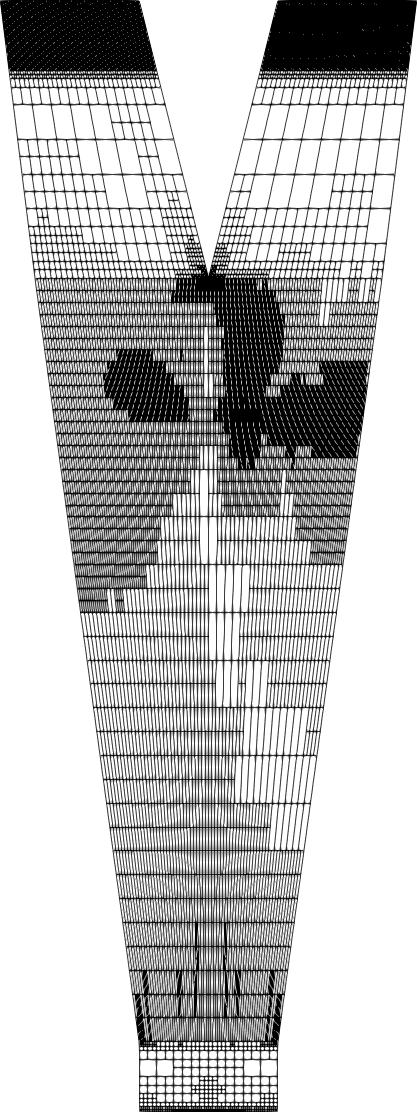}
		\end{minipage}%
		\hfill
		\begin{minipage}{0.16\linewidth}	
			\centering
			\includegraphics[scale=0.15]{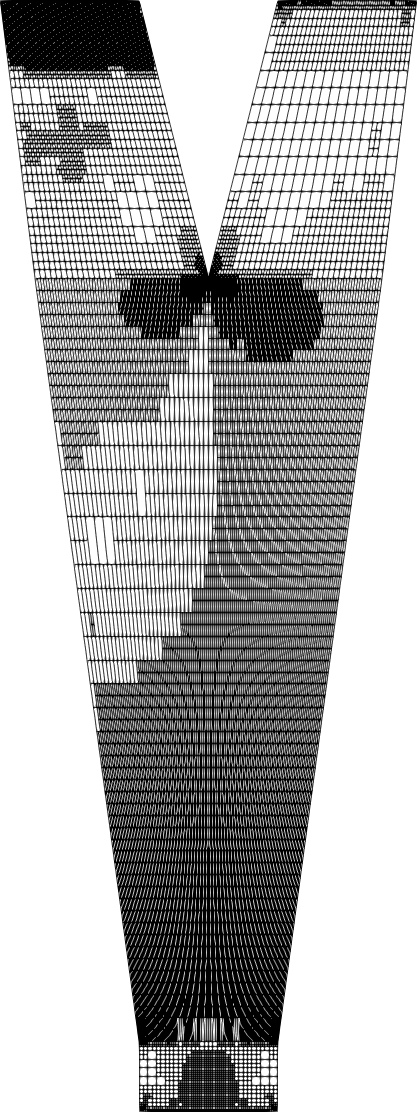}
		\end{minipage}
		\begin{minipage}{0.16\linewidth}
			\centering
			\includegraphics[scale=0.193055555555]{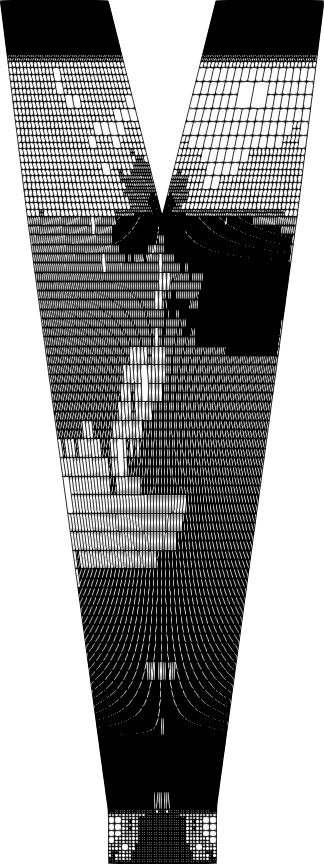}
		\end{minipage}%
		\hfill
		
	}
	\caption{The inital and final mesh 
		for Configuration 1-6 from left to right.
	}
	\label{fig_test_3_b}
\end{figure}

The magnitute of the velocity and the temperature for different configurations
are displayed in Figure \ref{fig_test_3_a}. Furthermore, 
different locally refined meshes are shown in Figure \ref{fig_test_3_b}. These show 
refinement in geometric singularities such as the kink where the splitter branches, but 
as well local refinement due to the goal functionals. 

In the following, we 
discuss our quantitative findings. First, the reference values are provided 
in Table \ref{table_test_3_a}. The errors and error estimators 
versus the degrees of freedom are displayed in the 
Figures \ref{fig_test_3_ca} to \ref{fig_test_3_e}.
The various effectivity indices are shown in the 
Figures \ref{fig_test_3_fa} to \ref{fig_test_3_h},
from which we observe excellent 
performances. This is in particular remarkable due to the different configurations 
and the nonlinear behavior of the coupled PDE system. It can be inferred that 
our multigoal framework is robust and yields a cost-efficient numerical procedure.

\begin{table}[H]
	\centering
	\begin{tabular}{|c|c|c|c|}\hline
		Config. &    1         &         2         &      3           \\\hline
		$J_1$	& -5.1630481e-06 &-2.1224898e-05& -3.36248616e-05  \\ \hline
		$J_2$	& 2.5815241e-06  &1.0612449e-05& 1.68124308e-05  \\ \hline
		$J_3$	& 2.5815241e-06  &1.0612449e-05& 1.68124308e-05 \\ \hline
		$J_4$	& 2.9511378e+02  &2.9511878e+02& 2.95118918e+02 \\ \hline
		$J_5$	& 2.9704640e+02  &3.0918292e+02& 3.20253892e+02\\ \hline
		$J_6$	& 2.9704640e+02  &3.0918292+02& 3.20253892e+02 \\ \hline
		$J_7$	& 0 & 0 & 0  \\ \hline
	\end{tabular}\\
	\vspace{0.5cm}
	\begin{tabular}{|c|c|c|c|}\hline
		Config. &         4             &             5     &  6    \\\hline
		$J_1$	 & -3.3629129e-05 & -3.7522141e-05 &-7.2595565e-05  \\ \hline
		$J_2$	& 1.7183604e-05  & 1.3366980e-05  &3.2117505e-05\\ \hline
		$J_3$	& 1.6445525e-05  & 2.4155161e-05  &4.0478060e-05\\ \hline
		$J_4$& 2.9511892e+02  & 2.9511902e+02  &2.9670043e+02\\ \hline
		$J_5$	& 3.2052058e+02  & 3.2170400e+02  &3.4082013e+02\\ \hline
		$J_6$	& 3.2001790e+02  & 3.3282825e+02  &3.4671220e+02\\ \hline
		$J_7$	& 2.5268810e-01  & 1.237489e+02 &3.4716466e+01  \\ \hline
	\end{tabular}
	\caption{Reference values
	}
	\label{table_test_3_a}
\end{table}

\begin{figure}[H]
	\begin{minipage}{0.49\linewidth}
		\ifMAKEPICS
		\begin{gnuplot}[terminal=epslatex,terminaloptions=color]
			set output "Figures/YSplitterErrorConfig1.tex"
			set title "Configuration 1"
			set rmargin 15
			set key outside
			set key right top
			set key opaque
			set datafile separator "|"
			set logscale x
			set logscale y
			set grid ytics lc rgb "#bbbbbb" lw 1 lt 0
			set grid xtics lc rgb "#bbbbbb" lw 1 lt 0
			set xlabel '\text{DOFs}'
			set format '
			plot \
			'< sqlite3 Data/Ybeam/ysplitterConfig1.db "SELECT DISTINCT DOFS_primal, abs(Exact_Error) from data WHERE DOFS_primal>100"' u 1:2 w  lp lt 1 lw 3 title ' \footnotesize Error in $J_c$', \
			'< sqlite3 Data/Ybeam/ysplitterConfig1.db "SELECT DISTINCT DOFS_primal, abs(Estimated_Error) from data WHERE DOFS_primal>100"' u 1:2 w  lp  lw 3 title ' \footnotesize  $\eta_h$', \
			'< sqlite3 Data/Ybeam/ysplitterConfig1.db "SELECT DISTINCT DOFS_primal, abs(relativeError0) from data WHERE DOFS_primal>100"' u 1:2 w  lp  lw 3 title ' \footnotesize Error in $J_1$', \
			'< sqlite3 Data/Ybeam/ysplitterConfig1.db "SELECT DISTINCT DOFS_primal, abs(relativeError1) from data WHERE DOFS_primal>100"' u 1:2 w  lp  lw 3 title ' \footnotesize Error in $J_2$', \
			'< sqlite3 Data/Ybeam/ysplitterConfig1.db "SELECT DISTINCT DOFS_primal, abs(relativeError3) from data WHERE DOFS_primal>100"' u 1:2 w  lp  lw 3 title ' \footnotesize Error in $J_4$', \
			'< sqlite3 Data/Ybeam/ysplitterConfig1.db "SELECT DISTINCT DOFS_primal, abs(relativeError5) from data WHERE DOFS_primal>100"' u 1:2 w  lp lc rgb 'orange' lw 3 title ' \footnotesize Error in $J_5$', \
			10/x**1.0 lw  10 lt 10	 title ' \footnotesize $\text{DOFs}^{-1}$'											
			#					 '< sqlite3 Data/Multigoalp4/Higher_Order/dataHigherOrderJE.db "SELECT DISTINCT DOFs, abs(Exact_Error) from data "' u 1:2 w  lp lw 3 title ' \footnotesize Error in $J_\mathfrak{E}$', \
		\end{gnuplot}
		\fi
		{		\scalebox{0.55}{\input{Figures/YSplitterErrorConfig1.tex}}}
		\caption{The error and error estimator for Configuration 1.}
		\label{fig_test_3_ca}
	\end{minipage}%
	\hfill
	\begin{minipage}{0.49\linewidth}
		\ifMAKEPICS
		\begin{gnuplot}[terminal=epslatex,terminaloptions=color]
			set output "Figures/YSplitterErrorConfig2.tex"
			set title "Configuration 2"
			set rmargin 15
			set key outside
			set key right top
			set key opaque
			set datafile separator "|"
			set logscale x
			set logscale y
			set grid ytics lc rgb "#bbbbbb" lw 1 lt 0
			set grid xtics lc rgb "#bbbbbb" lw 1 lt 0
			set xlabel '\text{DOFs}'
			set format '
			plot \
			'< sqlite3 Data/Ybeam/ysplitterConfig5.db "SELECT DISTINCT DOFS_primal, abs(Exact_Error) from data WHERE DOFS_primal>100"' u 1:2 w  lp lt 1 lw 3 title ' \footnotesize Error in $J_c$', \
			'< sqlite3 Data/Ybeam/ysplitterConfig5.db "SELECT DISTINCT DOFS_primal, abs(Estimated_Error) from data WHERE DOFS_primal>100"' u 1:2 w  lp  lw 3 title ' \footnotesize  $\eta_h$', \
			'< sqlite3 Data/Ybeam/ysplitterConfig5.db "SELECT DISTINCT DOFS_primal, abs(relativeError0) from data WHERE DOFS_primal>100"' u 1:2 w  lp  lw 3 title ' \footnotesize Error in $J_1$', \
			'< sqlite3 Data/Ybeam/ysplitterConfig5.db "SELECT DISTINCT DOFS_primal, abs(relativeError1) from data WHERE DOFS_primal>100"' u 1:2 w  lp  lw 3 title ' \footnotesize Error in $J_2$', \
			'< sqlite3 Data/Ybeam/ysplitterConfig5.db "SELECT DISTINCT DOFS_primal, abs(relativeError3) from data WHERE DOFS_primal>100"' u 1:2 w  lp  lw 3 title ' \footnotesize Error in $J_4$', \
			'< sqlite3 Data/Ybeam/ysplitterConfig5.db "SELECT DISTINCT DOFS_primal, abs(relativeError5) from data WHERE DOFS_primal>100"' u 1:2 w  lp  lc rgb 'orange' lw 3 title ' \footnotesize Error in $J_5$', \
			10/x**1.0 lw  10 lt 10	 title ' \footnotesize $\text{DOFs}^{-1}$'									
			#					 '< sqlite3 Data/Multigoalp4/Higher_Order/dataHigherOrderJE.db "SELECT DISTINCT DOFs, abs(Exact_Error) from data "' u 1:2 w  lp lw 3 title ' \footnotesize Error in $J_\mathfrak{E}$', \
		\end{gnuplot}
		\fi
		{		\scalebox{0.55}{\input{Figures/YSplitterErrorConfig2.tex}}}
		\caption{The error and error estimator for Configuration 2.}
		\label{fig_test_3_c}
	\end{minipage}
\end{figure}

\begin{figure}[H]
	\begin{minipage}{0.49\linewidth}
		\ifMAKEPICS
		\begin{gnuplot}[terminal=epslatex,terminaloptions=color]
			set output "Figures/YSplitterErrorConfig3.tex"
			set title "Configuration 3"
			set rmargin 15
			set key outside
			set key right top
			set key opaque
			set datafile separator "|"
			set logscale x
			set logscale y
			set grid ytics lc rgb "#bbbbbb" lw 1 lt 0
			set grid xtics lc rgb "#bbbbbb" lw 1 lt 0
			set xlabel '\text{DOFs}'
			set format '
			plot \
			'< sqlite3 Data/Ybeam/ysplitterConfig2.db "SELECT DISTINCT DOFS_primal, abs(Exact_Error) from data WHERE DOFS_primal>100"' u 1:2 w  lp lt 1 lw 3 title ' \footnotesize Error in $J_c$', \
			'< sqlite3 Data/Ybeam/ysplitterConfig2.db "SELECT DISTINCT DOFS_primal, abs(Estimated_Error) from data WHERE DOFS_primal>100"' u 1:2 w  lp  lw 3 title ' \footnotesize  $\eta_h$', \
			'< sqlite3 Data/Ybeam/ysplitterConfig2.db "SELECT DISTINCT DOFS_primal, abs(relativeError0) from data WHERE DOFS_primal>100"' u 1:2 w  lp  lw 3 title ' \footnotesize Error in $J_1$', \
			'< sqlite3 Data/Ybeam/ysplitterConfig2.db "SELECT DISTINCT DOFS_primal, abs(relativeError1) from data WHERE DOFS_primal>100"' u 1:2 w  lp  lw 3 title ' \footnotesize Error in $J_2$', \
			'< sqlite3 Data/Ybeam/ysplitterConfig2.db "SELECT DISTINCT DOFS_primal, abs(relativeError2) from data WHERE DOFS_primal>100"' u 1:2 w  lp  lw 3 title ' \footnotesize Error in $J_3$', \
			'< sqlite3 Data/Ybeam/ysplitterConfig2.db "SELECT DISTINCT DOFS_primal, abs(relativeError3) from data WHERE DOFS_primal>100"' u 1:2 w  lp lc rgb 'orange'  lw 3 title ' \footnotesize Error in $J_4$', \
			'< sqlite3 Data/Ybeam/ysplitterConfig2.db "SELECT DISTINCT DOFS_primal, abs(relativeError4) from data WHERE DOFS_primal>100"' u 1:2 w  lp  lw 3 title ' \footnotesize Error in $J_5$', \
			'< sqlite3 Data/Ybeam/ysplitterConfig2.db "SELECT DISTINCT DOFS_primal, abs(relativeError5) from data WHERE DOFS_primal>100"' u 1:2 w  lp  lw 3 title ' \footnotesize Error in $J_6$', \
			10/x**1.0 lw  10	 title ' \footnotesize $\text{DOFs}^{-1}$'										
			#					 '< sqlite3 Data/Multigoalp4/Higher_Order/dataHigherOrderJE.db "SELECT DISTINCT DOFs, abs(Exact_Error) from data "' u 1:2 w  lp lw 3 title ' \footnotesize Error in $J_\mathfrak{E}$', \
		\end{gnuplot}
		\fi
		{		\scalebox{0.55}{\input{Figures/YSplitterErrorConfig3.tex}}}
		\caption{The errors and error estimator for Configuration 3.}
	\end{minipage}%
	\hfill
	\begin{minipage}{0.49\linewidth}
		\ifMAKEPICS
		\begin{gnuplot}[terminal=epslatex,terminaloptions=color]
			set output "Figures/YSplitterErrorConfig4.tex"
			set title "Configuration 4"
			set rmargin 15
			set key outside
			set key right top
			set key opaque
			set datafile separator "|"
			set logscale x
			set logscale y
			set grid ytics lc rgb "#bbbbbb" lw 1 lt 0
			set grid xtics lc rgb "#bbbbbb" lw 1 lt 0
			set xlabel '\text{DOFs}'
			set format '
			plot \
			'< sqlite3 Data/Ybeam/ysplitterConfig3.db "SELECT DISTINCT DOFS_primal, abs(Exact_Error) from data WHERE DOFS_primal>100"' u 1:2 w  lp lt 1 lw 3 title ' \footnotesize Error in $J_c$', \
			'< sqlite3 Data/Ybeam/ysplitterConfig3.db "SELECT DISTINCT DOFS_primal, abs(Estimated_Error) from data WHERE DOFS_primal>100"' u 1:2 w  lp  lw 3 title ' \footnotesize  $\eta_h$', \
			'< sqlite3 Data/Ybeam/ysplitterConfig3.db "SELECT DISTINCT DOFS_primal, abs(relativeError0) from data WHERE DOFS_primal>100"' u 1:2 w  lp  lw 3 title ' \footnotesize Error in $J_1$', \
			'< sqlite3 Data/Ybeam/ysplitterConfig3.db "SELECT DISTINCT DOFS_primal, abs(relativeError1) from data WHERE DOFS_primal>100"' u 1:2 w  lp  lw 3 title ' \footnotesize Error in $J_2$', \
			'< sqlite3 Data/Ybeam/ysplitterConfig3.db "SELECT DISTINCT DOFS_primal, abs(relativeError2) from data WHERE DOFS_primal>100"' u 1:2 w  lp  lw 3 title ' \footnotesize Error in $J_3$', \
			'< sqlite3 Data/Ybeam/ysplitterConfig3.db "SELECT DISTINCT DOFS_primal, abs(relativeError3) from data WHERE DOFS_primal>100"' u 1:2 w  lp lc rgb 'orange'  lw 3 title ' \footnotesize Error in $J_4$', \
			'< sqlite3 Data/Ybeam/ysplitterConfig3.db "SELECT DISTINCT DOFS_primal, abs(relativeError4) from data WHERE DOFS_primal>100"' u 1:2 w  lp  lw 3 title ' \footnotesize Error in $J_5$', \
			'< sqlite3 Data/Ybeam/ysplitterConfig3.db "SELECT DISTINCT DOFS_primal, abs(relativeError5) from data WHERE DOFS_primal>100"' u 1:2 w  lp  lw 3 title ' \footnotesize Error in $J_6$', \
			'< sqlite3 Data/Ybeam/ysplitterConfig3.db "SELECT DISTINCT DOFS_primal, abs(relativeError6) from data WHERE DOFS_primal>100"' u 1:2 w  lp  lw 3 title ' \footnotesize Error in $J_7$', \
			10/x**1.0 lw  10	 title ' \footnotesize $\text{DOFs}^{-1}$'										
			#					 '< sqlite3 Data/Multigoalp4/Higher_Order/dataHigherOrderJE.db "SELECT DISTINCT DOFs, abs(Exact_Error) from data "' u 1:2 w  lp lw 3 title ' \footnotesize Error in $J_\mathfrak{E}$', \
		\end{gnuplot}
		\fi
		{		\scalebox{0.55}{\input{Figures/YSplitterErrorConfig4.tex}}}
		\caption{The errors and error estimator for Configuration 4.}
		\label{fig_test_3_d}
	\end{minipage}
	
\end{figure}

\begin{figure}[H]
	\begin{minipage}{0.49\linewidth}
		\ifMAKEPICS
		\begin{gnuplot}[terminal=epslatex,terminaloptions=color]
			set output "Figures/YSplitterErrorConfig5.tex"
			set title "Configuration 5"
			set rmargin 15
			set key outside
			set key right top
			set key opaque
			set datafile separator "|"
			set logscale x
			set logscale y
			set grid ytics lc rgb "#bbbbbb" lw 1 lt 0
			set grid xtics lc rgb "#bbbbbb" lw 1 lt 0
			set xlabel '\text{DOFs}'
			set format '
			plot \
			'< sqlite3 Data/Ybeam/ysplitterConfig4.db "SELECT DISTINCT DOFS_primal, abs(Exact_Error) from data WHERE DOFS_primal>100"' u 1:2 w  lp lt 1 lw 3 title ' \footnotesize Error in $J_c$', \
			'< sqlite3 Data/Ybeam/ysplitterConfig4.db "SELECT DISTINCT DOFS_primal, abs(Estimated_Error) from data WHERE DOFS_primal>100"' u 1:2 w  lp  lw 3 title ' \footnotesize  $\eta_h$', \
			'< sqlite3 Data/Ybeam/ysplitterConfig4.db "SELECT DISTINCT DOFS_primal, abs(relativeError0) from data WHERE DOFS_primal>100"' u 1:2 w  lp  lw 3 title ' \footnotesize Error in $J_1$', \
			'< sqlite3 Data/Ybeam/ysplitterConfig4.db "SELECT DISTINCT DOFS_primal, abs(relativeError1) from data WHERE DOFS_primal>100"' u 1:2 w  lp  lw 3 title ' \footnotesize Error in $J_2$', \
			'< sqlite3 Data/Ybeam/ysplitterConfig4.db "SELECT DISTINCT DOFS_primal, abs(relativeError2) from data WHERE DOFS_primal>100"' u 1:2 w  lp  lw 3 title ' \footnotesize Error in $J_3$', \
			'< sqlite3 Data/Ybeam/ysplitterConfig4.db "SELECT DISTINCT DOFS_primal, abs(relativeError3) from data WHERE DOFS_primal>100"' u 1:2 w  lp lc rgb 'orange'  lw 3 title ' \footnotesize Error in $J_4$', \
			'< sqlite3 Data/Ybeam/ysplitterConfig4.db "SELECT DISTINCT DOFS_primal, abs(relativeError4) from data WHERE DOFS_primal>100"' u 1:2 w  lp  lw 3 title ' \footnotesize Error in $J_5$', \
			'< sqlite3 Data/Ybeam/ysplitterConfig4.db "SELECT DISTINCT DOFS_primal, abs(relativeError5) from data WHERE DOFS_primal>100"' u 1:2 w  lp  lw 3 title ' \footnotesize Error in $J_6$', \
			'< sqlite3 Data/Ybeam/ysplitterConfig4.db "SELECT DISTINCT DOFS_primal, abs(relativeError6) from data WHERE DOFS_primal>100"' u 1:2 w  lp  lw 3 title ' \footnotesize Error in $J_7$', \
			100/x**1.5 lw  10  title '$\mathcal{O}(\text{DOFs}^-\frac{3}{2})$'											
			#					 '< sqlite3 Data/Multigoalp4/Higher_Order/dataHigherOrderJE.db "SELECT DISTINCT DOFs, abs(Exact_Error) from data "' u 1:2 w  lp lw 3 title ' \footnotesize Error in $J_\mathfrak{E}$', \
		\end{gnuplot}
		\fi
		{		\scalebox{0.55}{\input{Figures/YSplitterErrorConfig5.tex}}}
		\caption{The errors and error estimator for Configuration 5.}
	\end{minipage}%
	\hfill
	\begin{minipage}{0.49\linewidth}
		\ifMAKEPICS
		\begin{gnuplot}[terminal=epslatex,terminaloptions=color]
			set output "Figures/YSplitterErrorConfig6.tex"
			set title "Configuration 6"
			set rmargin 15
			set key outside
			set key right top
			set key opaque
			set datafile separator "|"
			set logscale x
			set logscale y
			set grid ytics lc rgb "#bbbbbb" lw 1 lt 0
			set grid xtics lc rgb "#bbbbbb" lw 1 lt 0
			set xlabel '\text{DOFs}'
			set format '
			plot \
			'< sqlite3 Data/Ybeam/ysplitterConfig6.db "SELECT DISTINCT DOFS_primal, abs(Exact_Error) from data WHERE DOFS_primal>100"' u 1:2 w  lp lt 1 lw 3 title ' \footnotesize Error in $J_c$', \
			'< sqlite3 Data/Ybeam/ysplitterConfig6.db "SELECT DISTINCT DOFS_primal, abs(Estimated_Error) from data WHERE DOFS_primal>100"' u 1:2 w  lp  lw 3 title ' \footnotesize  $\eta_h$', \
			'< sqlite3 Data/Ybeam/ysplitterConfig6.db "SELECT DISTINCT DOFS_primal, abs(relativeError0) from data WHERE DOFS_primal>100"' u 1:2 w  lp  lw 3 title ' \footnotesize Error in $J_1$', \
			'< sqlite3 Data/Ybeam/ysplitterConfig6.db "SELECT DISTINCT DOFS_primal, abs(relativeError1) from data WHERE DOFS_primal>100"' u 1:2 w  lp  lw 3 title ' \footnotesize Error in $J_2$', \
			'< sqlite3 Data/Ybeam/ysplitterConfig6.db "SELECT DISTINCT DOFS_primal, abs(relativeError2) from data WHERE DOFS_primal>100"' u 1:2 w  lp  lw 3 title ' \footnotesize Error in $J_3$', \
			'< sqlite3 Data/Ybeam/ysplitterConfig6.db "SELECT DISTINCT DOFS_primal, abs(relativeError3) from data WHERE DOFS_primal>100"' u 1:2 w  lp lc rgb 'orange'  lw 3 title ' \footnotesize Error in $J_4$', \
			'< sqlite3 Data/Ybeam/ysplitterConfig6.db "SELECT DISTINCT DOFS_primal, abs(relativeError4) from data WHERE DOFS_primal>100"' u 1:2 w  lp  lw 3 title ' \footnotesize Error in $J_5$', \
			'< sqlite3 Data/Ybeam/ysplitterConfig6.db "SELECT DISTINCT DOFS_primal, abs(relativeError5) from data WHERE DOFS_primal>100"' u 1:2 w  lp  lw 3 title ' \footnotesize Error in $J_6$', \
			'< sqlite3 Data/Ybeam/ysplitterConfig6.db "SELECT DISTINCT DOFS_primal, abs(relativeError6) from data WHERE DOFS_primal>100"' u 1:2 w  lp  lw 3 title ' \footnotesize Error in $J_7$', \
			100/x**1.5 lw  10 title '$\mathcal{O}(\text{DOFs}^-\frac{3}{2})$'										
			#					 '< sqlite3 Data/Multigoalp4/Higher_Order/dataHigherOrderJE.db "SELECT DISTINCT DOFs, abs(Exact_Error) from data "' u 1:2 w  lp lw 3 title ' \footnotesize Error in $J_\mathfrak{E}$', \
		\end{gnuplot}
		\fi
		{		\scalebox{0.55}{\input{Figures/YSplitterErrorConfig6.tex}}}
		
		\caption{The errors and error estimator for Configuration 6.}
		\label{fig_test_3_e}
	\end{minipage}
\end{figure}

\begin{figure}[H]
	\begin{minipage}{0.49\linewidth}
		\ifMAKEPICS
		\begin{gnuplot}[terminal=epslatex,terminaloptions=color]
			set output "Figures/YSplitterIeffConfig1.tex"
			set title "Configuration 1"
			set rmargin 15
			set key outside
			set key right top
			set key opaque
			set datafile separator "|"
			set logscale x
			set logscale y
			set grid ytics lc rgb "#bbbbbb" lw 1 lt 0
			set grid xtics lc rgb "#bbbbbb" lw 1 lt 0
			set xlabel '\text{DOFs}'
			set format '
			plot \
			'< sqlite3 Data/Ybeam/ysplitterConfig1.db "SELECT DISTINCT DOFS_primal, Ieff from data WHERE DOFS_primal>100"' u 1:2 w  lp lt 1 lw 3 title ' \footnotesize $I_{eff}$', \
			'< sqlite3 Data/Ybeam/ysplitterConfig1.db "SELECT DISTINCT DOFS_primal, Ieff_adjoint from data WHERE DOFS_primal>100"' u 1:2 w  lp  lw 3 title ' \footnotesize  $I_{eff,a}$', \
			'< sqlite3 Data/Ybeam/ysplitterConfig1.db "SELECT DISTINCT DOFS_primal, Ieff_primal from data WHERE DOFS_primal>100"' u 1:2 w  lp  lw 3 title ' \footnotesize $I_{eff,p}$', \
			1
		\end{gnuplot}
		\fi
		{		\scalebox{0.55}{\input{Figures/YSplitterIeffConfig1.tex}}}
		\caption{Effectivity indices for Configuration 1.}
		\label{fig_test_3_fa}
	\end{minipage}%
	\hfill
	\begin{minipage}{0.49\linewidth}
		\ifMAKEPICS
		\begin{gnuplot}[terminal=epslatex,terminaloptions=color]
			set output "Figures/YSplitterIeffConfig2.tex"
			set title "Configuration 2"
			set rmargin 15
			set key outside
			set key right top
			set key opaque
			set datafile separator "|"
			set logscale x
			set logscale y
			set grid ytics lc rgb "#bbbbbb" lw 1 lt 0
			set grid xtics lc rgb "#bbbbbb" lw 1 lt 0
			set xlabel '\text{DOFs}'
			set format '
			plot \
			'< sqlite3 Data/Ybeam/ysplitterConfig5.db "SELECT DISTINCT DOFS_primal, Ieff from data WHERE DOFS_primal>100"' u 1:2 w  lp lt 1 lw 3 title ' \footnotesize $I_{eff}$', \
			'< sqlite3 Data/Ybeam/ysplitterConfig5.db "SELECT DISTINCT DOFS_primal, Ieff_adjoint from data WHERE DOFS_primal>100"' u 1:2 w  lp  lw 3 title ' \footnotesize  $I_{eff,a}$', \
			'< sqlite3 Data/Ybeam/ysplitterConfig5.db "SELECT DISTINCT DOFS_primal, Ieff_primal from data WHERE DOFS_primal>100"' u 1:2 w  lp  lw 3 title ' \footnotesize $I_{eff,p}$', \
			1
		\end{gnuplot}
		\fi
		{		\scalebox{0.55}{\input{Figures/YSplitterIeffConfig2.tex}}}
		\caption{Effectivity indices for Configuration 2.}
		\label{fig_test_3_f}
	\end{minipage}
\end{figure}

\begin{figure}[H]
	\begin{minipage}{0.49\linewidth}
		\ifMAKEPICS
		\begin{gnuplot}[terminal=epslatex,terminaloptions=color]
			set output "Figures/YSplitterIeffConfig3.tex"
			set title "Configuration 3"
			set rmargin 15
			set key outside
			set key right top
			set key opaque
			set datafile separator "|"
			set logscale x
			set logscale y
			set grid ytics lc rgb "#bbbbbb" lw 1 lt 0
			set grid xtics lc rgb "#bbbbbb" lw 1 lt 0
			set xlabel '\text{DOFs}'
			set format '
			plot \
			'< sqlite3 Data/Ybeam/ysplitterConfig2.db "SELECT DISTINCT DOFS_primal, Ieff from data WHERE DOFS_primal>100"' u 1:2 w  lp lt 1 lw 3 title ' \footnotesize $I_{eff}$', \
			'< sqlite3 Data/Ybeam/ysplitterConfig2.db "SELECT DISTINCT DOFS_primal, Ieff_adjoint from data WHERE DOFS_primal>100"' u 1:2 w  lp  lw 3 title ' \footnotesize  $I_{eff,a}$', \
			'< sqlite3 Data/Ybeam/ysplitterConfig2.db "SELECT DISTINCT DOFS_primal, Ieff_primal from data WHERE DOFS_primal>100"' u 1:2 w  lp  lw 3 title ' \footnotesize $I_{eff,p}$', \
			1
		\end{gnuplot}
		\fi
		{		\scalebox{0.55}{\input{Figures/YSplitterIeffConfig3.tex}}}
		\caption{Effectivity indices for Configuration 3.}
	\end{minipage}%
	\hfill
	\begin{minipage}{0.49\linewidth}
		\ifMAKEPICS
		\begin{gnuplot}[terminal=epslatex,terminaloptions=color]
			set output "Figures/YSplitterIeffConfig4.tex"
			set title "Configuration 4"
			set rmargin 15
			set key outside
			set key right top
			set key opaque
			set datafile separator "|"
			set logscale x
			set logscale y
			set grid ytics lc rgb "#bbbbbb" lw 1 lt 0
			set grid xtics lc rgb "#bbbbbb" lw 1 lt 0
			set xlabel '\text{DOFs}'
			set format '
			plot \
			'< sqlite3 Data/Ybeam/ysplitterConfig3.db "SELECT DISTINCT DOFS_primal, Ieff from data WHERE DOFS_primal>100"' u 1:2 w  lp lt 1 lw 3 title ' \footnotesize $I_{eff}$', \
			'< sqlite3 Data/Ybeam/ysplitterConfig3.db "SELECT DISTINCT DOFS_primal, Ieff_adjoint from data WHERE DOFS_primal>100"' u 1:2 w  lp  lw 3 title ' \footnotesize  $I_{eff,a}$', \
			'< sqlite3 Data/Ybeam/ysplitterConfig3.db "SELECT DISTINCT DOFS_primal, Ieff_primal from data WHERE DOFS_primal>100"' u 1:2 w  lp  lw 3 title ' \footnotesize $I_{eff,p}$', \
			1
		\end{gnuplot}
		\fi
		{		\scalebox{0.55}{\input{Figures/YSplitterIeffConfig4.tex}}}
		\caption{Effectivity indices for Configuration 4.}
		\label{fig_test_3_g}
	\end{minipage}
\end{figure}

\begin{figure}[H]
	\begin{minipage}{0.49\linewidth}
		\ifMAKEPICS
		\begin{gnuplot}[terminal=epslatex,terminaloptions=color]
			set output "Figures/YSplitterIeffConfig5.tex"
			set title "Configuration 5"
			set rmargin 15
			set key outside
			set key right top
			set key opaque
			set datafile separator "|"
			set logscale x
			set logscale y
			set grid ytics lc rgb "#bbbbbb" lw 1 lt 0
			set grid xtics lc rgb "#bbbbbb" lw 1 lt 0
			set xlabel '\text{DOFs}'
			set format '
			plot \
			'< sqlite3 Data/Ybeam/ysplitterConfig4.db "SELECT DISTINCT DOFS_primal, Ieff from data WHERE DOFS_primal>100"' u 1:2 w  lp lt 1 lw 3 title ' \footnotesize $I_{eff}$', \
			'< sqlite3 Data/Ybeam/ysplitterConfig4.db "SELECT DISTINCT DOFS_primal, Ieff_adjoint from data WHERE DOFS_primal>100"' u 1:2 w  lp  lw 3 title ' \footnotesize  $I_{eff,a}$', \
			'< sqlite3 Data/Ybeam/ysplitterConfig4.db "SELECT DISTINCT DOFS_primal, Ieff_primal from data WHERE DOFS_primal>100"' u 1:2 w  lp  lw 3 title ' \footnotesize $I_{eff,p}$', \
			1
		\end{gnuplot}
		\fi
		{		\scalebox{0.55}{\input{Figures/YSplitterIeffConfig5.tex}}}
		\caption{Effectivity indices for Configuration 5.}
	\end{minipage}%
	\hfill
	\begin{minipage}{0.49\linewidth}
		\ifMAKEPICS
		\begin{gnuplot}[terminal=epslatex,terminaloptions=color]
			set output "Figures/YSplitterIeffConfig6.tex"
			set title "Configuration 6"
			set rmargin 15
			set key outside
			set key right top
			set key opaque
			set datafile separator "|"
			set logscale x
			set logscale y
			set grid ytics lc rgb "#bbbbbb" lw 1 lt 0
			set grid xtics lc rgb "#bbbbbb" lw 1 lt 0
			set xlabel '\text{DOFs}'
			set format '
			plot \
			'< sqlite3 Data/Ybeam/ysplitterConfig6.db "SELECT DISTINCT DOFS_primal, Ieff from data WHERE DOFS_primal>100"' u 1:2 w  lp lt 1 lw 3 title ' \footnotesize $I_{eff}$', \
			'< sqlite3 Data/Ybeam/ysplitterConfig6.db "SELECT DISTINCT DOFS_primal, Ieff_adjoint from data WHERE DOFS_primal>100"' u 1:2 w  lp  lw 3 title ' \footnotesize  $I_{eff,a}$', \
			'< sqlite3 Data/Ybeam/ysplitterConfig6.db "SELECT DISTINCT DOFS_primal, Ieff_primal from data WHERE DOFS_primal>100"' u 1:2 w  lp  lw 3 title ' \footnotesize $I_{eff,p}$', \
			1
		\end{gnuplot}
		\fi
		{		\scalebox{0.55}{\input{Figures/YSplitterIeffConfig6.tex}}}
		
		\caption{Effectivity indices for Configuration 6}
		\label{fig_test_3_h}
	\end{minipage}
\end{figure}

	\section{Conclusions}
	\label{sec_conclusions}
	In this work, we modeled laser material processing with the help
	of a generalized Boussinesq model. The resulting PDE system is nonlinear and we considered
	a monolithic coupling scheme. The focus was on multi-goal a posteriori 
	error estimation and local mesh adaptivity. Since the Boussinesq system 
	consists of two coupled PDEs (i.e., incompressible Navier-Stokes coupled to 
	a stationary heat equation ) several quantities of interest (i.e., goal functionals) might be 
	of interest simultaneously. In our multigoal-framework, a combined goal 
	functional is defined and serves as right hand side in the adjoint problem from which
	local sensitivity measures enter the error estimator. Three numerical experiments 
	were conducted: One classical benchmark problem and two configurations that are 
	motivated from interdisciplinary collaborations in our Cluster of Excellence. 
	In all tests, we observed error reductions and effectivity indices. The latter 
	show excellent performance and indicate that we have a robust 
	and efficient adaptive framework at hand. In future work, we plan 
	\black{to extend this framework to three-dimensional situations.
	Theorem 2.2 and the multigoal framework in general cover three-dimensional domains 
	(see also \cite{EndtLaWi18_pamm,Endt21}). A suitable analogue for the Y-beam splitter of 
	Section \ref{sec_Y_beam} would be $\Omega\times I$, with $I=(0,0.2)$, $\Gamma_{\nu N}=(\Gamma_1\cup \Gamma_2\cup \Gamma_3)\times I$ etc. (where $\Omega$, $\Gamma_1$, $\Gamma_2$, $\Gamma_3$ are as in section 4.3), and the points $A,\ldots,F$ obtaining an additional third component of $0.1$ each; a choice of $g=(0,0,-9.81)^T$ in place of $g=(0,-9.81)^T$ would be interesting. The form of the goal functionals of Section \ref{sec_Y_beam_goal} remains unchanged.
	However, the computational extension requires some work, 
	with the main bottleneck being the linear solver and preconditioners that need to be developed 
	due to memory consumptions and computational cost. 
	Finally,} another future extension are time-dependent cases by using a space-time framework.

	\section*{Acknowledgments}
	This work has been supported by 
	the Cluster of Excellence PhoenixD (EXC 2122, Project ID 390833453).
	Furthermore, the second author is funded by an Humboldt Postdoctoral 
	Fellowship. 
	
	
	\bibliographystyle{abbrv}
	\bibliography{lit}

\end{document}

%% file: header.tex
\usepackage{amsmath,amsthm,amsfonts,amssymb,amscd}
\usepackage{mathtools}
\usepackage{pgf,tikz}
\usetikzlibrary{arrows}
\usepackage{mathrsfs}
\usepackage{float}
\usepackage{csvsimple}
\usepackage{enumerate}
\usepackage{booktabs}
\usepackage{pgfplots}
\usepackage[noend]{algpseudocode}
\usepackage{graphicx,bm,xcolor}
\usepackage{algorithm}
\usepackage{pdfpages}
\usepackage{algpseudocode}
\usepackage{stmaryrd}
\usepackage[square,sort,comma,numbers]{natbib}
\usepackage{hyperref}
\usetikzlibrary{arrows}
\usepackage{caption}

\usepackage{t1enc}
\usepackage[german,english]{babel}
\usepackage{verbatim} 

\usepackage{url}

\newif\ifMAKEPICS
\MAKEPICSfalse

\ifMAKEPICS
\usepackage[cleanup,subfolder]{gnuplottex}
\usepackage{xparse}

\ExplSyntaxOn
\DeclareExpandableDocumentCommand{\convertlen}{ O{cm} m }
{
	\dim_to_decimal_in_unit:nn { #2 } { 1 #1 } cm
}
\ExplSyntaxOff
\fi

%% file: Figures/Examplecoldtowarm.tex
\begingroup
  \makeatletter
  \providecommand\color[2][]{%
    \GenericError{(gnuplot) \space\space\space\@spaces}{%
      Package color not loaded in conjunction with
      terminal option `colourtext'%
    }{See the gnuplot documentation for explanation.%
    }{Either use 'blacktext' in gnuplot or load the package
      color.sty in LaTeX.}%
    \renewcommand\color[2][]{}%
  }%
  \providecommand\includegraphics[2][]{%
    \GenericError{(gnuplot) \space\space\space\@spaces}{%
      Package graphicx or graphics not loaded%
    }{See the gnuplot documentation for explanation.%
    }{The gnuplot epslatex terminal needs graphicx.sty or graphics.sty.}%
    \renewcommand\includegraphics[2][]{}%
  }%
  \providecommand\rotatebox[2]{#2}%
  \@ifundefined{ifGPcolor}{%
    \newif\ifGPcolor
    \GPcolortrue
  }{}%
  \@ifundefined{ifGPblacktext}{%
    \newif\ifGPblacktext
    \GPblacktexttrue
  }{}%
  \let\gplgaddtomacro\g@addto@macro
  \gdef\gplbacktext{}%
  \gdef\gplfronttext{}%
  \makeatother
  \ifGPblacktext
    \def\colorrgb#1{}%
    \def\colorgray#1{}%
  \else
    \ifGPcolor
      \def\colorrgb#1{\color[rgb]{#1}}%
      \def\colorgray#1{\color[gray]{#1}}%
      \expandafter\def\csname LTw\endcsname{\color{white}}%
      \expandafter\def\csname LTb\endcsname{\color{black}}%
      \expandafter\def\csname LTa\endcsname{\color{black}}%
      \expandafter\def\csname LT0\endcsname{\color[rgb]{1,0,0}}%
      \expandafter\def\csname LT1\endcsname{\color[rgb]{0,1,0}}%
      \expandafter\def\csname LT2\endcsname{\color[rgb]{0,0,1}}%
      \expandafter\def\csname LT3\endcsname{\color[rgb]{1,0,1}}%
      \expandafter\def\csname LT4\endcsname{\color[rgb]{0,1,1}}%
      \expandafter\def\csname LT5\endcsname{\color[rgb]{1,1,0}}%
      \expandafter\def\csname LT6\endcsname{\color[rgb]{0,0,0}}%
      \expandafter\def\csname LT7\endcsname{\color[rgb]{1,0.3,0}}%
      \expandafter\def\csname LT8\endcsname{\color[rgb]{0.5,0.5,0.5}}%
    \else
      \def\colorrgb#1{\color{black}}%
      \def\colorgray#1{\color[gray]{#1}}%
      \expandafter\def\csname LTw\endcsname{\color{white}}%
      \expandafter\def\csname LTb\endcsname{\color{black}}%
      \expandafter\def\csname LTa\endcsname{\color{black}}%
      \expandafter\def\csname LT0\endcsname{\color{black}}%
      \expandafter\def\csname LT1\endcsname{\color{black}}%
      \expandafter\def\csname LT2\endcsname{\color{black}}%
      \expandafter\def\csname LT3\endcsname{\color{black}}%
      \expandafter\def\csname LT4\endcsname{\color{black}}%
      \expandafter\def\csname LT5\endcsname{\color{black}}%
      \expandafter\def\csname LT6\endcsname{\color{black}}%
      \expandafter\def\csname LT7\endcsname{\color{black}}%
      \expandafter\def\csname LT8\endcsname{\color{black}}%
    \fi
  \fi
  \setlength{\unitlength}{0.0500bp}%
  \begin{picture}(7200.00,5040.00)%
    \gplgaddtomacro\gplbacktext{%
      \csname LTb\endcsname%
      \put(990,704){\makebox(0,0)[r]{\strut{}0.0001}}%
      \csname LTb\endcsname%
      \put(990,1439){\makebox(0,0)[r]{\strut{}0.001}}%
      \csname LTb\endcsname%
      \put(990,2174){\makebox(0,0)[r]{\strut{}0.01}}%
      \csname LTb\endcsname%
      \put(990,2909){\makebox(0,0)[r]{\strut{}0.1}}%
      \csname LTb\endcsname%
      \put(990,3644){\makebox(0,0)[r]{\strut{}1}}%
      \csname LTb\endcsname%
      \put(990,4379){\makebox(0,0)[r]{\strut{}10}}%
      \csname LTb\endcsname%
      \put(1122,484){\makebox(0,0){\strut{}1000}}%
      \csname LTb\endcsname%
      \put(3016,484){\makebox(0,0){\strut{}10000}}%
      \csname LTb\endcsname%
      \put(4909,484){\makebox(0,0){\strut{}100000}}%
      \csname LTb\endcsname%
      \put(6803,484){\makebox(0,0){\strut{}1e+06}}%
      \put(3962,154){\makebox(0,0){\strut{}\text{DOFs}}}%
      \put(3962,4709){\makebox(0,0){\strut{}cold to warm}}%
    }%
    \gplgaddtomacro\gplfronttext{%
      \csname LTb\endcsname%
      \put(4950,4206){\makebox(0,0)[r]{\strut{} \footnotesize Error in $p_{\text{diff}}$}}%
      \csname LTb\endcsname%
      \put(4950,3986){\makebox(0,0)[r]{\strut{} \footnotesize Error Estimator of $p_{\text{diff}}$}}%
      \csname LTb\endcsname%
      \put(4950,3766){\makebox(0,0)[r]{\strut{}$\mathcal{O}(\text{DOFs}^-1)$}}%
    }%
    \gplbacktext
    \put(0,0){\includegraphics{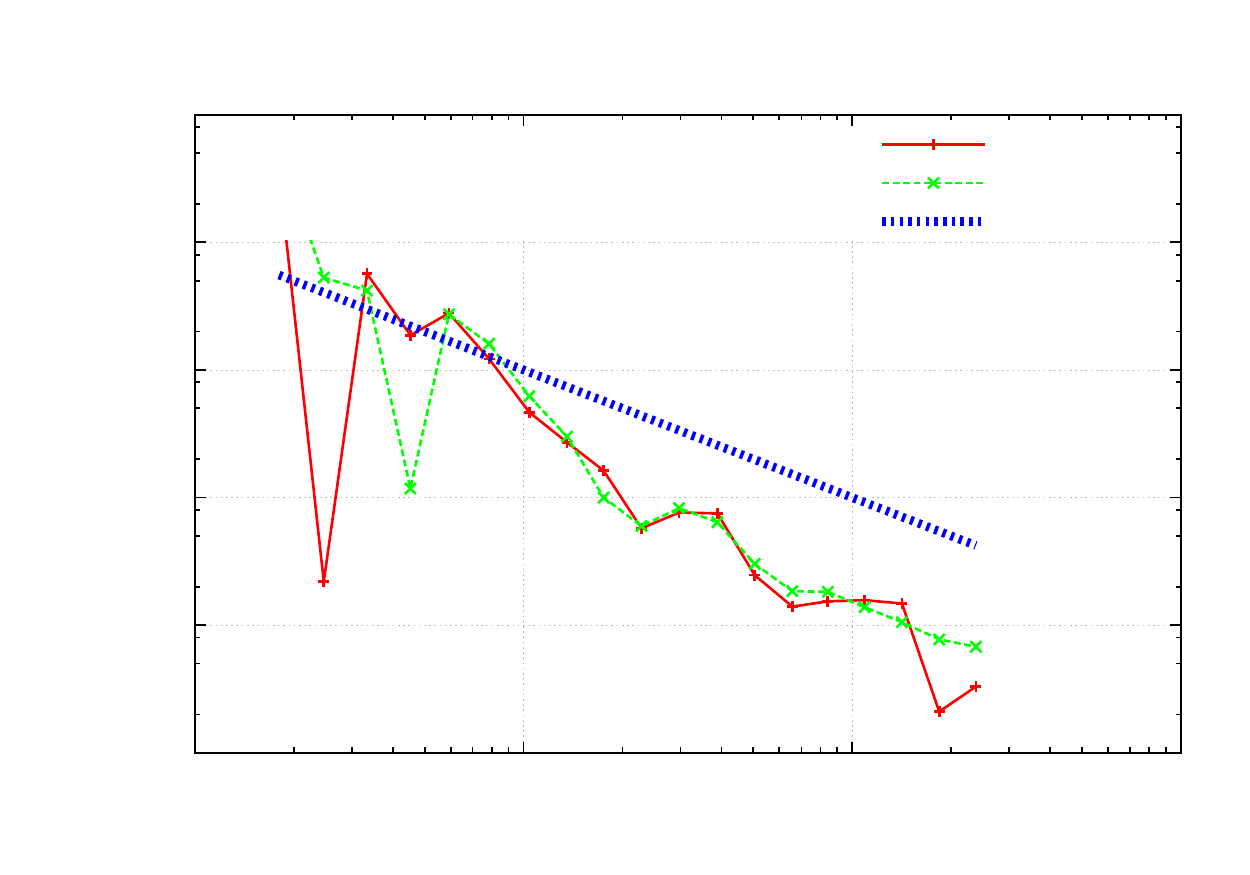}}%
    \gplfronttext
  \end{picture}%
\endgroup

%% file: Figures/Examplewarmtocold.tex
\begingroup
  \makeatletter
  \providecommand\color[2][]{%
    \GenericError{(gnuplot) \space\space\space\@spaces}{%
      Package color not loaded in conjunction with
      terminal option `colourtext'%
    }{See the gnuplot documentation for explanation.%
    }{Either use 'blacktext' in gnuplot or load the package
      color.sty in LaTeX.}%
    \renewcommand\color[2][]{}%
  }%
  \providecommand\includegraphics[2][]{%
    \GenericError{(gnuplot) \space\space\space\@spaces}{%
      Package graphicx or graphics not loaded%
    }{See the gnuplot documentation for explanation.%
    }{The gnuplot epslatex terminal needs graphicx.sty or graphics.sty.}%
    \renewcommand\includegraphics[2][]{}%
  }%
  \providecommand\rotatebox[2]{#2}%
  \@ifundefined{ifGPcolor}{%
    \newif\ifGPcolor
    \GPcolortrue
  }{}%
  \@ifundefined{ifGPblacktext}{%
    \newif\ifGPblacktext
    \GPblacktexttrue
  }{}%
  \let\gplgaddtomacro\g@addto@macro
  \gdef\gplbacktext{}%
  \gdef\gplfronttext{}%
  \makeatother
  \ifGPblacktext
    \def\colorrgb#1{}%
    \def\colorgray#1{}%
  \else
    \ifGPcolor
      \def\colorrgb#1{\color[rgb]{#1}}%
      \def\colorgray#1{\color[gray]{#1}}%
      \expandafter\def\csname LTw\endcsname{\color{white}}%
      \expandafter\def\csname LTb\endcsname{\color{black}}%
      \expandafter\def\csname LTa\endcsname{\color{black}}%
      \expandafter\def\csname LT0\endcsname{\color[rgb]{1,0,0}}%
      \expandafter\def\csname LT1\endcsname{\color[rgb]{0,1,0}}%
      \expandafter\def\csname LT2\endcsname{\color[rgb]{0,0,1}}%
      \expandafter\def\csname LT3\endcsname{\color[rgb]{1,0,1}}%
      \expandafter\def\csname LT4\endcsname{\color[rgb]{0,1,1}}%
      \expandafter\def\csname LT5\endcsname{\color[rgb]{1,1,0}}%
      \expandafter\def\csname LT6\endcsname{\color[rgb]{0,0,0}}%
      \expandafter\def\csname LT7\endcsname{\color[rgb]{1,0.3,0}}%
      \expandafter\def\csname LT8\endcsname{\color[rgb]{0.5,0.5,0.5}}%
    \else
      \def\colorrgb#1{\color{black}}%
      \def\colorgray#1{\color[gray]{#1}}%
      \expandafter\def\csname LTw\endcsname{\color{white}}%
      \expandafter\def\csname LTb\endcsname{\color{black}}%
      \expandafter\def\csname LTa\endcsname{\color{black}}%
      \expandafter\def\csname LT0\endcsname{\color{black}}%
      \expandafter\def\csname LT1\endcsname{\color{black}}%
      \expandafter\def\csname LT2\endcsname{\color{black}}%
      \expandafter\def\csname LT3\endcsname{\color{black}}%
      \expandafter\def\csname LT4\endcsname{\color{black}}%
      \expandafter\def\csname LT5\endcsname{\color{black}}%
      \expandafter\def\csname LT6\endcsname{\color{black}}%
      \expandafter\def\csname LT7\endcsname{\color{black}}%
      \expandafter\def\csname LT8\endcsname{\color{black}}%
    \fi
  \fi
  \setlength{\unitlength}{0.0500bp}%
  \begin{picture}(7200.00,5040.00)%
    \gplgaddtomacro\gplbacktext{%
      \csname LTb\endcsname%
      \put(990,704){\makebox(0,0)[r]{\strut{}0.0001}}%
      \csname LTb\endcsname%
      \put(990,1439){\makebox(0,0)[r]{\strut{}0.001}}%
      \csname LTb\endcsname%
      \put(990,2174){\makebox(0,0)[r]{\strut{}0.01}}%
      \csname LTb\endcsname%
      \put(990,2909){\makebox(0,0)[r]{\strut{}0.1}}%
      \csname LTb\endcsname%
      \put(990,3644){\makebox(0,0)[r]{\strut{}1}}%
      \csname LTb\endcsname%
      \put(990,4379){\makebox(0,0)[r]{\strut{}10}}%
      \csname LTb\endcsname%
      \put(1122,484){\makebox(0,0){\strut{}1000}}%
      \csname LTb\endcsname%
      \put(3016,484){\makebox(0,0){\strut{}10000}}%
      \csname LTb\endcsname%
      \put(4909,484){\makebox(0,0){\strut{}100000}}%
      \csname LTb\endcsname%
      \put(6803,484){\makebox(0,0){\strut{}1e+06}}%
      \put(3962,154){\makebox(0,0){\strut{}\text{DOFs}}}%
      \put(3962,4709){\makebox(0,0){\strut{}warm to cold}}%
    }%
    \gplgaddtomacro\gplfronttext{%
      \csname LTb\endcsname%
      \put(4950,4206){\makebox(0,0)[r]{\strut{} \footnotesize Error in $p_{\text{diff}}$}}%
      \csname LTb\endcsname%
      \put(4950,3986){\makebox(0,0)[r]{\strut{} \footnotesize Error Estimator of $p_{\text{diff}}$}}%
      \csname LTb\endcsname%
      \put(4950,3766){\makebox(0,0)[r]{\strut{}$\mathcal{O}(\text{DOFs}^-1)$}}%
    }%
    \gplbacktext
    \put(0,0){\includegraphics{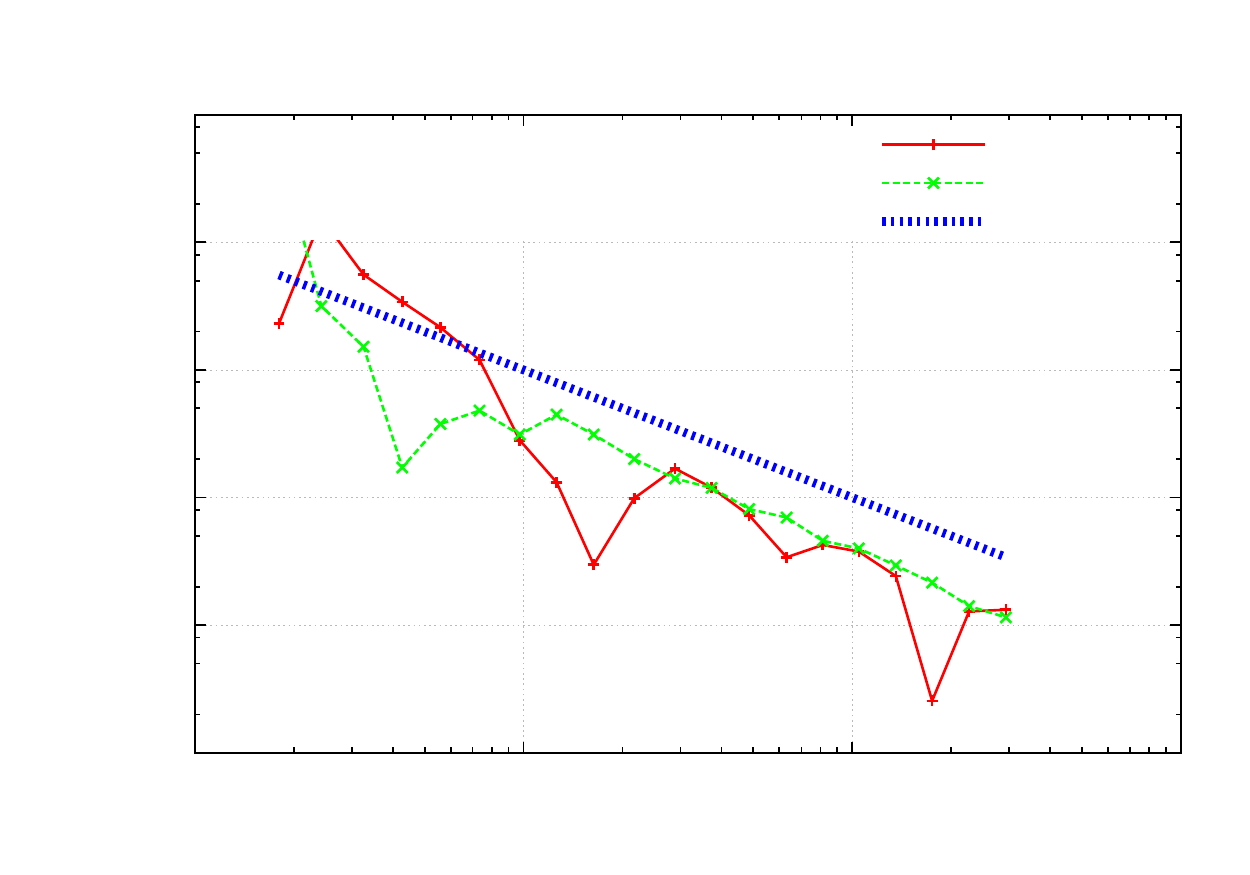}}%
    \gplfronttext
  \end{picture}%
\endgroup

%% file: Figures/LaseronSquareMeanError.tex
\begingroup
  \makeatletter
  \providecommand\color[2][]{%
    \GenericError{(gnuplot) \space\space\space\@spaces}{%
      Package color not loaded in conjunction with
      terminal option `colourtext'%
    }{See the gnuplot documentation for explanation.%
    }{Either use 'blacktext' in gnuplot or load the package
      color.sty in LaTeX.}%
    \renewcommand\color[2][]{}%
  }%
  \providecommand\includegraphics[2][]{%
    \GenericError{(gnuplot) \space\space\space\@spaces}{%
      Package graphicx or graphics not loaded%
    }{See the gnuplot documentation for explanation.%
    }{The gnuplot epslatex terminal needs graphicx.sty or graphics.sty.}%
    \renewcommand\includegraphics[2][]{}%
  }%
  \providecommand\rotatebox[2]{#2}%
  \@ifundefined{ifGPcolor}{%
    \newif\ifGPcolor
    \GPcolortrue
  }{}%
  \@ifundefined{ifGPblacktext}{%
    \newif\ifGPblacktext
    \GPblacktexttrue
  }{}%
  \let\gplgaddtomacro\g@addto@macro
  \gdef\gplbacktext{}%
  \gdef\gplfronttext{}%
  \makeatother
  \ifGPblacktext
    \def\colorrgb#1{}%
    \def\colorgray#1{}%
  \else
    \ifGPcolor
      \def\colorrgb#1{\color[rgb]{#1}}%
      \def\colorgray#1{\color[gray]{#1}}%
      \expandafter\def\csname LTw\endcsname{\color{white}}%
      \expandafter\def\csname LTb\endcsname{\color{black}}%
      \expandafter\def\csname LTa\endcsname{\color{black}}%
      \expandafter\def\csname LT0\endcsname{\color[rgb]{1,0,0}}%
      \expandafter\def\csname LT1\endcsname{\color[rgb]{0,1,0}}%
      \expandafter\def\csname LT2\endcsname{\color[rgb]{0,0,1}}%
      \expandafter\def\csname LT3\endcsname{\color[rgb]{1,0,1}}%
      \expandafter\def\csname LT4\endcsname{\color[rgb]{0,1,1}}%
      \expandafter\def\csname LT5\endcsname{\color[rgb]{1,1,0}}%
      \expandafter\def\csname LT6\endcsname{\color[rgb]{0,0,0}}%
      \expandafter\def\csname LT7\endcsname{\color[rgb]{1,0.3,0}}%
      \expandafter\def\csname LT8\endcsname{\color[rgb]{0.5,0.5,0.5}}%
    \else
      \def\colorrgb#1{\color{black}}%
      \def\colorgray#1{\color[gray]{#1}}%
      \expandafter\def\csname LTw\endcsname{\color{white}}%
      \expandafter\def\csname LTb\endcsname{\color{black}}%
      \expandafter\def\csname LTa\endcsname{\color{black}}%
      \expandafter\def\csname LT0\endcsname{\color{black}}%
      \expandafter\def\csname LT1\endcsname{\color{black}}%
      \expandafter\def\csname LT2\endcsname{\color{black}}%
      \expandafter\def\csname LT3\endcsname{\color{black}}%
      \expandafter\def\csname LT4\endcsname{\color{black}}%
      \expandafter\def\csname LT5\endcsname{\color{black}}%
      \expandafter\def\csname LT6\endcsname{\color{black}}%
      \expandafter\def\csname LT7\endcsname{\color{black}}%
      \expandafter\def\csname LT8\endcsname{\color{black}}%
    \fi
  \fi
  \setlength{\unitlength}{0.0500bp}%
  \begin{picture}(7200.00,5040.00)%
    \gplgaddtomacro\gplbacktext{%
      \csname LTb\endcsname%
      \put(990,704){\makebox(0,0)[r]{\strut{}1e-07}}%
      \csname LTb\endcsname%
      \put(990,1163){\makebox(0,0)[r]{\strut{}1e-06}}%
      \csname LTb\endcsname%
      \put(990,1623){\makebox(0,0)[r]{\strut{}1e-05}}%
      \csname LTb\endcsname%
      \put(990,2082){\makebox(0,0)[r]{\strut{}0.0001}}%
      \csname LTb\endcsname%
      \put(990,2542){\makebox(0,0)[r]{\strut{}0.001}}%
      \csname LTb\endcsname%
      \put(990,3001){\makebox(0,0)[r]{\strut{}0.01}}%
      \csname LTb\endcsname%
      \put(990,3460){\makebox(0,0)[r]{\strut{}0.1}}%
      \csname LTb\endcsname%
      \put(990,3920){\makebox(0,0)[r]{\strut{}1}}%
      \csname LTb\endcsname%
      \put(990,4379){\makebox(0,0)[r]{\strut{}10}}%
      \csname LTb\endcsname%
      \put(1122,484){\makebox(0,0){\strut{}100}}%
      \csname LTb\endcsname%
      \put(2542,484){\makebox(0,0){\strut{}1000}}%
      \csname LTb\endcsname%
      \put(3963,484){\makebox(0,0){\strut{}10000}}%
      \csname LTb\endcsname%
      \put(5383,484){\makebox(0,0){\strut{}100000}}%
      \csname LTb\endcsname%
      \put(6803,484){\makebox(0,0){\strut{}1e+06}}%
      \put(3962,154){\makebox(0,0){\strut{}\text{DOFs}}}%
      \put(3962,4709){\makebox(0,0){\strut{}Laser on Square}}%
    }%
    \gplgaddtomacro\gplfronttext{%
      \csname LTb\endcsname%
      \put(5816,4206){\makebox(0,0)[r]{\strut{} \footnotesize Error in $\bar{\T}$}}%
      \csname LTb\endcsname%
      \put(5816,3986){\makebox(0,0)[r]{\strut{} \footnotesize Error Estimator of $\bar{\T}$}}%
      \csname LTb\endcsname%
      \put(5816,3766){\makebox(0,0)[r]{\strut{}$\mathcal{O}(\text{DOFs}^-\frac{3}{2})$}}%
    }%
    \gplbacktext
    \put(0,0){\includegraphics{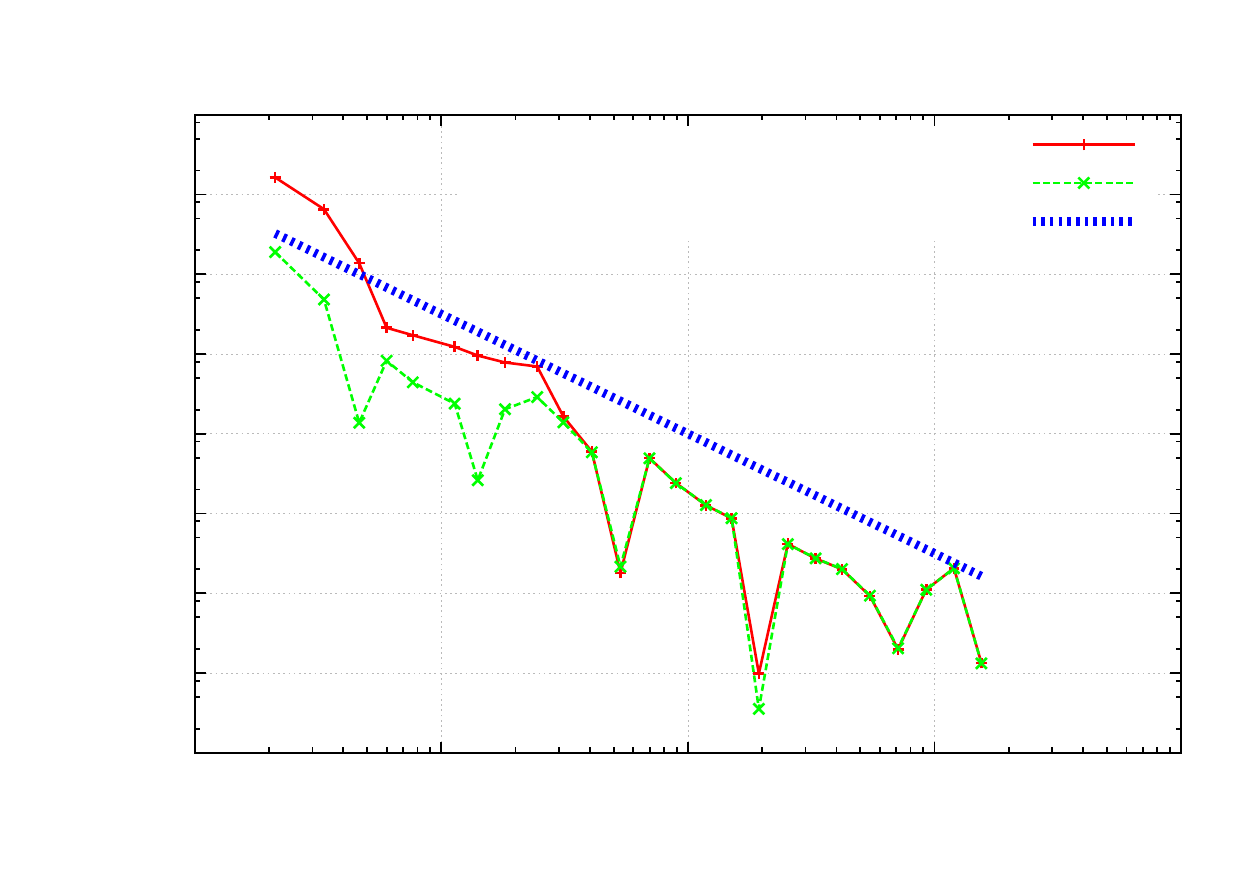}}%
    \gplfronttext
  \end{picture}%
\endgroup

%% file: Figures/LaseronSquareMeanIeff.tex
\begingroup
  \makeatletter
  \providecommand\color[2][]{%
    \GenericError{(gnuplot) \space\space\space\@spaces}{%
      Package color not loaded in conjunction with
      terminal option `colourtext'%
    }{See the gnuplot documentation for explanation.%
    }{Either use 'blacktext' in gnuplot or load the package
      color.sty in LaTeX.}%
    \renewcommand\color[2][]{}%
  }%
  \providecommand\includegraphics[2][]{%
    \GenericError{(gnuplot) \space\space\space\@spaces}{%
      Package graphicx or graphics not loaded%
    }{See the gnuplot documentation for explanation.%
    }{The gnuplot epslatex terminal needs graphicx.sty or graphics.sty.}%
    \renewcommand\includegraphics[2][]{}%
  }%
  \providecommand\rotatebox[2]{#2}%
  \@ifundefined{ifGPcolor}{%
    \newif\ifGPcolor
    \GPcolortrue
  }{}%
  \@ifundefined{ifGPblacktext}{%
    \newif\ifGPblacktext
    \GPblacktexttrue
  }{}%
  \let\gplgaddtomacro\g@addto@macro
  \gdef\gplbacktext{}%
  \gdef\gplfronttext{}%
  \makeatother
  \ifGPblacktext
    \def\colorrgb#1{}%
    \def\colorgray#1{}%
  \else
    \ifGPcolor
      \def\colorrgb#1{\color[rgb]{#1}}%
      \def\colorgray#1{\color[gray]{#1}}%
      \expandafter\def\csname LTw\endcsname{\color{white}}%
      \expandafter\def\csname LTb\endcsname{\color{black}}%
      \expandafter\def\csname LTa\endcsname{\color{black}}%
      \expandafter\def\csname LT0\endcsname{\color[rgb]{1,0,0}}%
      \expandafter\def\csname LT1\endcsname{\color[rgb]{0,1,0}}%
      \expandafter\def\csname LT2\endcsname{\color[rgb]{0,0,1}}%
      \expandafter\def\csname LT3\endcsname{\color[rgb]{1,0,1}}%
      \expandafter\def\csname LT4\endcsname{\color[rgb]{0,1,1}}%
      \expandafter\def\csname LT5\endcsname{\color[rgb]{1,1,0}}%
      \expandafter\def\csname LT6\endcsname{\color[rgb]{0,0,0}}%
      \expandafter\def\csname LT7\endcsname{\color[rgb]{1,0.3,0}}%
      \expandafter\def\csname LT8\endcsname{\color[rgb]{0.5,0.5,0.5}}%
    \else
      \def\colorrgb#1{\color{black}}%
      \def\colorgray#1{\color[gray]{#1}}%
      \expandafter\def\csname LTw\endcsname{\color{white}}%
      \expandafter\def\csname LTb\endcsname{\color{black}}%
      \expandafter\def\csname LTa\endcsname{\color{black}}%
      \expandafter\def\csname LT0\endcsname{\color{black}}%
      \expandafter\def\csname LT1\endcsname{\color{black}}%
      \expandafter\def\csname LT2\endcsname{\color{black}}%
      \expandafter\def\csname LT3\endcsname{\color{black}}%
      \expandafter\def\csname LT4\endcsname{\color{black}}%
      \expandafter\def\csname LT5\endcsname{\color{black}}%
      \expandafter\def\csname LT6\endcsname{\color{black}}%
      \expandafter\def\csname LT7\endcsname{\color{black}}%
      \expandafter\def\csname LT8\endcsname{\color{black}}%
    \fi
  \fi
  \setlength{\unitlength}{0.0500bp}%
  \begin{picture}(7200.00,5040.00)%
    \gplgaddtomacro\gplbacktext{%
      \csname LTb\endcsname%
      \put(726,704){\makebox(0,0)[r]{\strut{}0.01}}%
      \csname LTb\endcsname%
      \put(726,1929){\makebox(0,0)[r]{\strut{}0.1}}%
      \csname LTb\endcsname%
      \put(726,3154){\makebox(0,0)[r]{\strut{}1}}%
      \csname LTb\endcsname%
      \put(726,4379){\makebox(0,0)[r]{\strut{}10}}%
      \csname LTb\endcsname%
      \put(858,484){\makebox(0,0){\strut{}1000}}%
      \csname LTb\endcsname%
      \put(2840,484){\makebox(0,0){\strut{}10000}}%
      \csname LTb\endcsname%
      \put(4821,484){\makebox(0,0){\strut{}100000}}%
      \csname LTb\endcsname%
      \put(6803,484){\makebox(0,0){\strut{}1e+06}}%
      \put(3830,154){\makebox(0,0){\strut{}\text{DOFs}}}%
      \put(3830,4709){\makebox(0,0){\strut{}Laser on Square}}%
    }%
    \gplgaddtomacro\gplfronttext{%
      \csname LTb\endcsname%
      \put(2574,4206){\makebox(0,0)[r]{\strut{} \footnotesize $I_{eff}$ $\bar{\T}$}}%
      \csname LTb\endcsname%
      \put(2574,3986){\makebox(0,0)[r]{\strut{} \footnotesize $I_{eff,p}$ $\bar{\T}$}}%
      \csname LTb\endcsname%
      \put(2574,3766){\makebox(0,0)[r]{\strut{} \footnotesize $I_{eff,a}$ $\bar{\T}$}}%
      \csname LTb\endcsname%
      \put(2574,3546){\makebox(0,0)[r]{\strut{}1}}%
    }%
    \gplbacktext
    \put(0,0){\includegraphics{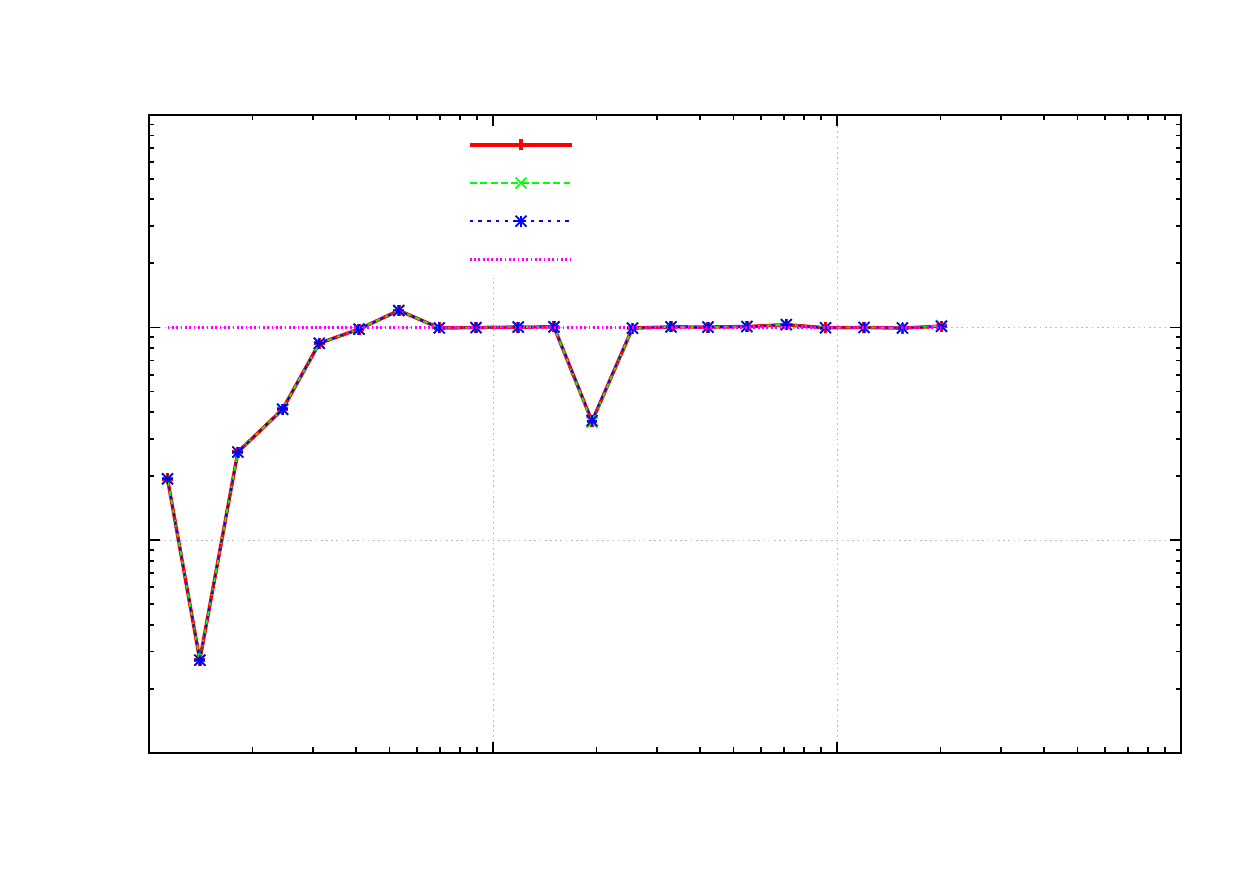}}%
    \gplfronttext
  \end{picture}%
\endgroup

%% file: Figures/LaseronSquareMeanVelociityError.tex
\begingroup
  \makeatletter
  \providecommand\color[2][]{%
    \GenericError{(gnuplot) \space\space\space\@spaces}{%
      Package color not loaded in conjunction with
      terminal option `colourtext'%
    }{See the gnuplot documentation for explanation.%
    }{Either use 'blacktext' in gnuplot or load the package
      color.sty in LaTeX.}%
    \renewcommand\color[2][]{}%
  }%
  \providecommand\includegraphics[2][]{%
    \GenericError{(gnuplot) \space\space\space\@spaces}{%
      Package graphicx or graphics not loaded%
    }{See the gnuplot documentation for explanation.%
    }{The gnuplot epslatex terminal needs graphicx.sty or graphics.sty.}%
    \renewcommand\includegraphics[2][]{}%
  }%
  \providecommand\rotatebox[2]{#2}%
  \@ifundefined{ifGPcolor}{%
    \newif\ifGPcolor
    \GPcolortrue
  }{}%
  \@ifundefined{ifGPblacktext}{%
    \newif\ifGPblacktext
    \GPblacktexttrue
  }{}%
  \let\gplgaddtomacro\g@addto@macro
  \gdef\gplbacktext{}%
  \gdef\gplfronttext{}%
  \makeatother
  \ifGPblacktext
    \def\colorrgb#1{}%
    \def\colorgray#1{}%
  \else
    \ifGPcolor
      \def\colorrgb#1{\color[rgb]{#1}}%
      \def\colorgray#1{\color[gray]{#1}}%
      \expandafter\def\csname LTw\endcsname{\color{white}}%
      \expandafter\def\csname LTb\endcsname{\color{black}}%
      \expandafter\def\csname LTa\endcsname{\color{black}}%
      \expandafter\def\csname LT0\endcsname{\color[rgb]{1,0,0}}%
      \expandafter\def\csname LT1\endcsname{\color[rgb]{0,1,0}}%
      \expandafter\def\csname LT2\endcsname{\color[rgb]{0,0,1}}%
      \expandafter\def\csname LT3\endcsname{\color[rgb]{1,0,1}}%
      \expandafter\def\csname LT4\endcsname{\color[rgb]{0,1,1}}%
      \expandafter\def\csname LT5\endcsname{\color[rgb]{1,1,0}}%
      \expandafter\def\csname LT6\endcsname{\color[rgb]{0,0,0}}%
      \expandafter\def\csname LT7\endcsname{\color[rgb]{1,0.3,0}}%
      \expandafter\def\csname LT8\endcsname{\color[rgb]{0.5,0.5,0.5}}%
    \else
      \def\colorrgb#1{\color{black}}%
      \def\colorgray#1{\color[gray]{#1}}%
      \expandafter\def\csname LTw\endcsname{\color{white}}%
      \expandafter\def\csname LTb\endcsname{\color{black}}%
      \expandafter\def\csname LTa\endcsname{\color{black}}%
      \expandafter\def\csname LT0\endcsname{\color{black}}%
      \expandafter\def\csname LT1\endcsname{\color{black}}%
      \expandafter\def\csname LT2\endcsname{\color{black}}%
      \expandafter\def\csname LT3\endcsname{\color{black}}%
      \expandafter\def\csname LT4\endcsname{\color{black}}%
      \expandafter\def\csname LT5\endcsname{\color{black}}%
      \expandafter\def\csname LT6\endcsname{\color{black}}%
      \expandafter\def\csname LT7\endcsname{\color{black}}%
      \expandafter\def\csname LT8\endcsname{\color{black}}%
    \fi
  \fi
  \setlength{\unitlength}{0.0500bp}%
  \begin{picture}(7200.00,5040.00)%
    \gplgaddtomacro\gplbacktext{%
      \csname LTb\endcsname%
      \put(990,704){\makebox(0,0)[r]{\strut{}1e-10}}%
      \csname LTb\endcsname%
      \put(990,1229){\makebox(0,0)[r]{\strut{}1e-09}}%
      \csname LTb\endcsname%
      \put(990,1754){\makebox(0,0)[r]{\strut{}1e-08}}%
      \csname LTb\endcsname%
      \put(990,2279){\makebox(0,0)[r]{\strut{}1e-07}}%
      \csname LTb\endcsname%
      \put(990,2804){\makebox(0,0)[r]{\strut{}1e-06}}%
      \csname LTb\endcsname%
      \put(990,3329){\makebox(0,0)[r]{\strut{}1e-05}}%
      \csname LTb\endcsname%
      \put(990,3854){\makebox(0,0)[r]{\strut{}0.0001}}%
      \csname LTb\endcsname%
      \put(990,4379){\makebox(0,0)[r]{\strut{}0.001}}%
      \csname LTb\endcsname%
      \put(1122,484){\makebox(0,0){\strut{}100}}%
      \csname LTb\endcsname%
      \put(2542,484){\makebox(0,0){\strut{}1000}}%
      \csname LTb\endcsname%
      \put(3963,484){\makebox(0,0){\strut{}10000}}%
      \csname LTb\endcsname%
      \put(5383,484){\makebox(0,0){\strut{}100000}}%
      \csname LTb\endcsname%
      \put(6803,484){\makebox(0,0){\strut{}1e+06}}%
      \put(3962,154){\makebox(0,0){\strut{}\text{DOFs}}}%
      \put(3962,4709){\makebox(0,0){\strut{}Laser on Square}}%
    }%
    \gplgaddtomacro\gplfronttext{%
      \csname LTb\endcsname%
      \put(5816,4206){\makebox(0,0)[r]{\strut{} \footnotesize Error in $\bar{|v|}$}}%
      \csname LTb\endcsname%
      \put(5816,3986){\makebox(0,0)[r]{\strut{} \footnotesize Error Estimator of $\bar{|v|}$}}%
      \csname LTb\endcsname%
      \put(5816,3766){\makebox(0,0)[r]{\strut{}$\mathcal{O}(\text{DOFs}^-\frac{3}{2})$}}%
    }%
    \gplbacktext
    \put(0,0){\includegraphics{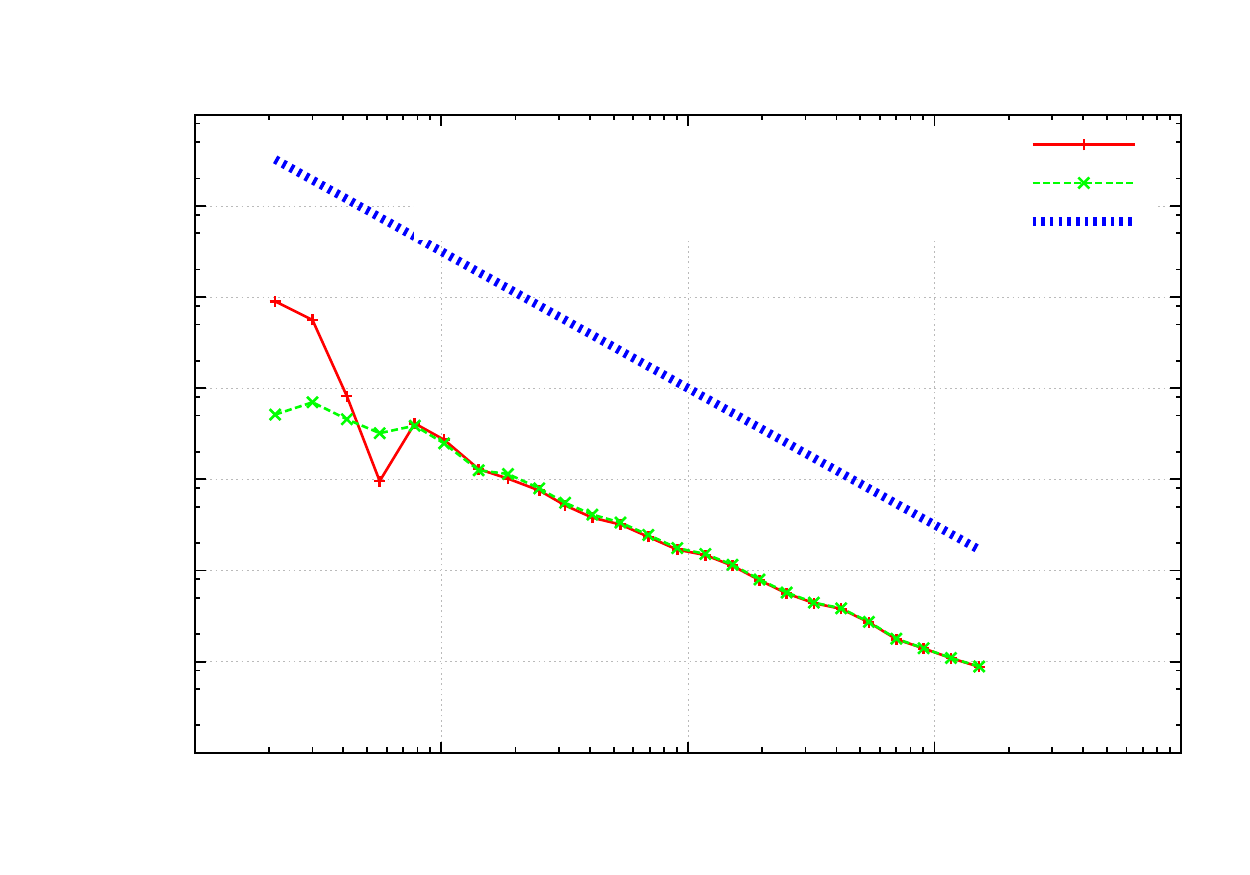}}%
    \gplfronttext
  \end{picture}%
\endgroup

%% file: Figures/LaseronSquareMeanVelociityIeff.tex
\begingroup
  \makeatletter
  \providecommand\color[2][]{%
    \GenericError{(gnuplot) \space\space\space\@spaces}{%
      Package color not loaded in conjunction with
      terminal option `colourtext'%
    }{See the gnuplot documentation for explanation.%
    }{Either use 'blacktext' in gnuplot or load the package
      color.sty in LaTeX.}%
    \renewcommand\color[2][]{}%
  }%
  \providecommand\includegraphics[2][]{%
    \GenericError{(gnuplot) \space\space\space\@spaces}{%
      Package graphicx or graphics not loaded%
    }{See the gnuplot documentation for explanation.%
    }{The gnuplot epslatex terminal needs graphicx.sty or graphics.sty.}%
    \renewcommand\includegraphics[2][]{}%
  }%
  \providecommand\rotatebox[2]{#2}%
  \@ifundefined{ifGPcolor}{%
    \newif\ifGPcolor
    \GPcolortrue
  }{}%
  \@ifundefined{ifGPblacktext}{%
    \newif\ifGPblacktext
    \GPblacktexttrue
  }{}%
  \let\gplgaddtomacro\g@addto@macro
  \gdef\gplbacktext{}%
  \gdef\gplfronttext{}%
  \makeatother
  \ifGPblacktext
    \def\colorrgb#1{}%
    \def\colorgray#1{}%
  \else
    \ifGPcolor
      \def\colorrgb#1{\color[rgb]{#1}}%
      \def\colorgray#1{\color[gray]{#1}}%
      \expandafter\def\csname LTw\endcsname{\color{white}}%
      \expandafter\def\csname LTb\endcsname{\color{black}}%
      \expandafter\def\csname LTa\endcsname{\color{black}}%
      \expandafter\def\csname LT0\endcsname{\color[rgb]{1,0,0}}%
      \expandafter\def\csname LT1\endcsname{\color[rgb]{0,1,0}}%
      \expandafter\def\csname LT2\endcsname{\color[rgb]{0,0,1}}%
      \expandafter\def\csname LT3\endcsname{\color[rgb]{1,0,1}}%
      \expandafter\def\csname LT4\endcsname{\color[rgb]{0,1,1}}%
      \expandafter\def\csname LT5\endcsname{\color[rgb]{1,1,0}}%
      \expandafter\def\csname LT6\endcsname{\color[rgb]{0,0,0}}%
      \expandafter\def\csname LT7\endcsname{\color[rgb]{1,0.3,0}}%
      \expandafter\def\csname LT8\endcsname{\color[rgb]{0.5,0.5,0.5}}%
    \else
      \def\colorrgb#1{\color{black}}%
      \def\colorgray#1{\color[gray]{#1}}%
      \expandafter\def\csname LTw\endcsname{\color{white}}%
      \expandafter\def\csname LTb\endcsname{\color{black}}%
      \expandafter\def\csname LTa\endcsname{\color{black}}%
      \expandafter\def\csname LT0\endcsname{\color{black}}%
      \expandafter\def\csname LT1\endcsname{\color{black}}%
      \expandafter\def\csname LT2\endcsname{\color{black}}%
      \expandafter\def\csname LT3\endcsname{\color{black}}%
      \expandafter\def\csname LT4\endcsname{\color{black}}%
      \expandafter\def\csname LT5\endcsname{\color{black}}%
      \expandafter\def\csname LT6\endcsname{\color{black}}%
      \expandafter\def\csname LT7\endcsname{\color{black}}%
      \expandafter\def\csname LT8\endcsname{\color{black}}%
    \fi
  \fi
  \setlength{\unitlength}{0.0500bp}%
  \begin{picture}(7200.00,5040.00)%
    \gplgaddtomacro\gplbacktext{%
      \csname LTb\endcsname%
      \put(594,704){\makebox(0,0)[r]{\strut{}0.1}}%
      \csname LTb\endcsname%
      \put(594,2542){\makebox(0,0)[r]{\strut{}1}}%
      \csname LTb\endcsname%
      \put(594,4379){\makebox(0,0)[r]{\strut{}10}}%
      \csname LTb\endcsname%
      \put(726,484){\makebox(0,0){\strut{}1000}}%
      \csname LTb\endcsname%
      \put(2752,484){\makebox(0,0){\strut{}10000}}%
      \csname LTb\endcsname%
      \put(4777,484){\makebox(0,0){\strut{}100000}}%
      \csname LTb\endcsname%
      \put(6803,484){\makebox(0,0){\strut{}1e+06}}%
      \put(3764,154){\makebox(0,0){\strut{}\text{DOFs}}}%
      \put(3764,4709){\makebox(0,0){\strut{}Laser on Square}}%
    }%
    \gplgaddtomacro\gplfronttext{%
      \csname LTb\endcsname%
      \put(2706,4206){\makebox(0,0)[r]{\strut{} \footnotesize $I_{eff}$ $\bar{|v|}$}}%
      \csname LTb\endcsname%
      \put(2706,3986){\makebox(0,0)[r]{\strut{} \footnotesize $I_{eff,p}$ $\bar{|v|}$}}%
      \csname LTb\endcsname%
      \put(2706,3766){\makebox(0,0)[r]{\strut{} \footnotesize $I_{eff,a}$ $\bar{|v|}$}}%
      \csname LTb\endcsname%
      \put(2706,3546){\makebox(0,0)[r]{\strut{}1}}%
    }%
    \gplbacktext
    \put(0,0){\includegraphics{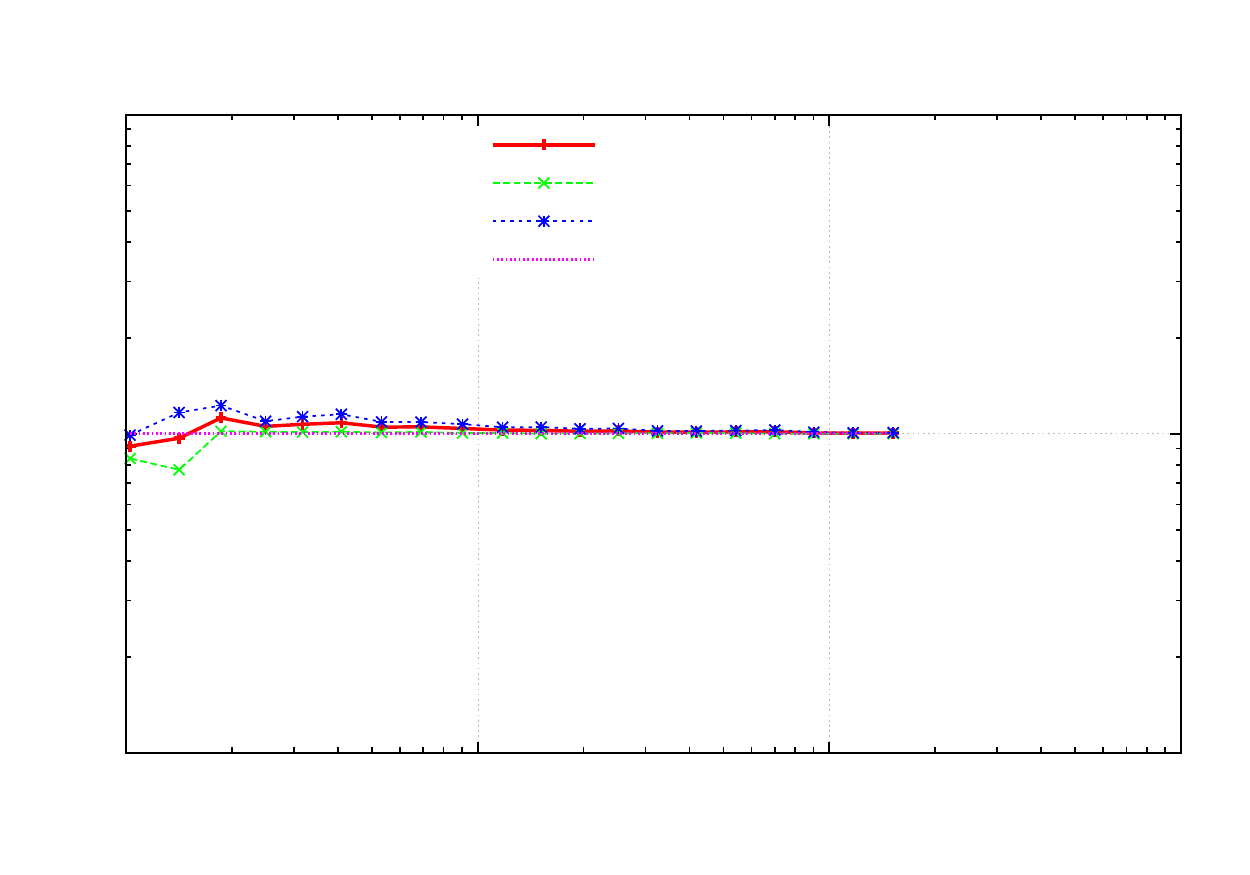}}%
    \gplfronttext
  \end{picture}%
\endgroup

%% file: Figures/LaseronSquareCombinedError.tex
\begingroup
  \makeatletter
  \providecommand\color[2][]{%
    \GenericError{(gnuplot) \space\space\space\@spaces}{%
      Package color not loaded in conjunction with
      terminal option `colourtext'%
    }{See the gnuplot documentation for explanation.%
    }{Either use 'blacktext' in gnuplot or load the package
      color.sty in LaTeX.}%
    \renewcommand\color[2][]{}%
  }%
  \providecommand\includegraphics[2][]{%
    \GenericError{(gnuplot) \space\space\space\@spaces}{%
      Package graphicx or graphics not loaded%
    }{See the gnuplot documentation for explanation.%
    }{The gnuplot epslatex terminal needs graphicx.sty or graphics.sty.}%
    \renewcommand\includegraphics[2][]{}%
  }%
  \providecommand\rotatebox[2]{#2}%
  \@ifundefined{ifGPcolor}{%
    \newif\ifGPcolor
    \GPcolortrue
  }{}%
  \@ifundefined{ifGPblacktext}{%
    \newif\ifGPblacktext
    \GPblacktexttrue
  }{}%
  \let\gplgaddtomacro\g@addto@macro
  \gdef\gplbacktext{}%
  \gdef\gplfronttext{}%
  \makeatother
  \ifGPblacktext
    \def\colorrgb#1{}%
    \def\colorgray#1{}%
  \else
    \ifGPcolor
      \def\colorrgb#1{\color[rgb]{#1}}%
      \def\colorgray#1{\color[gray]{#1}}%
      \expandafter\def\csname LTw\endcsname{\color{white}}%
      \expandafter\def\csname LTb\endcsname{\color{black}}%
      \expandafter\def\csname LTa\endcsname{\color{black}}%
      \expandafter\def\csname LT0\endcsname{\color[rgb]{1,0,0}}%
      \expandafter\def\csname LT1\endcsname{\color[rgb]{0,1,0}}%
      \expandafter\def\csname LT2\endcsname{\color[rgb]{0,0,1}}%
      \expandafter\def\csname LT3\endcsname{\color[rgb]{1,0,1}}%
      \expandafter\def\csname LT4\endcsname{\color[rgb]{0,1,1}}%
      \expandafter\def\csname LT5\endcsname{\color[rgb]{1,1,0}}%
      \expandafter\def\csname LT6\endcsname{\color[rgb]{0,0,0}}%
      \expandafter\def\csname LT7\endcsname{\color[rgb]{1,0.3,0}}%
      \expandafter\def\csname LT8\endcsname{\color[rgb]{0.5,0.5,0.5}}%
    \else
      \def\colorrgb#1{\color{black}}%
      \def\colorgray#1{\color[gray]{#1}}%
      \expandafter\def\csname LTw\endcsname{\color{white}}%
      \expandafter\def\csname LTb\endcsname{\color{black}}%
      \expandafter\def\csname LTa\endcsname{\color{black}}%
      \expandafter\def\csname LT0\endcsname{\color{black}}%
      \expandafter\def\csname LT1\endcsname{\color{black}}%
      \expandafter\def\csname LT2\endcsname{\color{black}}%
      \expandafter\def\csname LT3\endcsname{\color{black}}%
      \expandafter\def\csname LT4\endcsname{\color{black}}%
      \expandafter\def\csname LT5\endcsname{\color{black}}%
      \expandafter\def\csname LT6\endcsname{\color{black}}%
      \expandafter\def\csname LT7\endcsname{\color{black}}%
      \expandafter\def\csname LT8\endcsname{\color{black}}%
    \fi
  \fi
  \setlength{\unitlength}{0.0500bp}%
  \begin{picture}(7200.00,5040.00)%
    \gplgaddtomacro\gplbacktext{%
      \csname LTb\endcsname%
      \put(990,704){\makebox(0,0)[r]{\strut{}1e-10}}%
      \csname LTb\endcsname%
      \put(990,1229){\makebox(0,0)[r]{\strut{}1e-08}}%
      \csname LTb\endcsname%
      \put(990,1754){\makebox(0,0)[r]{\strut{}1e-06}}%
      \csname LTb\endcsname%
      \put(990,2279){\makebox(0,0)[r]{\strut{}0.0001}}%
      \csname LTb\endcsname%
      \put(990,2804){\makebox(0,0)[r]{\strut{}0.01}}%
      \csname LTb\endcsname%
      \put(990,3329){\makebox(0,0)[r]{\strut{}1}}%
      \csname LTb\endcsname%
      \put(990,3854){\makebox(0,0)[r]{\strut{}100}}%
      \csname LTb\endcsname%
      \put(990,4379){\makebox(0,0)[r]{\strut{}10000}}%
      \csname LTb\endcsname%
      \put(1122,484){\makebox(0,0){\strut{}100}}%
      \csname LTb\endcsname%
      \put(2542,484){\makebox(0,0){\strut{}1000}}%
      \csname LTb\endcsname%
      \put(3963,484){\makebox(0,0){\strut{}10000}}%
      \csname LTb\endcsname%
      \put(5383,484){\makebox(0,0){\strut{}100000}}%
      \csname LTb\endcsname%
      \put(6803,484){\makebox(0,0){\strut{}1e+06}}%
      \put(3962,154){\makebox(0,0){\strut{}\text{DOFs}}}%
      \put(3962,4709){\makebox(0,0){\strut{}Laser on Square}}%
    }%
    \gplgaddtomacro\gplfronttext{%
      \csname LTb\endcsname%
      \put(5816,4206){\makebox(0,0)[r]{\strut{} \footnotesize Error in $J_\mathfrak{E}$}}%
      \csname LTb\endcsname%
      \put(5816,3986){\makebox(0,0)[r]{\strut{} \footnotesize Error Estimator of $J_\mathfrak{E}$}}%
      \csname LTb\endcsname%
      \put(5816,3766){\makebox(0,0)[r]{\strut{} \footnotesize Error in $\bar{|v|}$}}%
      \csname LTb\endcsname%
      \put(5816,3546){\makebox(0,0)[r]{\strut{} \footnotesize Error in $\bar{T}$}}%
      \csname LTb\endcsname%
      \put(5816,3326){\makebox(0,0)[r]{\strut{}$\mathcal{O}(\text{DOFs}^-\frac{3}{2})$}}%
    }%
    \gplbacktext
    \put(0,0){\includegraphics{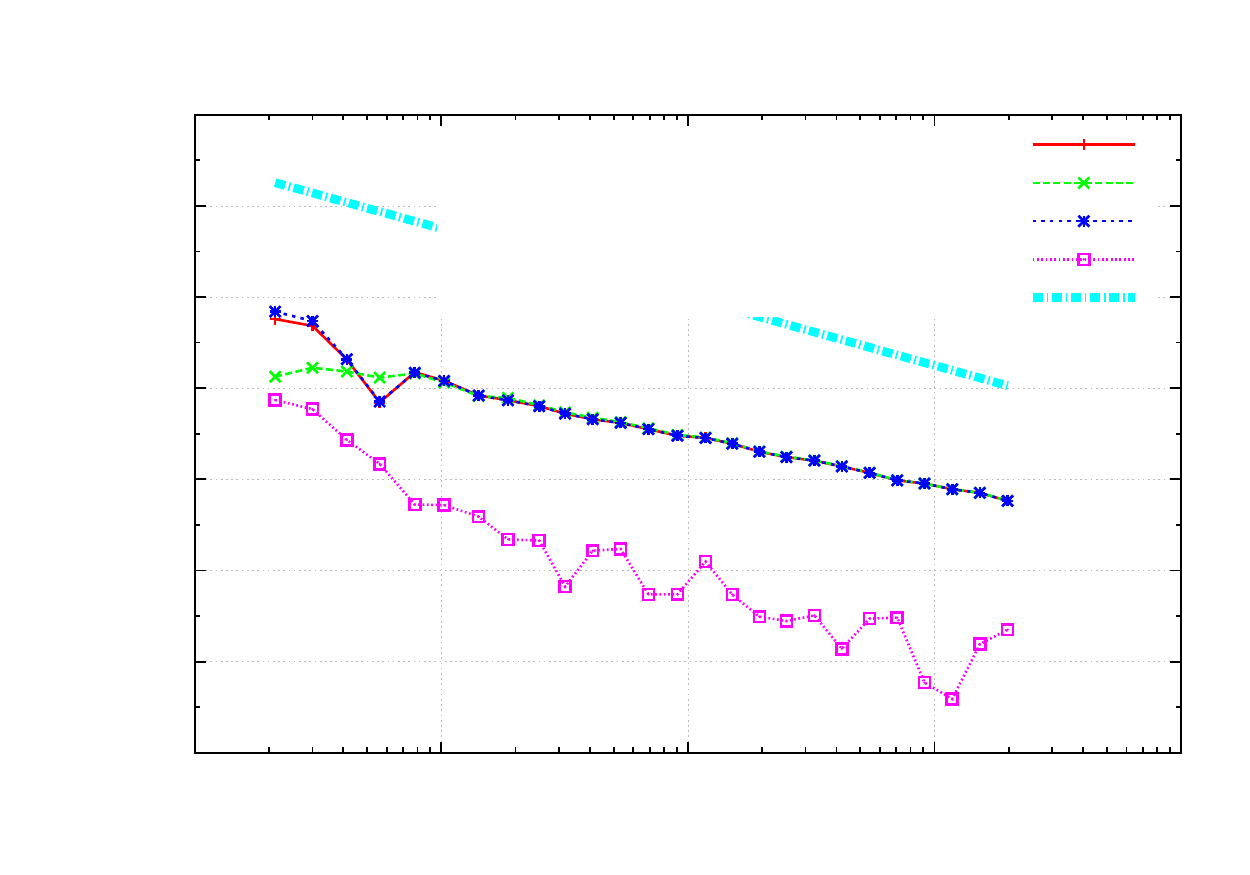}}%
    \gplfronttext
  \end{picture}%
\endgroup

%% file: Figures/LaseronSquareCombinedIeff.tex
\begingroup
  \makeatletter
  \providecommand\color[2][]{%
    \GenericError{(gnuplot) \space\space\space\@spaces}{%
      Package color not loaded in conjunction with
      terminal option `colourtext'%
    }{See the gnuplot documentation for explanation.%
    }{Either use 'blacktext' in gnuplot or load the package
      color.sty in LaTeX.}%
    \renewcommand\color[2][]{}%
  }%
  \providecommand\includegraphics[2][]{%
    \GenericError{(gnuplot) \space\space\space\@spaces}{%
      Package graphicx or graphics not loaded%
    }{See the gnuplot documentation for explanation.%
    }{The gnuplot epslatex terminal needs graphicx.sty or graphics.sty.}%
    \renewcommand\includegraphics[2][]{}%
  }%
  \providecommand\rotatebox[2]{#2}%
  \@ifundefined{ifGPcolor}{%
    \newif\ifGPcolor
    \GPcolortrue
  }{}%
  \@ifundefined{ifGPblacktext}{%
    \newif\ifGPblacktext
    \GPblacktexttrue
  }{}%
  \let\gplgaddtomacro\g@addto@macro
  \gdef\gplbacktext{}%
  \gdef\gplfronttext{}%
  \makeatother
  \ifGPblacktext
    \def\colorrgb#1{}%
    \def\colorgray#1{}%
  \else
    \ifGPcolor
      \def\colorrgb#1{\color[rgb]{#1}}%
      \def\colorgray#1{\color[gray]{#1}}%
      \expandafter\def\csname LTw\endcsname{\color{white}}%
      \expandafter\def\csname LTb\endcsname{\color{black}}%
      \expandafter\def\csname LTa\endcsname{\color{black}}%
      \expandafter\def\csname LT0\endcsname{\color[rgb]{1,0,0}}%
      \expandafter\def\csname LT1\endcsname{\color[rgb]{0,1,0}}%
      \expandafter\def\csname LT2\endcsname{\color[rgb]{0,0,1}}%
      \expandafter\def\csname LT3\endcsname{\color[rgb]{1,0,1}}%
      \expandafter\def\csname LT4\endcsname{\color[rgb]{0,1,1}}%
      \expandafter\def\csname LT5\endcsname{\color[rgb]{1,1,0}}%
      \expandafter\def\csname LT6\endcsname{\color[rgb]{0,0,0}}%
      \expandafter\def\csname LT7\endcsname{\color[rgb]{1,0.3,0}}%
      \expandafter\def\csname LT8\endcsname{\color[rgb]{0.5,0.5,0.5}}%
    \else
      \def\colorrgb#1{\color{black}}%
      \def\colorgray#1{\color[gray]{#1}}%
      \expandafter\def\csname LTw\endcsname{\color{white}}%
      \expandafter\def\csname LTb\endcsname{\color{black}}%
      \expandafter\def\csname LTa\endcsname{\color{black}}%
      \expandafter\def\csname LT0\endcsname{\color{black}}%
      \expandafter\def\csname LT1\endcsname{\color{black}}%
      \expandafter\def\csname LT2\endcsname{\color{black}}%
      \expandafter\def\csname LT3\endcsname{\color{black}}%
      \expandafter\def\csname LT4\endcsname{\color{black}}%
      \expandafter\def\csname LT5\endcsname{\color{black}}%
      \expandafter\def\csname LT6\endcsname{\color{black}}%
      \expandafter\def\csname LT7\endcsname{\color{black}}%
      \expandafter\def\csname LT8\endcsname{\color{black}}%
    \fi
  \fi
  \setlength{\unitlength}{0.0500bp}%
  \begin{picture}(7200.00,5040.00)%
    \gplgaddtomacro\gplbacktext{%
      \csname LTb\endcsname%
      \put(726,704){\makebox(0,0)[r]{\strut{}0.01}}%
      \csname LTb\endcsname%
      \put(726,1929){\makebox(0,0)[r]{\strut{}0.1}}%
      \csname LTb\endcsname%
      \put(726,3154){\makebox(0,0)[r]{\strut{}1}}%
      \csname LTb\endcsname%
      \put(726,4379){\makebox(0,0)[r]{\strut{}10}}%
      \csname LTb\endcsname%
      \put(858,484){\makebox(0,0){\strut{}100}}%
      \csname LTb\endcsname%
      \put(2344,484){\makebox(0,0){\strut{}1000}}%
      \csname LTb\endcsname%
      \put(3831,484){\makebox(0,0){\strut{}10000}}%
      \csname LTb\endcsname%
      \put(5317,484){\makebox(0,0){\strut{}100000}}%
      \csname LTb\endcsname%
      \put(6803,484){\makebox(0,0){\strut{}1e+06}}%
      \put(3830,154){\makebox(0,0){\strut{}\text{DOFs}}}%
      \put(3830,4709){\makebox(0,0){\strut{}Laser on Square}}%
    }%
    \gplgaddtomacro\gplfronttext{%
      \csname LTb\endcsname%
      \put(5816,4206){\makebox(0,0)[r]{\strut{} \footnotesize  $I_{eff}$ $J_\mathfrak{E}$}}%
      \csname LTb\endcsname%
      \put(5816,3986){\makebox(0,0)[r]{\strut{} \footnotesize  $I_{eff,p}$ $J_\mathfrak{E}$}}%
      \csname LTb\endcsname%
      \put(5816,3766){\makebox(0,0)[r]{\strut{} \footnotesize  $I_{eff,a}$ $J_\mathfrak{E}$}}%
      \csname LTb\endcsname%
      \put(5816,3546){\makebox(0,0)[r]{\strut{}1}}%
    }%
    \gplbacktext
    \put(0,0){\includegraphics{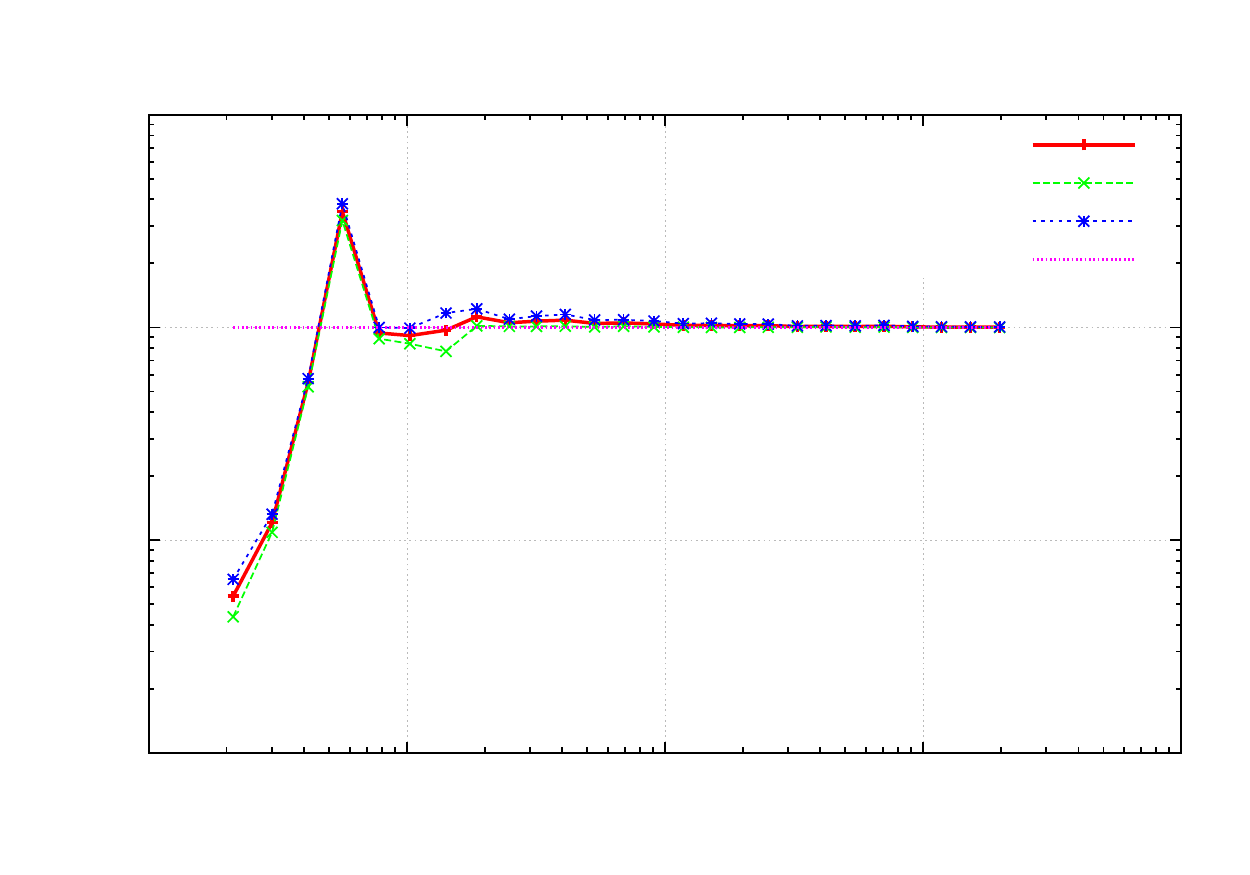}}%
    \gplfronttext
  \end{picture}%
\endgroup

%% file: Figures/YSplitterErrorConfig1.tex
\begingroup
  \makeatletter
  \providecommand\color[2][]{%
    \GenericError{(gnuplot) \space\space\space\@spaces}{%
      Package color not loaded in conjunction with
      terminal option `colourtext'%
    }{See the gnuplot documentation for explanation.%
    }{Either use 'blacktext' in gnuplot or load the package
      color.sty in LaTeX.}%
    \renewcommand\color[2][]{}%
  }%
  \providecommand\includegraphics[2][]{%
    \GenericError{(gnuplot) \space\space\space\@spaces}{%
      Package graphicx or graphics not loaded%
    }{See the gnuplot documentation for explanation.%
    }{The gnuplot epslatex terminal needs graphicx.sty or graphics.sty.}%
    \renewcommand\includegraphics[2][]{}%
  }%
  \providecommand\rotatebox[2]{#2}%
  \@ifundefined{ifGPcolor}{%
    \newif\ifGPcolor
    \GPcolortrue
  }{}%
  \@ifundefined{ifGPblacktext}{%
    \newif\ifGPblacktext
    \GPblacktexttrue
  }{}%
  \let\gplgaddtomacro\g@addto@macro
  \gdef\gplbacktext{}%
  \gdef\gplfronttext{}%
  \makeatother
  \ifGPblacktext
    \def\colorrgb#1{}%
    \def\colorgray#1{}%
  \else
    \ifGPcolor
      \def\colorrgb#1{\color[rgb]{#1}}%
      \def\colorgray#1{\color[gray]{#1}}%
      \expandafter\def\csname LTw\endcsname{\color{white}}%
      \expandafter\def\csname LTb\endcsname{\color{black}}%
      \expandafter\def\csname LTa\endcsname{\color{black}}%
      \expandafter\def\csname LT0\endcsname{\color[rgb]{1,0,0}}%
      \expandafter\def\csname LT1\endcsname{\color[rgb]{0,1,0}}%
      \expandafter\def\csname LT2\endcsname{\color[rgb]{0,0,1}}%
      \expandafter\def\csname LT3\endcsname{\color[rgb]{1,0,1}}%
      \expandafter\def\csname LT4\endcsname{\color[rgb]{0,1,1}}%
      \expandafter\def\csname LT5\endcsname{\color[rgb]{1,1,0}}%
      \expandafter\def\csname LT6\endcsname{\color[rgb]{0,0,0}}%
      \expandafter\def\csname LT7\endcsname{\color[rgb]{1,0.3,0}}%
      \expandafter\def\csname LT8\endcsname{\color[rgb]{0.5,0.5,0.5}}%
    \else
      \def\colorrgb#1{\color{black}}%
      \def\colorgray#1{\color[gray]{#1}}%
      \expandafter\def\csname LTw\endcsname{\color{white}}%
      \expandafter\def\csname LTb\endcsname{\color{black}}%
      \expandafter\def\csname LTa\endcsname{\color{black}}%
      \expandafter\def\csname LT0\endcsname{\color{black}}%
      \expandafter\def\csname LT1\endcsname{\color{black}}%
      \expandafter\def\csname LT2\endcsname{\color{black}}%
      \expandafter\def\csname LT3\endcsname{\color{black}}%
      \expandafter\def\csname LT4\endcsname{\color{black}}%
      \expandafter\def\csname LT5\endcsname{\color{black}}%
      \expandafter\def\csname LT6\endcsname{\color{black}}%
      \expandafter\def\csname LT7\endcsname{\color{black}}%
      \expandafter\def\csname LT8\endcsname{\color{black}}%
    \fi
  \fi
  \setlength{\unitlength}{0.0500bp}%
  \begin{picture}(7200.00,5040.00)%
    \gplgaddtomacro\gplbacktext{%
      \csname LTb\endcsname%
      \put(990,704){\makebox(0,0)[r]{\strut{}1e-08}}%
      \csname LTb\endcsname%
      \put(990,1163){\makebox(0,0)[r]{\strut{}1e-07}}%
      \csname LTb\endcsname%
      \put(990,1623){\makebox(0,0)[r]{\strut{}1e-06}}%
      \csname LTb\endcsname%
      \put(990,2082){\makebox(0,0)[r]{\strut{}1e-05}}%
      \csname LTb\endcsname%
      \put(990,2542){\makebox(0,0)[r]{\strut{}0.0001}}%
      \csname LTb\endcsname%
      \put(990,3001){\makebox(0,0)[r]{\strut{}0.001}}%
      \csname LTb\endcsname%
      \put(990,3460){\makebox(0,0)[r]{\strut{}0.01}}%
      \csname LTb\endcsname%
      \put(990,3920){\makebox(0,0)[r]{\strut{}0.1}}%
      \csname LTb\endcsname%
      \put(990,4379){\makebox(0,0)[r]{\strut{}1}}%
      \csname LTb\endcsname%
      \put(1122,484){\makebox(0,0){\strut{}100}}%
      \csname LTb\endcsname%
      \put(2146,484){\makebox(0,0){\strut{}1000}}%
      \csname LTb\endcsname%
      \put(3171,484){\makebox(0,0){\strut{}10000}}%
      \csname LTb\endcsname%
      \put(4195,484){\makebox(0,0){\strut{}100000}}%
      \csname LTb\endcsname%
      \put(5219,484){\makebox(0,0){\strut{}1e+06}}%
      \put(3170,154){\makebox(0,0){\strut{}\text{DOFs}}}%
      \put(3170,4709){\makebox(0,0){\strut{}Configuration 1}}%
    }%
    \gplgaddtomacro\gplfronttext{%
      \csname LTb\endcsname%
      \put(6212,4269){\makebox(0,0)[r]{\strut{} \footnotesize Error in $J_c$}}%
      \csname LTb\endcsname%
      \put(6212,4049){\makebox(0,0)[r]{\strut{} \footnotesize  $\eta_h$}}%
      \csname LTb\endcsname%
      \put(6212,3829){\makebox(0,0)[r]{\strut{} \footnotesize Error in $J_1$}}%
      \csname LTb\endcsname%
      \put(6212,3609){\makebox(0,0)[r]{\strut{} \footnotesize Error in $J_2$}}%
      \csname LTb\endcsname%
      \put(6212,3389){\makebox(0,0)[r]{\strut{} \footnotesize Error in $J_4$}}%
      \csname LTb\endcsname%
      \put(6212,3169){\makebox(0,0)[r]{\strut{} \footnotesize Error in $J_5$}}%
      \csname LTb\endcsname%
      \put(6212,2949){\makebox(0,0)[r]{\strut{} \footnotesize $\text{DOFs}^{-1}$}}%
    }%
    \gplbacktext
    \put(0,0){\includegraphics{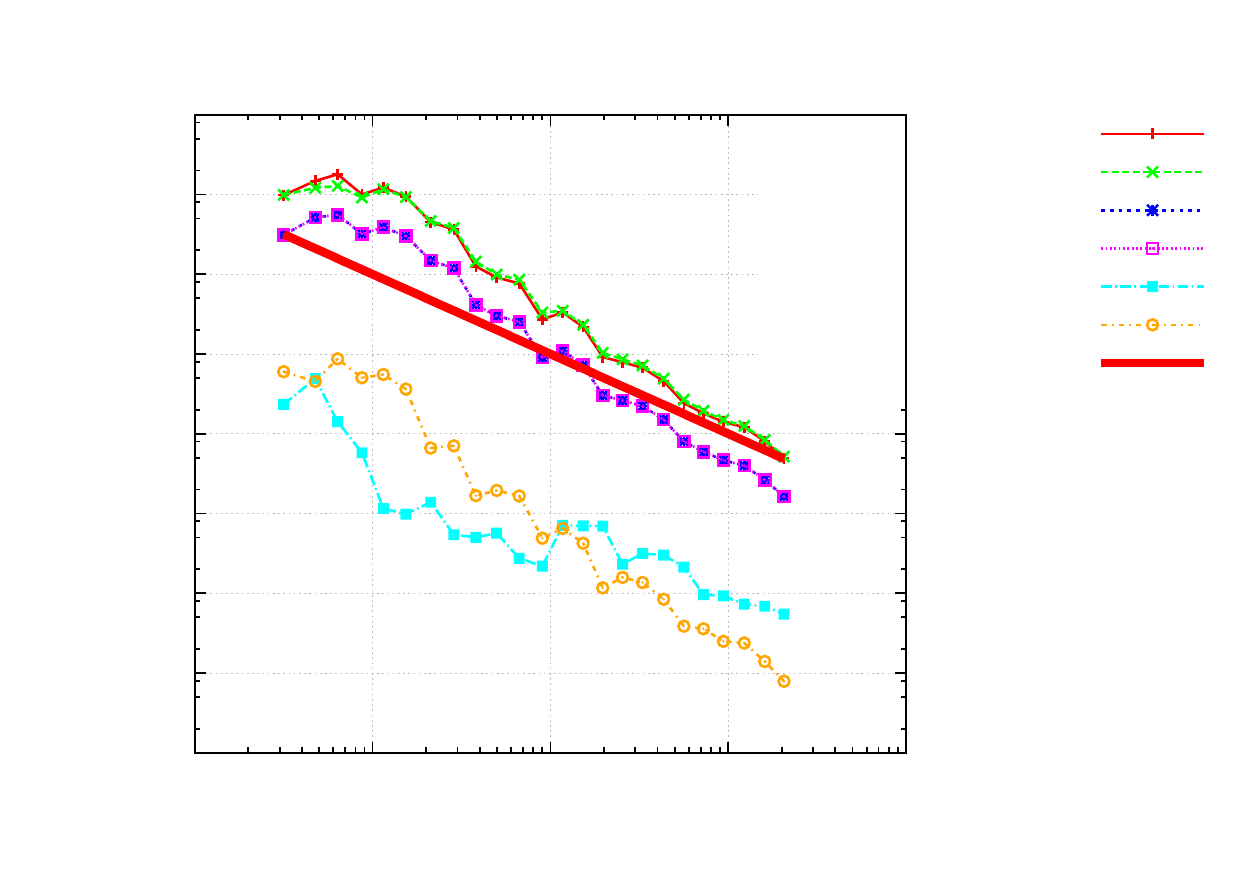}}%
    \gplfronttext
  \end{picture}%
\endgroup

%% file: Figures/YSplitterErrorConfig2.tex
\begingroup
  \makeatletter
  \providecommand\color[2][]{%
    \GenericError{(gnuplot) \space\space\space\@spaces}{%
      Package color not loaded in conjunction with
      terminal option `colourtext'%
    }{See the gnuplot documentation for explanation.%
    }{Either use 'blacktext' in gnuplot or load the package
      color.sty in LaTeX.}%
    \renewcommand\color[2][]{}%
  }%
  \providecommand\includegraphics[2][]{%
    \GenericError{(gnuplot) \space\space\space\@spaces}{%
      Package graphicx or graphics not loaded%
    }{See the gnuplot documentation for explanation.%
    }{The gnuplot epslatex terminal needs graphicx.sty or graphics.sty.}%
    \renewcommand\includegraphics[2][]{}%
  }%
  \providecommand\rotatebox[2]{#2}%
  \@ifundefined{ifGPcolor}{%
    \newif\ifGPcolor
    \GPcolortrue
  }{}%
  \@ifundefined{ifGPblacktext}{%
    \newif\ifGPblacktext
    \GPblacktexttrue
  }{}%
  \let\gplgaddtomacro\g@addto@macro
  \gdef\gplbacktext{}%
  \gdef\gplfronttext{}%
  \makeatother
  \ifGPblacktext
    \def\colorrgb#1{}%
    \def\colorgray#1{}%
  \else
    \ifGPcolor
      \def\colorrgb#1{\color[rgb]{#1}}%
      \def\colorgray#1{\color[gray]{#1}}%
      \expandafter\def\csname LTw\endcsname{\color{white}}%
      \expandafter\def\csname LTb\endcsname{\color{black}}%
      \expandafter\def\csname LTa\endcsname{\color{black}}%
      \expandafter\def\csname LT0\endcsname{\color[rgb]{1,0,0}}%
      \expandafter\def\csname LT1\endcsname{\color[rgb]{0,1,0}}%
      \expandafter\def\csname LT2\endcsname{\color[rgb]{0,0,1}}%
      \expandafter\def\csname LT3\endcsname{\color[rgb]{1,0,1}}%
      \expandafter\def\csname LT4\endcsname{\color[rgb]{0,1,1}}%
      \expandafter\def\csname LT5\endcsname{\color[rgb]{1,1,0}}%
      \expandafter\def\csname LT6\endcsname{\color[rgb]{0,0,0}}%
      \expandafter\def\csname LT7\endcsname{\color[rgb]{1,0.3,0}}%
      \expandafter\def\csname LT8\endcsname{\color[rgb]{0.5,0.5,0.5}}%
    \else
      \def\colorrgb#1{\color{black}}%
      \def\colorgray#1{\color[gray]{#1}}%
      \expandafter\def\csname LTw\endcsname{\color{white}}%
      \expandafter\def\csname LTb\endcsname{\color{black}}%
      \expandafter\def\csname LTa\endcsname{\color{black}}%
      \expandafter\def\csname LT0\endcsname{\color{black}}%
      \expandafter\def\csname LT1\endcsname{\color{black}}%
      \expandafter\def\csname LT2\endcsname{\color{black}}%
      \expandafter\def\csname LT3\endcsname{\color{black}}%
      \expandafter\def\csname LT4\endcsname{\color{black}}%
      \expandafter\def\csname LT5\endcsname{\color{black}}%
      \expandafter\def\csname LT6\endcsname{\color{black}}%
      \expandafter\def\csname LT7\endcsname{\color{black}}%
      \expandafter\def\csname LT8\endcsname{\color{black}}%
    \fi
  \fi
  \setlength{\unitlength}{0.0500bp}%
  \begin{picture}(7200.00,5040.00)%
    \gplgaddtomacro\gplbacktext{%
      \csname LTb\endcsname%
      \put(990,704){\makebox(0,0)[r]{\strut{}1e-08}}%
      \csname LTb\endcsname%
      \put(990,1112){\makebox(0,0)[r]{\strut{}1e-07}}%
      \csname LTb\endcsname%
      \put(990,1521){\makebox(0,0)[r]{\strut{}1e-06}}%
      \csname LTb\endcsname%
      \put(990,1929){\makebox(0,0)[r]{\strut{}1e-05}}%
      \csname LTb\endcsname%
      \put(990,2337){\makebox(0,0)[r]{\strut{}0.0001}}%
      \csname LTb\endcsname%
      \put(990,2746){\makebox(0,0)[r]{\strut{}0.001}}%
      \csname LTb\endcsname%
      \put(990,3154){\makebox(0,0)[r]{\strut{}0.01}}%
      \csname LTb\endcsname%
      \put(990,3562){\makebox(0,0)[r]{\strut{}0.1}}%
      \csname LTb\endcsname%
      \put(990,3971){\makebox(0,0)[r]{\strut{}1}}%
      \csname LTb\endcsname%
      \put(990,4379){\makebox(0,0)[r]{\strut{}10}}%
      \csname LTb\endcsname%
      \put(1122,484){\makebox(0,0){\strut{}100}}%
      \csname LTb\endcsname%
      \put(2146,484){\makebox(0,0){\strut{}1000}}%
      \csname LTb\endcsname%
      \put(3171,484){\makebox(0,0){\strut{}10000}}%
      \csname LTb\endcsname%
      \put(4195,484){\makebox(0,0){\strut{}100000}}%
      \csname LTb\endcsname%
      \put(5219,484){\makebox(0,0){\strut{}1e+06}}%
      \put(3170,154){\makebox(0,0){\strut{}\text{DOFs}}}%
      \put(3170,4709){\makebox(0,0){\strut{}Configuration 2}}%
    }%
    \gplgaddtomacro\gplfronttext{%
      \csname LTb\endcsname%
      \put(6212,4269){\makebox(0,0)[r]{\strut{} \footnotesize Error in $J_c$}}%
      \csname LTb\endcsname%
      \put(6212,4049){\makebox(0,0)[r]{\strut{} \footnotesize  $\eta_h$}}%
      \csname LTb\endcsname%
      \put(6212,3829){\makebox(0,0)[r]{\strut{} \footnotesize Error in $J_1$}}%
      \csname LTb\endcsname%
      \put(6212,3609){\makebox(0,0)[r]{\strut{} \footnotesize Error in $J_2$}}%
      \csname LTb\endcsname%
      \put(6212,3389){\makebox(0,0)[r]{\strut{} \footnotesize Error in $J_4$}}%
      \csname LTb\endcsname%
      \put(6212,3169){\makebox(0,0)[r]{\strut{} \footnotesize Error in $J_5$}}%
      \csname LTb\endcsname%
      \put(6212,2949){\makebox(0,0)[r]{\strut{} \footnotesize $\text{DOFs}^{-1}$}}%
    }%
    \gplbacktext
    \put(0,0){\includegraphics{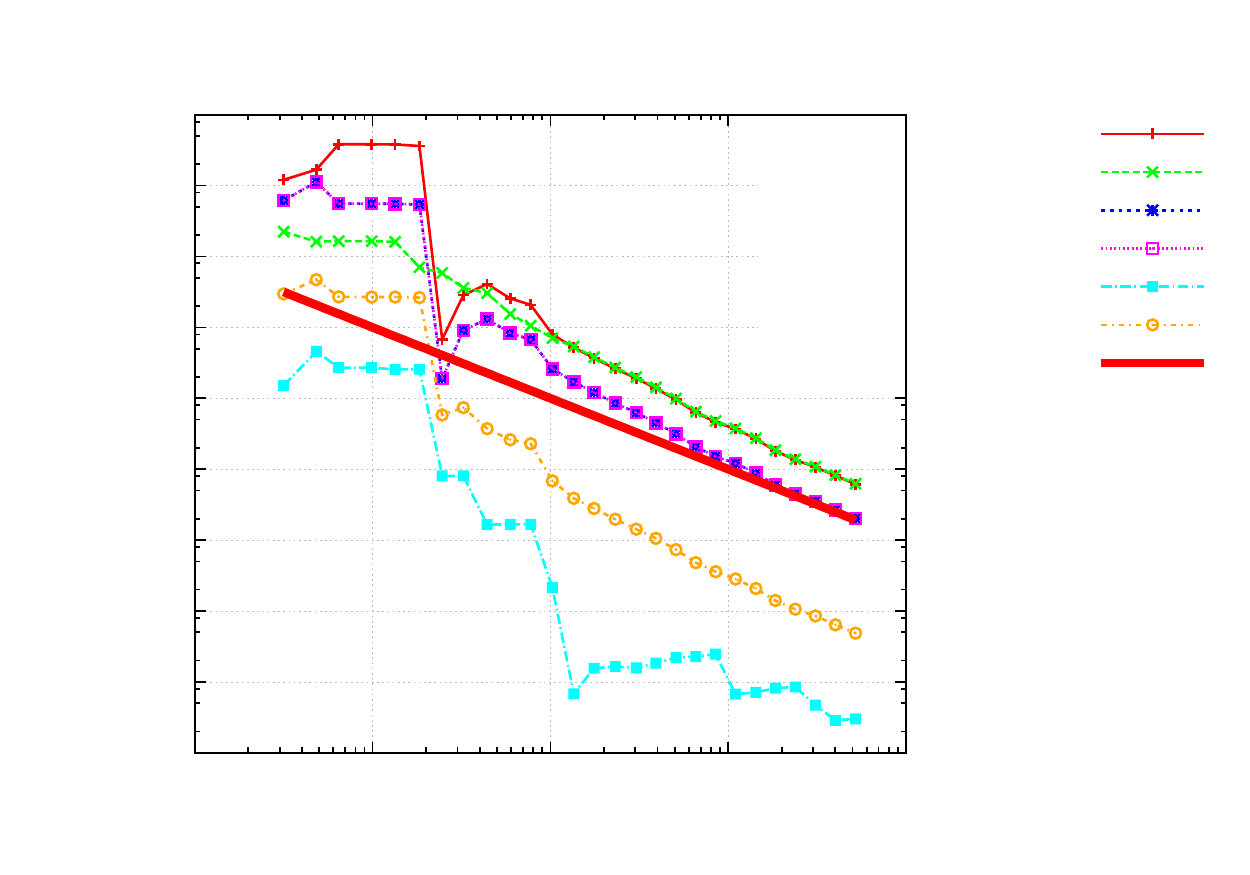}}%
    \gplfronttext
  \end{picture}%
\endgroup

%% file: Figures/YSplitterErrorConfig3.tex
\begingroup
  \makeatletter
  \providecommand\color[2][]{%
    \GenericError{(gnuplot) \space\space\space\@spaces}{%
      Package color not loaded in conjunction with
      terminal option `colourtext'%
    }{See the gnuplot documentation for explanation.%
    }{Either use 'blacktext' in gnuplot or load the package
      color.sty in LaTeX.}%
    \renewcommand\color[2][]{}%
  }%
  \providecommand\includegraphics[2][]{%
    \GenericError{(gnuplot) \space\space\space\@spaces}{%
      Package graphicx or graphics not loaded%
    }{See the gnuplot documentation for explanation.%
    }{The gnuplot epslatex terminal needs graphicx.sty or graphics.sty.}%
    \renewcommand\includegraphics[2][]{}%
  }%
  \providecommand\rotatebox[2]{#2}%
  \@ifundefined{ifGPcolor}{%
    \newif\ifGPcolor
    \GPcolortrue
  }{}%
  \@ifundefined{ifGPblacktext}{%
    \newif\ifGPblacktext
    \GPblacktexttrue
  }{}%
  \let\gplgaddtomacro\g@addto@macro
  \gdef\gplbacktext{}%
  \gdef\gplfronttext{}%
  \makeatother
  \ifGPblacktext
    \def\colorrgb#1{}%
    \def\colorgray#1{}%
  \else
    \ifGPcolor
      \def\colorrgb#1{\color[rgb]{#1}}%
      \def\colorgray#1{\color[gray]{#1}}%
      \expandafter\def\csname LTw\endcsname{\color{white}}%
      \expandafter\def\csname LTb\endcsname{\color{black}}%
      \expandafter\def\csname LTa\endcsname{\color{black}}%
      \expandafter\def\csname LT0\endcsname{\color[rgb]{1,0,0}}%
      \expandafter\def\csname LT1\endcsname{\color[rgb]{0,1,0}}%
      \expandafter\def\csname LT2\endcsname{\color[rgb]{0,0,1}}%
      \expandafter\def\csname LT3\endcsname{\color[rgb]{1,0,1}}%
      \expandafter\def\csname LT4\endcsname{\color[rgb]{0,1,1}}%
      \expandafter\def\csname LT5\endcsname{\color[rgb]{1,1,0}}%
      \expandafter\def\csname LT6\endcsname{\color[rgb]{0,0,0}}%
      \expandafter\def\csname LT7\endcsname{\color[rgb]{1,0.3,0}}%
      \expandafter\def\csname LT8\endcsname{\color[rgb]{0.5,0.5,0.5}}%
    \else
      \def\colorrgb#1{\color{black}}%
      \def\colorgray#1{\color[gray]{#1}}%
      \expandafter\def\csname LTw\endcsname{\color{white}}%
      \expandafter\def\csname LTb\endcsname{\color{black}}%
      \expandafter\def\csname LTa\endcsname{\color{black}}%
      \expandafter\def\csname LT0\endcsname{\color{black}}%
      \expandafter\def\csname LT1\endcsname{\color{black}}%
      \expandafter\def\csname LT2\endcsname{\color{black}}%
      \expandafter\def\csname LT3\endcsname{\color{black}}%
      \expandafter\def\csname LT4\endcsname{\color{black}}%
      \expandafter\def\csname LT5\endcsname{\color{black}}%
      \expandafter\def\csname LT6\endcsname{\color{black}}%
      \expandafter\def\csname LT7\endcsname{\color{black}}%
      \expandafter\def\csname LT8\endcsname{\color{black}}%
    \fi
  \fi
  \setlength{\unitlength}{0.0500bp}%
  \begin{picture}(7200.00,5040.00)%
    \gplgaddtomacro\gplbacktext{%
      \csname LTb\endcsname%
      \put(990,704){\makebox(0,0)[r]{\strut{}1e-08}}%
      \csname LTb\endcsname%
      \put(990,1038){\makebox(0,0)[r]{\strut{}1e-07}}%
      \csname LTb\endcsname%
      \put(990,1372){\makebox(0,0)[r]{\strut{}1e-06}}%
      \csname LTb\endcsname%
      \put(990,1706){\makebox(0,0)[r]{\strut{}1e-05}}%
      \csname LTb\endcsname%
      \put(990,2040){\makebox(0,0)[r]{\strut{}0.0001}}%
      \csname LTb\endcsname%
      \put(990,2374){\makebox(0,0)[r]{\strut{}0.001}}%
      \csname LTb\endcsname%
      \put(990,2709){\makebox(0,0)[r]{\strut{}0.01}}%
      \csname LTb\endcsname%
      \put(990,3043){\makebox(0,0)[r]{\strut{}0.1}}%
      \csname LTb\endcsname%
      \put(990,3377){\makebox(0,0)[r]{\strut{}1}}%
      \csname LTb\endcsname%
      \put(990,3711){\makebox(0,0)[r]{\strut{}10}}%
      \csname LTb\endcsname%
      \put(990,4045){\makebox(0,0)[r]{\strut{}100}}%
      \csname LTb\endcsname%
      \put(990,4379){\makebox(0,0)[r]{\strut{}1000}}%
      \csname LTb\endcsname%
      \put(1122,484){\makebox(0,0){\strut{}100}}%
      \csname LTb\endcsname%
      \put(2146,484){\makebox(0,0){\strut{}1000}}%
      \csname LTb\endcsname%
      \put(3171,484){\makebox(0,0){\strut{}10000}}%
      \csname LTb\endcsname%
      \put(4195,484){\makebox(0,0){\strut{}100000}}%
      \csname LTb\endcsname%
      \put(5219,484){\makebox(0,0){\strut{}1e+06}}%
      \put(3170,154){\makebox(0,0){\strut{}\text{DOFs}}}%
      \put(3170,4709){\makebox(0,0){\strut{}Configuration 3}}%
    }%
    \gplgaddtomacro\gplfronttext{%
      \csname LTb\endcsname%
      \put(6212,4269){\makebox(0,0)[r]{\strut{} \footnotesize Error in $J_c$}}%
      \csname LTb\endcsname%
      \put(6212,4049){\makebox(0,0)[r]{\strut{} \footnotesize  $\eta_h$}}%
      \csname LTb\endcsname%
      \put(6212,3829){\makebox(0,0)[r]{\strut{} \footnotesize Error in $J_1$}}%
      \csname LTb\endcsname%
      \put(6212,3609){\makebox(0,0)[r]{\strut{} \footnotesize Error in $J_2$}}%
      \csname LTb\endcsname%
      \put(6212,3389){\makebox(0,0)[r]{\strut{} \footnotesize Error in $J_3$}}%
      \csname LTb\endcsname%
      \put(6212,3169){\makebox(0,0)[r]{\strut{} \footnotesize Error in $J_4$}}%
      \csname LTb\endcsname%
      \put(6212,2949){\makebox(0,0)[r]{\strut{} \footnotesize Error in $J_5$}}%
      \csname LTb\endcsname%
      \put(6212,2729){\makebox(0,0)[r]{\strut{} \footnotesize Error in $J_6$}}%
      \csname LTb\endcsname%
      \put(6212,2509){\makebox(0,0)[r]{\strut{} \footnotesize $\text{DOFs}^{-1}$}}%
    }%
    \gplbacktext
    \put(0,0){\includegraphics{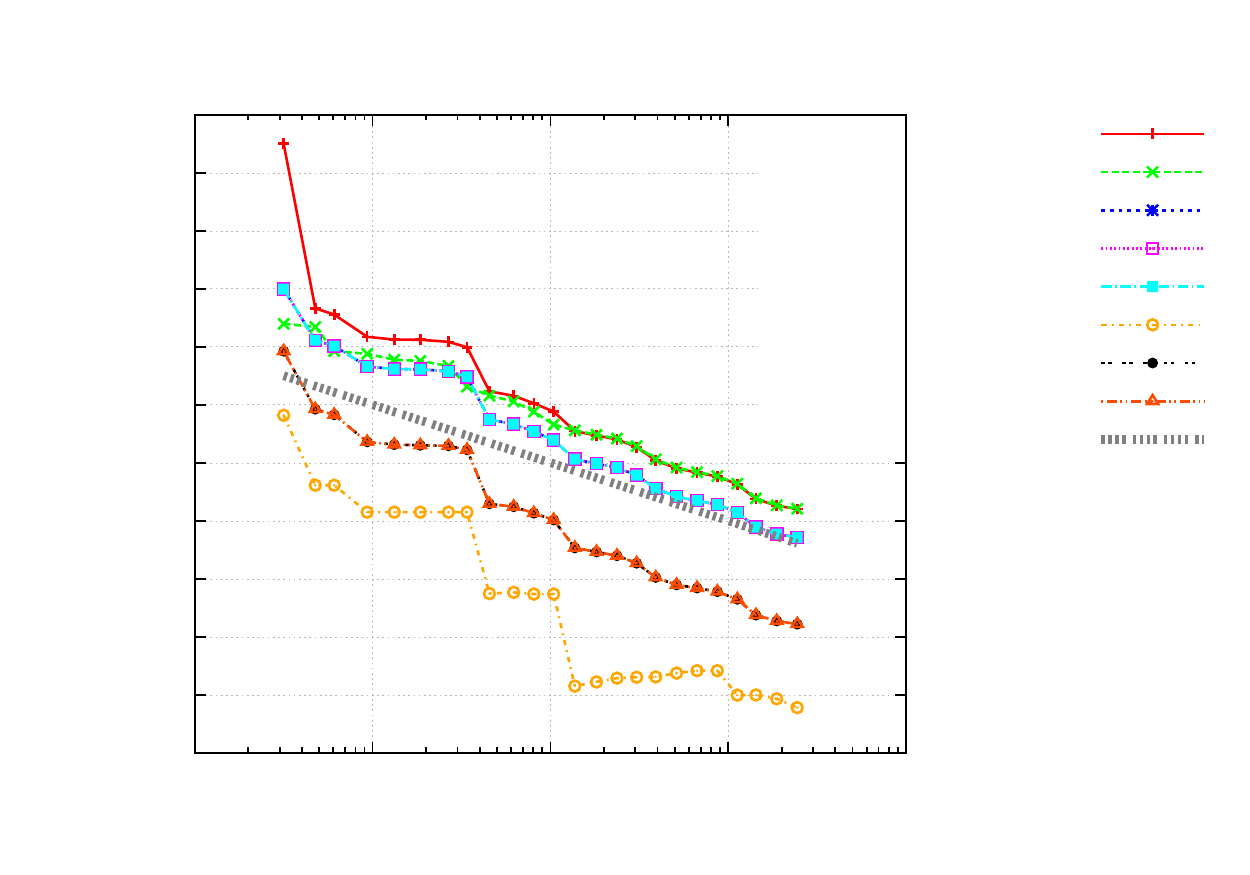}}%
    \gplfronttext
  \end{picture}%
\endgroup

%% file: Figures/YSplitterErrorConfig4.tex
\begingroup
  \makeatletter
  \providecommand\color[2][]{%
    \GenericError{(gnuplot) \space\space\space\@spaces}{%
      Package color not loaded in conjunction with
      terminal option `colourtext'%
    }{See the gnuplot documentation for explanation.%
    }{Either use 'blacktext' in gnuplot or load the package
      color.sty in LaTeX.}%
    \renewcommand\color[2][]{}%
  }%
  \providecommand\includegraphics[2][]{%
    \GenericError{(gnuplot) \space\space\space\@spaces}{%
      Package graphicx or graphics not loaded%
    }{See the gnuplot documentation for explanation.%
    }{The gnuplot epslatex terminal needs graphicx.sty or graphics.sty.}%
    \renewcommand\includegraphics[2][]{}%
  }%
  \providecommand\rotatebox[2]{#2}%
  \@ifundefined{ifGPcolor}{%
    \newif\ifGPcolor
    \GPcolortrue
  }{}%
  \@ifundefined{ifGPblacktext}{%
    \newif\ifGPblacktext
    \GPblacktexttrue
  }{}%
  \let\gplgaddtomacro\g@addto@macro
  \gdef\gplbacktext{}%
  \gdef\gplfronttext{}%
  \makeatother
  \ifGPblacktext
    \def\colorrgb#1{}%
    \def\colorgray#1{}%
  \else
    \ifGPcolor
      \def\colorrgb#1{\color[rgb]{#1}}%
      \def\colorgray#1{\color[gray]{#1}}%
      \expandafter\def\csname LTw\endcsname{\color{white}}%
      \expandafter\def\csname LTb\endcsname{\color{black}}%
      \expandafter\def\csname LTa\endcsname{\color{black}}%
      \expandafter\def\csname LT0\endcsname{\color[rgb]{1,0,0}}%
      \expandafter\def\csname LT1\endcsname{\color[rgb]{0,1,0}}%
      \expandafter\def\csname LT2\endcsname{\color[rgb]{0,0,1}}%
      \expandafter\def\csname LT3\endcsname{\color[rgb]{1,0,1}}%
      \expandafter\def\csname LT4\endcsname{\color[rgb]{0,1,1}}%
      \expandafter\def\csname LT5\endcsname{\color[rgb]{1,1,0}}%
      \expandafter\def\csname LT6\endcsname{\color[rgb]{0,0,0}}%
      \expandafter\def\csname LT7\endcsname{\color[rgb]{1,0.3,0}}%
      \expandafter\def\csname LT8\endcsname{\color[rgb]{0.5,0.5,0.5}}%
    \else
      \def\colorrgb#1{\color{black}}%
      \def\colorgray#1{\color[gray]{#1}}%
      \expandafter\def\csname LTw\endcsname{\color{white}}%
      \expandafter\def\csname LTb\endcsname{\color{black}}%
      \expandafter\def\csname LTa\endcsname{\color{black}}%
      \expandafter\def\csname LT0\endcsname{\color{black}}%
      \expandafter\def\csname LT1\endcsname{\color{black}}%
      \expandafter\def\csname LT2\endcsname{\color{black}}%
      \expandafter\def\csname LT3\endcsname{\color{black}}%
      \expandafter\def\csname LT4\endcsname{\color{black}}%
      \expandafter\def\csname LT5\endcsname{\color{black}}%
      \expandafter\def\csname LT6\endcsname{\color{black}}%
      \expandafter\def\csname LT7\endcsname{\color{black}}%
      \expandafter\def\csname LT8\endcsname{\color{black}}%
    \fi
  \fi
  \setlength{\unitlength}{0.0500bp}%
  \begin{picture}(7200.00,5040.00)%
    \gplgaddtomacro\gplbacktext{%
      \csname LTb\endcsname%
      \put(990,704){\makebox(0,0)[r]{\strut{}1e-10}}%
      \csname LTb\endcsname%
      \put(990,1229){\makebox(0,0)[r]{\strut{}1e-08}}%
      \csname LTb\endcsname%
      \put(990,1754){\makebox(0,0)[r]{\strut{}1e-06}}%
      \csname LTb\endcsname%
      \put(990,2279){\makebox(0,0)[r]{\strut{}0.0001}}%
      \csname LTb\endcsname%
      \put(990,2804){\makebox(0,0)[r]{\strut{}0.01}}%
      \csname LTb\endcsname%
      \put(990,3329){\makebox(0,0)[r]{\strut{}1}}%
      \csname LTb\endcsname%
      \put(990,3854){\makebox(0,0)[r]{\strut{}100}}%
      \csname LTb\endcsname%
      \put(990,4379){\makebox(0,0)[r]{\strut{}10000}}%
      \csname LTb\endcsname%
      \put(1122,484){\makebox(0,0){\strut{}100}}%
      \csname LTb\endcsname%
      \put(2146,484){\makebox(0,0){\strut{}1000}}%
      \csname LTb\endcsname%
      \put(3171,484){\makebox(0,0){\strut{}10000}}%
      \csname LTb\endcsname%
      \put(4195,484){\makebox(0,0){\strut{}100000}}%
      \csname LTb\endcsname%
      \put(5219,484){\makebox(0,0){\strut{}1e+06}}%
      \put(3170,154){\makebox(0,0){\strut{}\text{DOFs}}}%
      \put(3170,4709){\makebox(0,0){\strut{}Configuration 4}}%
    }%
    \gplgaddtomacro\gplfronttext{%
      \csname LTb\endcsname%
      \put(6212,4269){\makebox(0,0)[r]{\strut{} \footnotesize Error in $J_c$}}%
      \csname LTb\endcsname%
      \put(6212,4049){\makebox(0,0)[r]{\strut{} \footnotesize  $\eta_h$}}%
      \csname LTb\endcsname%
      \put(6212,3829){\makebox(0,0)[r]{\strut{} \footnotesize Error in $J_1$}}%
      \csname LTb\endcsname%
      \put(6212,3609){\makebox(0,0)[r]{\strut{} \footnotesize Error in $J_2$}}%
      \csname LTb\endcsname%
      \put(6212,3389){\makebox(0,0)[r]{\strut{} \footnotesize Error in $J_3$}}%
      \csname LTb\endcsname%
      \put(6212,3169){\makebox(0,0)[r]{\strut{} \footnotesize Error in $J_4$}}%
      \csname LTb\endcsname%
      \put(6212,2949){\makebox(0,0)[r]{\strut{} \footnotesize Error in $J_5$}}%
      \csname LTb\endcsname%
      \put(6212,2729){\makebox(0,0)[r]{\strut{} \footnotesize Error in $J_6$}}%
      \csname LTb\endcsname%
      \put(6212,2509){\makebox(0,0)[r]{\strut{} \footnotesize Error in $J_7$}}%
      \csname LTb\endcsname%
      \put(6212,2289){\makebox(0,0)[r]{\strut{} \footnotesize $\text{DOFs}^{-1}$}}%
    }%
    \gplbacktext
    \put(0,0){\includegraphics{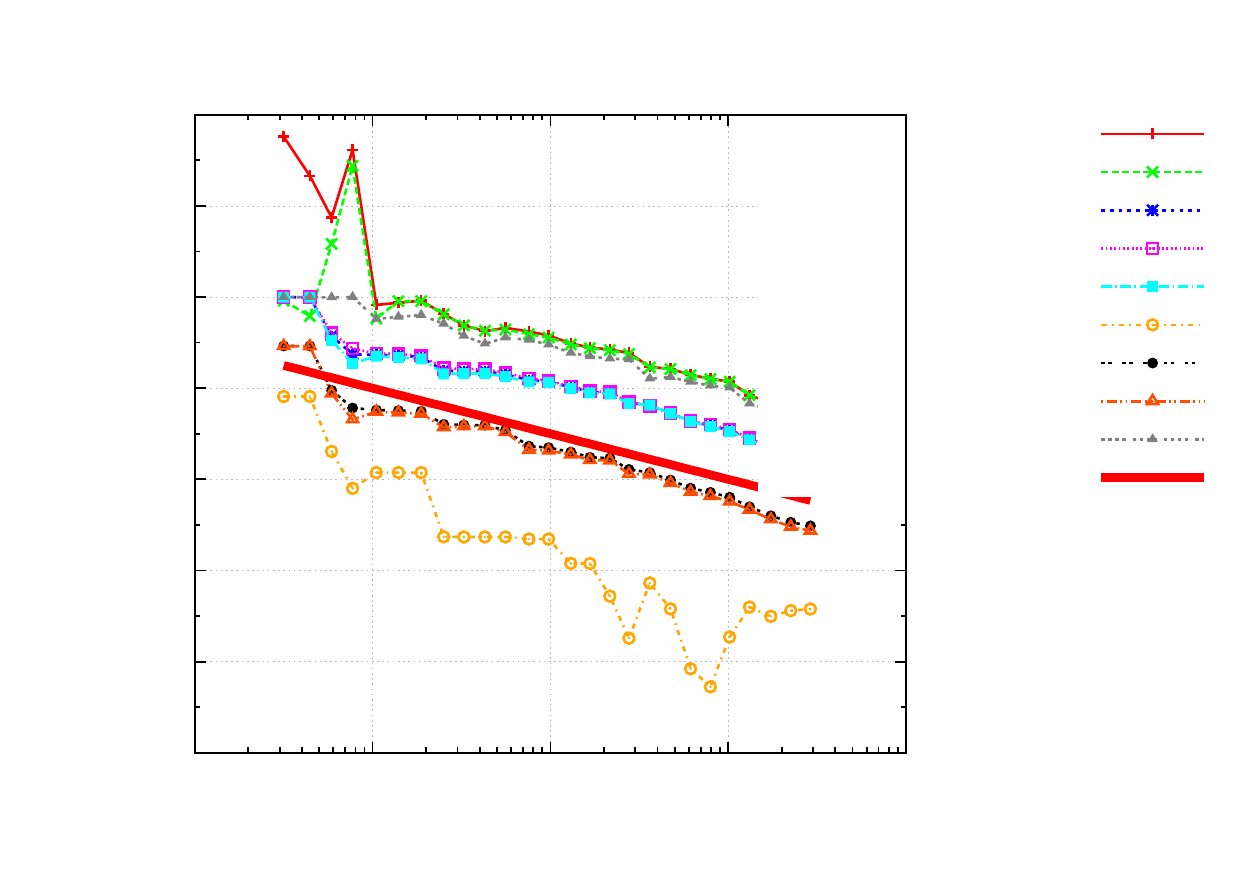}}%
    \gplfronttext
  \end{picture}%
\endgroup

%% file: Figures/YSplitterErrorConfig5.tex
\begingroup
  \makeatletter
  \providecommand\color[2][]{%
    \GenericError{(gnuplot) \space\space\space\@spaces}{%
      Package color not loaded in conjunction with
      terminal option `colourtext'%
    }{See the gnuplot documentation for explanation.%
    }{Either use 'blacktext' in gnuplot or load the package
      color.sty in LaTeX.}%
    \renewcommand\color[2][]{}%
  }%
  \providecommand\includegraphics[2][]{%
    \GenericError{(gnuplot) \space\space\space\@spaces}{%
      Package graphicx or graphics not loaded%
    }{See the gnuplot documentation for explanation.%
    }{The gnuplot epslatex terminal needs graphicx.sty or graphics.sty.}%
    \renewcommand\includegraphics[2][]{}%
  }%
  \providecommand\rotatebox[2]{#2}%
  \@ifundefined{ifGPcolor}{%
    \newif\ifGPcolor
    \GPcolortrue
  }{}%
  \@ifundefined{ifGPblacktext}{%
    \newif\ifGPblacktext
    \GPblacktexttrue
  }{}%
  \let\gplgaddtomacro\g@addto@macro
  \gdef\gplbacktext{}%
  \gdef\gplfronttext{}%
  \makeatother
  \ifGPblacktext
    \def\colorrgb#1{}%
    \def\colorgray#1{}%
  \else
    \ifGPcolor
      \def\colorrgb#1{\color[rgb]{#1}}%
      \def\colorgray#1{\color[gray]{#1}}%
      \expandafter\def\csname LTw\endcsname{\color{white}}%
      \expandafter\def\csname LTb\endcsname{\color{black}}%
      \expandafter\def\csname LTa\endcsname{\color{black}}%
      \expandafter\def\csname LT0\endcsname{\color[rgb]{1,0,0}}%
      \expandafter\def\csname LT1\endcsname{\color[rgb]{0,1,0}}%
      \expandafter\def\csname LT2\endcsname{\color[rgb]{0,0,1}}%
      \expandafter\def\csname LT3\endcsname{\color[rgb]{1,0,1}}%
      \expandafter\def\csname LT4\endcsname{\color[rgb]{0,1,1}}%
      \expandafter\def\csname LT5\endcsname{\color[rgb]{1,1,0}}%
      \expandafter\def\csname LT6\endcsname{\color[rgb]{0,0,0}}%
      \expandafter\def\csname LT7\endcsname{\color[rgb]{1,0.3,0}}%
      \expandafter\def\csname LT8\endcsname{\color[rgb]{0.5,0.5,0.5}}%
    \else
      \def\colorrgb#1{\color{black}}%
      \def\colorgray#1{\color[gray]{#1}}%
      \expandafter\def\csname LTw\endcsname{\color{white}}%
      \expandafter\def\csname LTb\endcsname{\color{black}}%
      \expandafter\def\csname LTa\endcsname{\color{black}}%
      \expandafter\def\csname LT0\endcsname{\color{black}}%
      \expandafter\def\csname LT1\endcsname{\color{black}}%
      \expandafter\def\csname LT2\endcsname{\color{black}}%
      \expandafter\def\csname LT3\endcsname{\color{black}}%
      \expandafter\def\csname LT4\endcsname{\color{black}}%
      \expandafter\def\csname LT5\endcsname{\color{black}}%
      \expandafter\def\csname LT6\endcsname{\color{black}}%
      \expandafter\def\csname LT7\endcsname{\color{black}}%
      \expandafter\def\csname LT8\endcsname{\color{black}}%
    \fi
  \fi
  \setlength{\unitlength}{0.0500bp}%
  \begin{picture}(7200.00,5040.00)%
    \gplgaddtomacro\gplbacktext{%
      \csname LTb\endcsname%
      \put(990,704){\makebox(0,0)[r]{\strut{}1e-08}}%
      \csname LTb\endcsname%
      \put(990,1112){\makebox(0,0)[r]{\strut{}1e-07}}%
      \csname LTb\endcsname%
      \put(990,1521){\makebox(0,0)[r]{\strut{}1e-06}}%
      \csname LTb\endcsname%
      \put(990,1929){\makebox(0,0)[r]{\strut{}1e-05}}%
      \csname LTb\endcsname%
      \put(990,2337){\makebox(0,0)[r]{\strut{}0.0001}}%
      \csname LTb\endcsname%
      \put(990,2746){\makebox(0,0)[r]{\strut{}0.001}}%
      \csname LTb\endcsname%
      \put(990,3154){\makebox(0,0)[r]{\strut{}0.01}}%
      \csname LTb\endcsname%
      \put(990,3562){\makebox(0,0)[r]{\strut{}0.1}}%
      \csname LTb\endcsname%
      \put(990,3971){\makebox(0,0)[r]{\strut{}1}}%
      \csname LTb\endcsname%
      \put(990,4379){\makebox(0,0)[r]{\strut{}10}}%
      \csname LTb\endcsname%
      \put(1122,484){\makebox(0,0){\strut{}100}}%
      \csname LTb\endcsname%
      \put(2146,484){\makebox(0,0){\strut{}1000}}%
      \csname LTb\endcsname%
      \put(3171,484){\makebox(0,0){\strut{}10000}}%
      \csname LTb\endcsname%
      \put(4195,484){\makebox(0,0){\strut{}100000}}%
      \csname LTb\endcsname%
      \put(5219,484){\makebox(0,0){\strut{}1e+06}}%
      \put(3170,154){\makebox(0,0){\strut{}\text{DOFs}}}%
      \put(3170,4709){\makebox(0,0){\strut{}Configuration 5}}%
    }%
    \gplgaddtomacro\gplfronttext{%
      \csname LTb\endcsname%
      \put(6212,4269){\makebox(0,0)[r]{\strut{} \footnotesize Error in $J_c$}}%
      \csname LTb\endcsname%
      \put(6212,4049){\makebox(0,0)[r]{\strut{} \footnotesize  $\eta_h$}}%
      \csname LTb\endcsname%
      \put(6212,3829){\makebox(0,0)[r]{\strut{} \footnotesize Error in $J_1$}}%
      \csname LTb\endcsname%
      \put(6212,3609){\makebox(0,0)[r]{\strut{} \footnotesize Error in $J_2$}}%
      \csname LTb\endcsname%
      \put(6212,3389){\makebox(0,0)[r]{\strut{} \footnotesize Error in $J_3$}}%
      \csname LTb\endcsname%
      \put(6212,3169){\makebox(0,0)[r]{\strut{} \footnotesize Error in $J_4$}}%
      \csname LTb\endcsname%
      \put(6212,2949){\makebox(0,0)[r]{\strut{} \footnotesize Error in $J_5$}}%
      \csname LTb\endcsname%
      \put(6212,2729){\makebox(0,0)[r]{\strut{} \footnotesize Error in $J_6$}}%
      \csname LTb\endcsname%
      \put(6212,2509){\makebox(0,0)[r]{\strut{} \footnotesize Error in $J_7$}}%
      \csname LTb\endcsname%
      \put(6212,2289){\makebox(0,0)[r]{\strut{}$\mathcal{O}(\text{DOFs}^-\frac{3}{2})$}}%
    }%
    \gplbacktext
    \put(0,0){\includegraphics{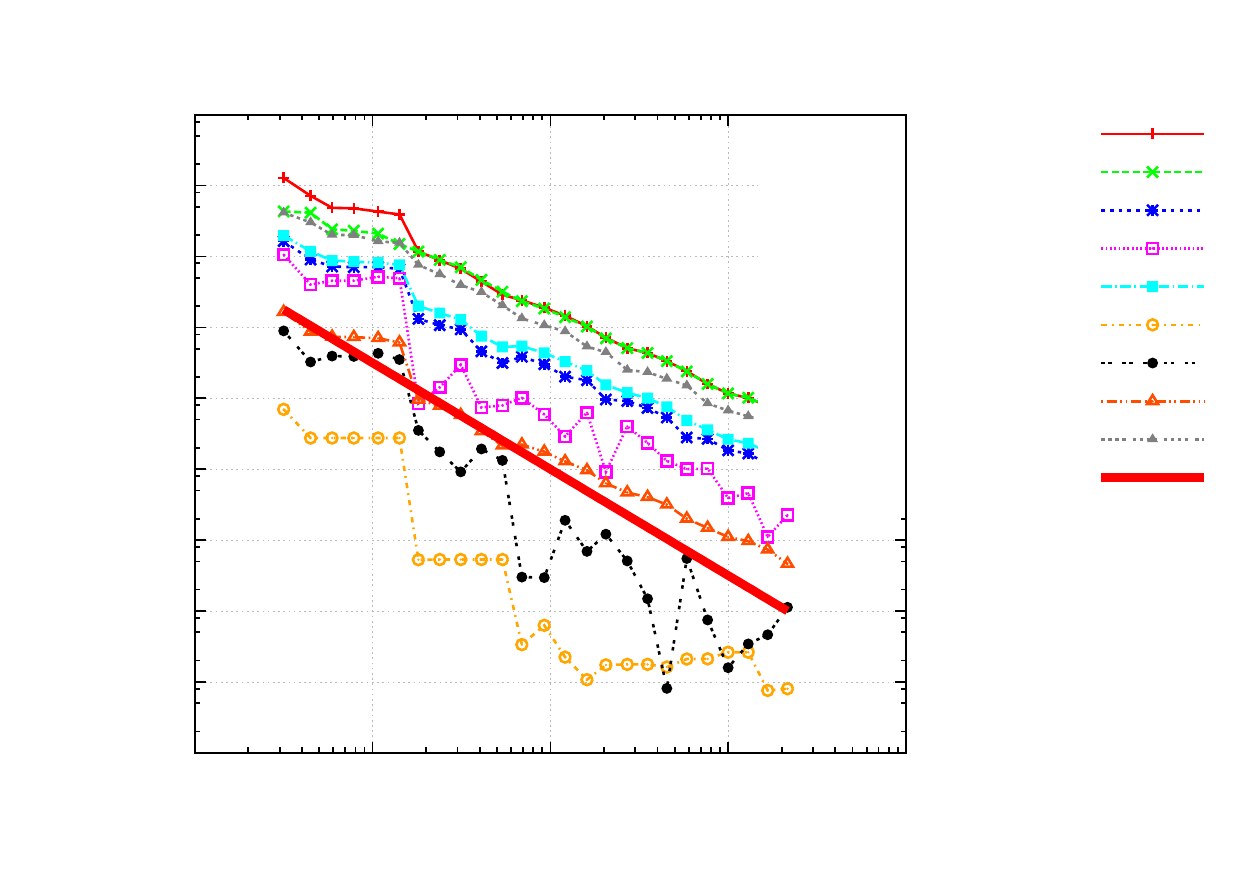}}%
    \gplfronttext
  \end{picture}%
\endgroup

%% file: Figures/YSplitterErrorConfig6.tex
\begingroup
  \makeatletter
  \providecommand\color[2][]{%
    \GenericError{(gnuplot) \space\space\space\@spaces}{%
      Package color not loaded in conjunction with
      terminal option `colourtext'%
    }{See the gnuplot documentation for explanation.%
    }{Either use 'blacktext' in gnuplot or load the package
      color.sty in LaTeX.}%
    \renewcommand\color[2][]{}%
  }%
  \providecommand\includegraphics[2][]{%
    \GenericError{(gnuplot) \space\space\space\@spaces}{%
      Package graphicx or graphics not loaded%
    }{See the gnuplot documentation for explanation.%
    }{The gnuplot epslatex terminal needs graphicx.sty or graphics.sty.}%
    \renewcommand\includegraphics[2][]{}%
  }%
  \providecommand\rotatebox[2]{#2}%
  \@ifundefined{ifGPcolor}{%
    \newif\ifGPcolor
    \GPcolortrue
  }{}%
  \@ifundefined{ifGPblacktext}{%
    \newif\ifGPblacktext
    \GPblacktexttrue
  }{}%
  \let\gplgaddtomacro\g@addto@macro
  \gdef\gplbacktext{}%
  \gdef\gplfronttext{}%
  \makeatother
  \ifGPblacktext
    \def\colorrgb#1{}%
    \def\colorgray#1{}%
  \else
    \ifGPcolor
      \def\colorrgb#1{\color[rgb]{#1}}%
      \def\colorgray#1{\color[gray]{#1}}%
      \expandafter\def\csname LTw\endcsname{\color{white}}%
      \expandafter\def\csname LTb\endcsname{\color{black}}%
      \expandafter\def\csname LTa\endcsname{\color{black}}%
      \expandafter\def\csname LT0\endcsname{\color[rgb]{1,0,0}}%
      \expandafter\def\csname LT1\endcsname{\color[rgb]{0,1,0}}%
      \expandafter\def\csname LT2\endcsname{\color[rgb]{0,0,1}}%
      \expandafter\def\csname LT3\endcsname{\color[rgb]{1,0,1}}%
      \expandafter\def\csname LT4\endcsname{\color[rgb]{0,1,1}}%
      \expandafter\def\csname LT5\endcsname{\color[rgb]{1,1,0}}%
      \expandafter\def\csname LT6\endcsname{\color[rgb]{0,0,0}}%
      \expandafter\def\csname LT7\endcsname{\color[rgb]{1,0.3,0}}%
      \expandafter\def\csname LT8\endcsname{\color[rgb]{0.5,0.5,0.5}}%
    \else
      \def\colorrgb#1{\color{black}}%
      \def\colorgray#1{\color[gray]{#1}}%
      \expandafter\def\csname LTw\endcsname{\color{white}}%
      \expandafter\def\csname LTb\endcsname{\color{black}}%
      \expandafter\def\csname LTa\endcsname{\color{black}}%
      \expandafter\def\csname LT0\endcsname{\color{black}}%
      \expandafter\def\csname LT1\endcsname{\color{black}}%
      \expandafter\def\csname LT2\endcsname{\color{black}}%
      \expandafter\def\csname LT3\endcsname{\color{black}}%
      \expandafter\def\csname LT4\endcsname{\color{black}}%
      \expandafter\def\csname LT5\endcsname{\color{black}}%
      \expandafter\def\csname LT6\endcsname{\color{black}}%
      \expandafter\def\csname LT7\endcsname{\color{black}}%
      \expandafter\def\csname LT8\endcsname{\color{black}}%
    \fi
  \fi
  \setlength{\unitlength}{0.0500bp}%
  \begin{picture}(7200.00,5040.00)%
    \gplgaddtomacro\gplbacktext{%
      \csname LTb\endcsname%
      \put(990,704){\makebox(0,0)[r]{\strut{}1e-08}}%
      \csname LTb\endcsname%
      \put(990,1163){\makebox(0,0)[r]{\strut{}1e-07}}%
      \csname LTb\endcsname%
      \put(990,1623){\makebox(0,0)[r]{\strut{}1e-06}}%
      \csname LTb\endcsname%
      \put(990,2082){\makebox(0,0)[r]{\strut{}1e-05}}%
      \csname LTb\endcsname%
      \put(990,2542){\makebox(0,0)[r]{\strut{}0.0001}}%
      \csname LTb\endcsname%
      \put(990,3001){\makebox(0,0)[r]{\strut{}0.001}}%
      \csname LTb\endcsname%
      \put(990,3460){\makebox(0,0)[r]{\strut{}0.01}}%
      \csname LTb\endcsname%
      \put(990,3920){\makebox(0,0)[r]{\strut{}0.1}}%
      \csname LTb\endcsname%
      \put(990,4379){\makebox(0,0)[r]{\strut{}1}}%
      \csname LTb\endcsname%
      \put(1122,484){\makebox(0,0){\strut{}100}}%
      \csname LTb\endcsname%
      \put(2146,484){\makebox(0,0){\strut{}1000}}%
      \csname LTb\endcsname%
      \put(3171,484){\makebox(0,0){\strut{}10000}}%
      \csname LTb\endcsname%
      \put(4195,484){\makebox(0,0){\strut{}100000}}%
      \csname LTb\endcsname%
      \put(5219,484){\makebox(0,0){\strut{}1e+06}}%
      \put(3170,154){\makebox(0,0){\strut{}\text{DOFs}}}%
      \put(3170,4709){\makebox(0,0){\strut{}Configuration 6}}%
    }%
    \gplgaddtomacro\gplfronttext{%
      \csname LTb\endcsname%
      \put(6212,4269){\makebox(0,0)[r]{\strut{} \footnotesize Error in $J_c$}}%
      \csname LTb\endcsname%
      \put(6212,4049){\makebox(0,0)[r]{\strut{} \footnotesize  $\eta_h$}}%
      \csname LTb\endcsname%
      \put(6212,3829){\makebox(0,0)[r]{\strut{} \footnotesize Error in $J_1$}}%
      \csname LTb\endcsname%
      \put(6212,3609){\makebox(0,0)[r]{\strut{} \footnotesize Error in $J_2$}}%
      \csname LTb\endcsname%
      \put(6212,3389){\makebox(0,0)[r]{\strut{} \footnotesize Error in $J_3$}}%
      \csname LTb\endcsname%
      \put(6212,3169){\makebox(0,0)[r]{\strut{} \footnotesize Error in $J_4$}}%
      \csname LTb\endcsname%
      \put(6212,2949){\makebox(0,0)[r]{\strut{} \footnotesize Error in $J_5$}}%
      \csname LTb\endcsname%
      \put(6212,2729){\makebox(0,0)[r]{\strut{} \footnotesize Error in $J_6$}}%
      \csname LTb\endcsname%
      \put(6212,2509){\makebox(0,0)[r]{\strut{} \footnotesize Error in $J_7$}}%
      \csname LTb\endcsname%
      \put(6212,2289){\makebox(0,0)[r]{\strut{}$\mathcal{O}(\text{DOFs}^-\frac{3}{2})$}}%
    }%
    \gplbacktext
    \put(0,0){\includegraphics{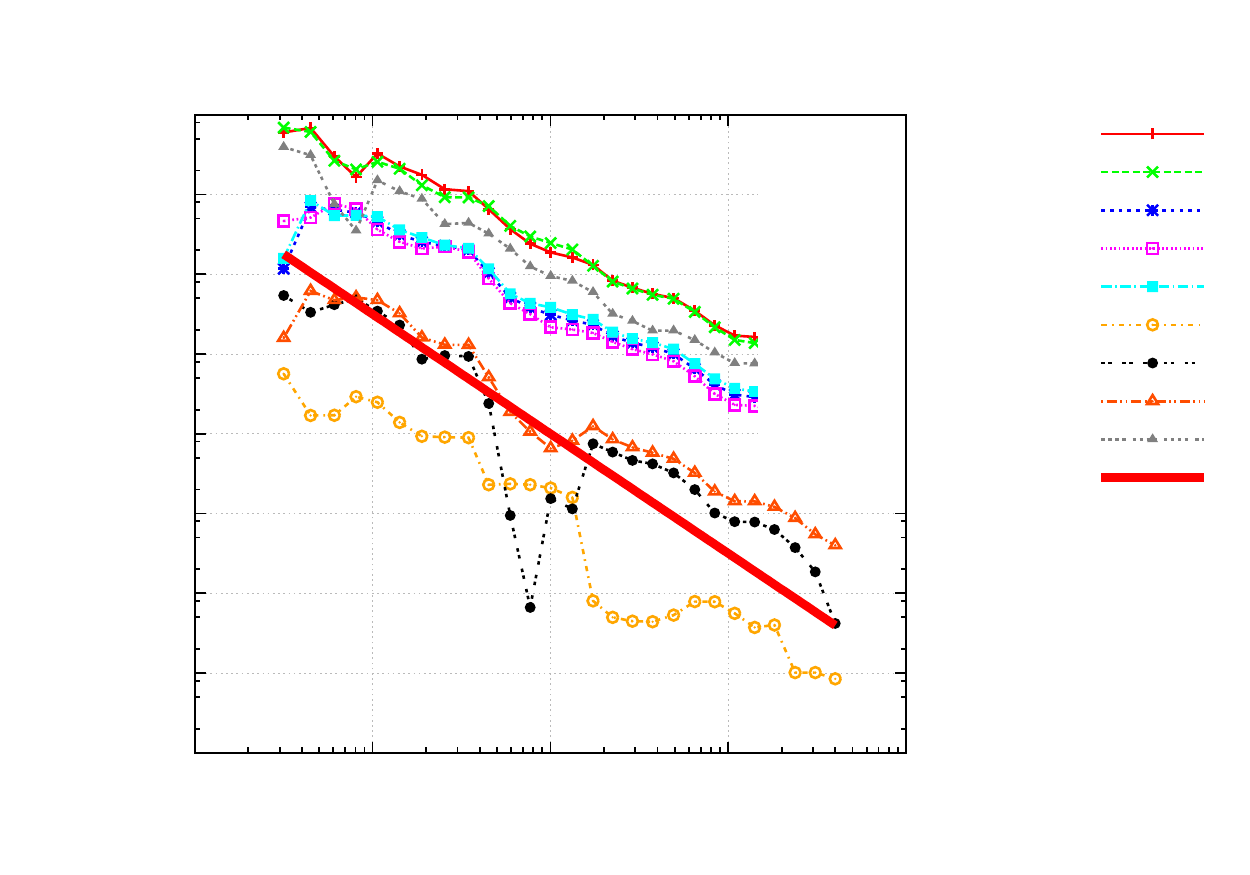}}%
    \gplfronttext
  \end{picture}%
\endgroup

%% file: Figures/YSplitterIeffConfig1.tex
\begingroup
  \makeatletter
  \providecommand\color[2][]{%
    \GenericError{(gnuplot) \space\space\space\@spaces}{%
      Package color not loaded in conjunction with
      terminal option `colourtext'%
    }{See the gnuplot documentation for explanation.%
    }{Either use 'blacktext' in gnuplot or load the package
      color.sty in LaTeX.}%
    \renewcommand\color[2][]{}%
  }%
  \providecommand\includegraphics[2][]{%
    \GenericError{(gnuplot) \space\space\space\@spaces}{%
      Package graphicx or graphics not loaded%
    }{See the gnuplot documentation for explanation.%
    }{The gnuplot epslatex terminal needs graphicx.sty or graphics.sty.}%
    \renewcommand\includegraphics[2][]{}%
  }%
  \providecommand\rotatebox[2]{#2}%
  \@ifundefined{ifGPcolor}{%
    \newif\ifGPcolor
    \GPcolortrue
  }{}%
  \@ifundefined{ifGPblacktext}{%
    \newif\ifGPblacktext
    \GPblacktexttrue
  }{}%
  \let\gplgaddtomacro\g@addto@macro
  \gdef\gplbacktext{}%
  \gdef\gplfronttext{}%
  \makeatother
  \ifGPblacktext
    \def\colorrgb#1{}%
    \def\colorgray#1{}%
  \else
    \ifGPcolor
      \def\colorrgb#1{\color[rgb]{#1}}%
      \def\colorgray#1{\color[gray]{#1}}%
      \expandafter\def\csname LTw\endcsname{\color{white}}%
      \expandafter\def\csname LTb\endcsname{\color{black}}%
      \expandafter\def\csname LTa\endcsname{\color{black}}%
      \expandafter\def\csname LT0\endcsname{\color[rgb]{1,0,0}}%
      \expandafter\def\csname LT1\endcsname{\color[rgb]{0,1,0}}%
      \expandafter\def\csname LT2\endcsname{\color[rgb]{0,0,1}}%
      \expandafter\def\csname LT3\endcsname{\color[rgb]{1,0,1}}%
      \expandafter\def\csname LT4\endcsname{\color[rgb]{0,1,1}}%
      \expandafter\def\csname LT5\endcsname{\color[rgb]{1,1,0}}%
      \expandafter\def\csname LT6\endcsname{\color[rgb]{0,0,0}}%
      \expandafter\def\csname LT7\endcsname{\color[rgb]{1,0.3,0}}%
      \expandafter\def\csname LT8\endcsname{\color[rgb]{0.5,0.5,0.5}}%
    \else
      \def\colorrgb#1{\color{black}}%
      \def\colorgray#1{\color[gray]{#1}}%
      \expandafter\def\csname LTw\endcsname{\color{white}}%
      \expandafter\def\csname LTb\endcsname{\color{black}}%
      \expandafter\def\csname LTa\endcsname{\color{black}}%
      \expandafter\def\csname LT0\endcsname{\color{black}}%
      \expandafter\def\csname LT1\endcsname{\color{black}}%
      \expandafter\def\csname LT2\endcsname{\color{black}}%
      \expandafter\def\csname LT3\endcsname{\color{black}}%
      \expandafter\def\csname LT4\endcsname{\color{black}}%
      \expandafter\def\csname LT5\endcsname{\color{black}}%
      \expandafter\def\csname LT6\endcsname{\color{black}}%
      \expandafter\def\csname LT7\endcsname{\color{black}}%
      \expandafter\def\csname LT8\endcsname{\color{black}}%
    \fi
  \fi
  \setlength{\unitlength}{0.0500bp}%
  \begin{picture}(7200.00,5040.00)%
    \gplgaddtomacro\gplbacktext{%
      \csname LTb\endcsname%
      \put(594,704){\makebox(0,0)[r]{\strut{}0.1}}%
      \csname LTb\endcsname%
      \put(594,2542){\makebox(0,0)[r]{\strut{}1}}%
      \csname LTb\endcsname%
      \put(594,4379){\makebox(0,0)[r]{\strut{}10}}%
      \csname LTb\endcsname%
      \put(726,484){\makebox(0,0){\strut{}100}}%
      \csname LTb\endcsname%
      \put(1849,484){\makebox(0,0){\strut{}1000}}%
      \csname LTb\endcsname%
      \put(2973,484){\makebox(0,0){\strut{}10000}}%
      \csname LTb\endcsname%
      \put(4096,484){\makebox(0,0){\strut{}100000}}%
      \csname LTb\endcsname%
      \put(5219,484){\makebox(0,0){\strut{}1e+06}}%
      \put(2972,154){\makebox(0,0){\strut{}\text{DOFs}}}%
      \put(2972,4709){\makebox(0,0){\strut{}Configuration 1}}%
    }%
    \gplgaddtomacro\gplfronttext{%
      \csname LTb\endcsname%
      \put(6212,4269){\makebox(0,0)[r]{\strut{} \footnotesize $I_{eff}$}}%
      \csname LTb\endcsname%
      \put(6212,4049){\makebox(0,0)[r]{\strut{} \footnotesize  $I_{eff,a}$}}%
      \csname LTb\endcsname%
      \put(6212,3829){\makebox(0,0)[r]{\strut{} \footnotesize $I_{eff,p}$}}%
      \csname LTb\endcsname%
      \put(6212,3609){\makebox(0,0)[r]{\strut{}1}}%
    }%
    \gplbacktext
    \put(0,0){\includegraphics{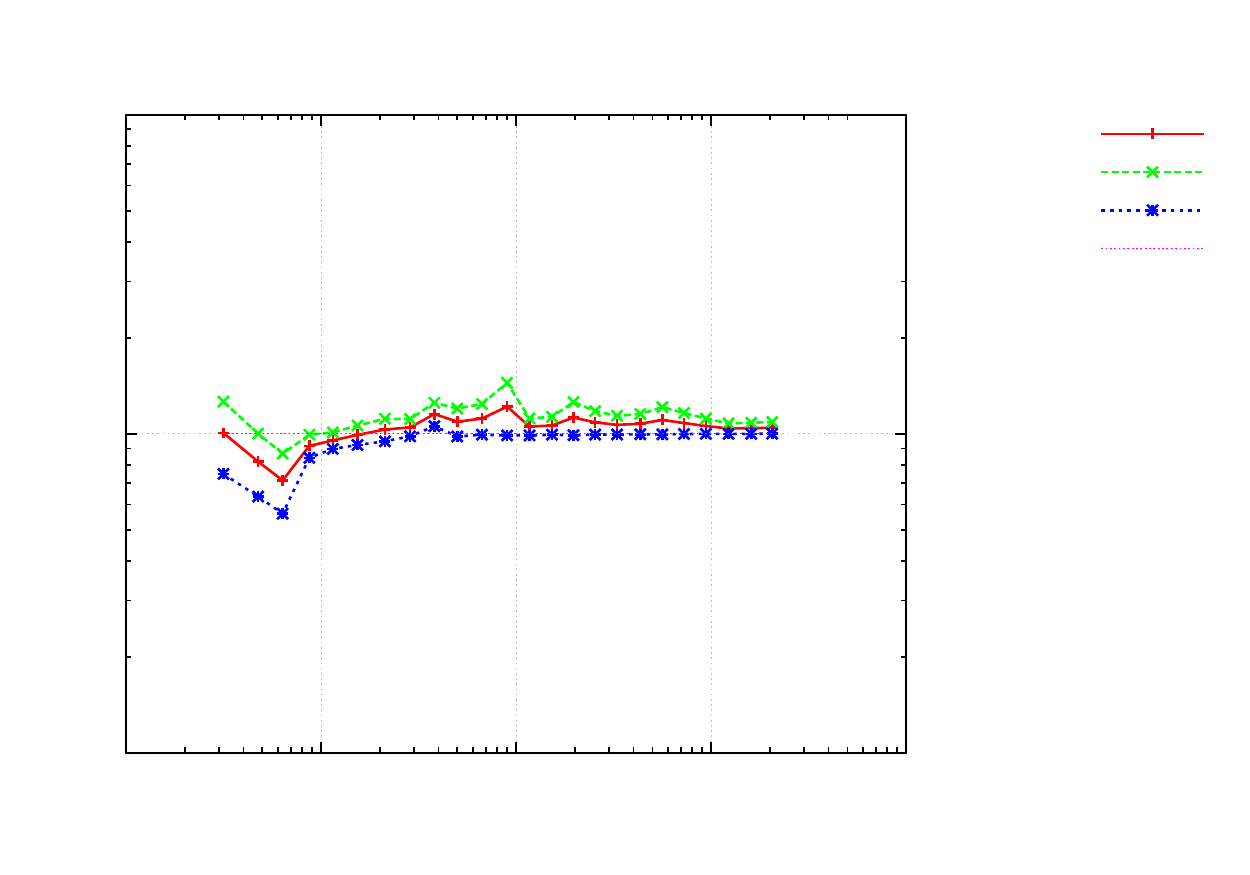}}%
    \gplfronttext
  \end{picture}%
\endgroup

%% file: Figures/YSplitterIeffConfig2.tex
\begingroup
  \makeatletter
  \providecommand\color[2][]{%
    \GenericError{(gnuplot) \space\space\space\@spaces}{%
      Package color not loaded in conjunction with
      terminal option `colourtext'%
    }{See the gnuplot documentation for explanation.%
    }{Either use 'blacktext' in gnuplot or load the package
      color.sty in LaTeX.}%
    \renewcommand\color[2][]{}%
  }%
  \providecommand\includegraphics[2][]{%
    \GenericError{(gnuplot) \space\space\space\@spaces}{%
      Package graphicx or graphics not loaded%
    }{See the gnuplot documentation for explanation.%
    }{The gnuplot epslatex terminal needs graphicx.sty or graphics.sty.}%
    \renewcommand\includegraphics[2][]{}%
  }%
  \providecommand\rotatebox[2]{#2}%
  \@ifundefined{ifGPcolor}{%
    \newif\ifGPcolor
    \GPcolortrue
  }{}%
  \@ifundefined{ifGPblacktext}{%
    \newif\ifGPblacktext
    \GPblacktexttrue
  }{}%
  \let\gplgaddtomacro\g@addto@macro
  \gdef\gplbacktext{}%
  \gdef\gplfronttext{}%
  \makeatother
  \ifGPblacktext
    \def\colorrgb#1{}%
    \def\colorgray#1{}%
  \else
    \ifGPcolor
      \def\colorrgb#1{\color[rgb]{#1}}%
      \def\colorgray#1{\color[gray]{#1}}%
      \expandafter\def\csname LTw\endcsname{\color{white}}%
      \expandafter\def\csname LTb\endcsname{\color{black}}%
      \expandafter\def\csname LTa\endcsname{\color{black}}%
      \expandafter\def\csname LT0\endcsname{\color[rgb]{1,0,0}}%
      \expandafter\def\csname LT1\endcsname{\color[rgb]{0,1,0}}%
      \expandafter\def\csname LT2\endcsname{\color[rgb]{0,0,1}}%
      \expandafter\def\csname LT3\endcsname{\color[rgb]{1,0,1}}%
      \expandafter\def\csname LT4\endcsname{\color[rgb]{0,1,1}}%
      \expandafter\def\csname LT5\endcsname{\color[rgb]{1,1,0}}%
      \expandafter\def\csname LT6\endcsname{\color[rgb]{0,0,0}}%
      \expandafter\def\csname LT7\endcsname{\color[rgb]{1,0.3,0}}%
      \expandafter\def\csname LT8\endcsname{\color[rgb]{0.5,0.5,0.5}}%
    \else
      \def\colorrgb#1{\color{black}}%
      \def\colorgray#1{\color[gray]{#1}}%
      \expandafter\def\csname LTw\endcsname{\color{white}}%
      \expandafter\def\csname LTb\endcsname{\color{black}}%
      \expandafter\def\csname LTa\endcsname{\color{black}}%
      \expandafter\def\csname LT0\endcsname{\color{black}}%
      \expandafter\def\csname LT1\endcsname{\color{black}}%
      \expandafter\def\csname LT2\endcsname{\color{black}}%
      \expandafter\def\csname LT3\endcsname{\color{black}}%
      \expandafter\def\csname LT4\endcsname{\color{black}}%
      \expandafter\def\csname LT5\endcsname{\color{black}}%
      \expandafter\def\csname LT6\endcsname{\color{black}}%
      \expandafter\def\csname LT7\endcsname{\color{black}}%
      \expandafter\def\csname LT8\endcsname{\color{black}}%
    \fi
  \fi
  \setlength{\unitlength}{0.0500bp}%
  \begin{picture}(7200.00,5040.00)%
    \gplgaddtomacro\gplbacktext{%
      \csname LTb\endcsname%
      \put(726,704){\makebox(0,0)[r]{\strut{}0.01}}%
      \csname LTb\endcsname%
      \put(726,1929){\makebox(0,0)[r]{\strut{}0.1}}%
      \csname LTb\endcsname%
      \put(726,3154){\makebox(0,0)[r]{\strut{}1}}%
      \csname LTb\endcsname%
      \put(726,4379){\makebox(0,0)[r]{\strut{}10}}%
      \csname LTb\endcsname%
      \put(858,484){\makebox(0,0){\strut{}100}}%
      \csname LTb\endcsname%
      \put(1948,484){\makebox(0,0){\strut{}1000}}%
      \csname LTb\endcsname%
      \put(3039,484){\makebox(0,0){\strut{}10000}}%
      \csname LTb\endcsname%
      \put(4129,484){\makebox(0,0){\strut{}100000}}%
      \csname LTb\endcsname%
      \put(5219,484){\makebox(0,0){\strut{}1e+06}}%
      \put(3038,154){\makebox(0,0){\strut{}\text{DOFs}}}%
      \put(3038,4709){\makebox(0,0){\strut{}Configuration 2}}%
    }%
    \gplgaddtomacro\gplfronttext{%
      \csname LTb\endcsname%
      \put(6212,4269){\makebox(0,0)[r]{\strut{} \footnotesize $I_{eff}$}}%
      \csname LTb\endcsname%
      \put(6212,4049){\makebox(0,0)[r]{\strut{} \footnotesize  $I_{eff,a}$}}%
      \csname LTb\endcsname%
      \put(6212,3829){\makebox(0,0)[r]{\strut{} \footnotesize $I_{eff,p}$}}%
      \csname LTb\endcsname%
      \put(6212,3609){\makebox(0,0)[r]{\strut{}1}}%
    }%
    \gplbacktext
    \put(0,0){\includegraphics{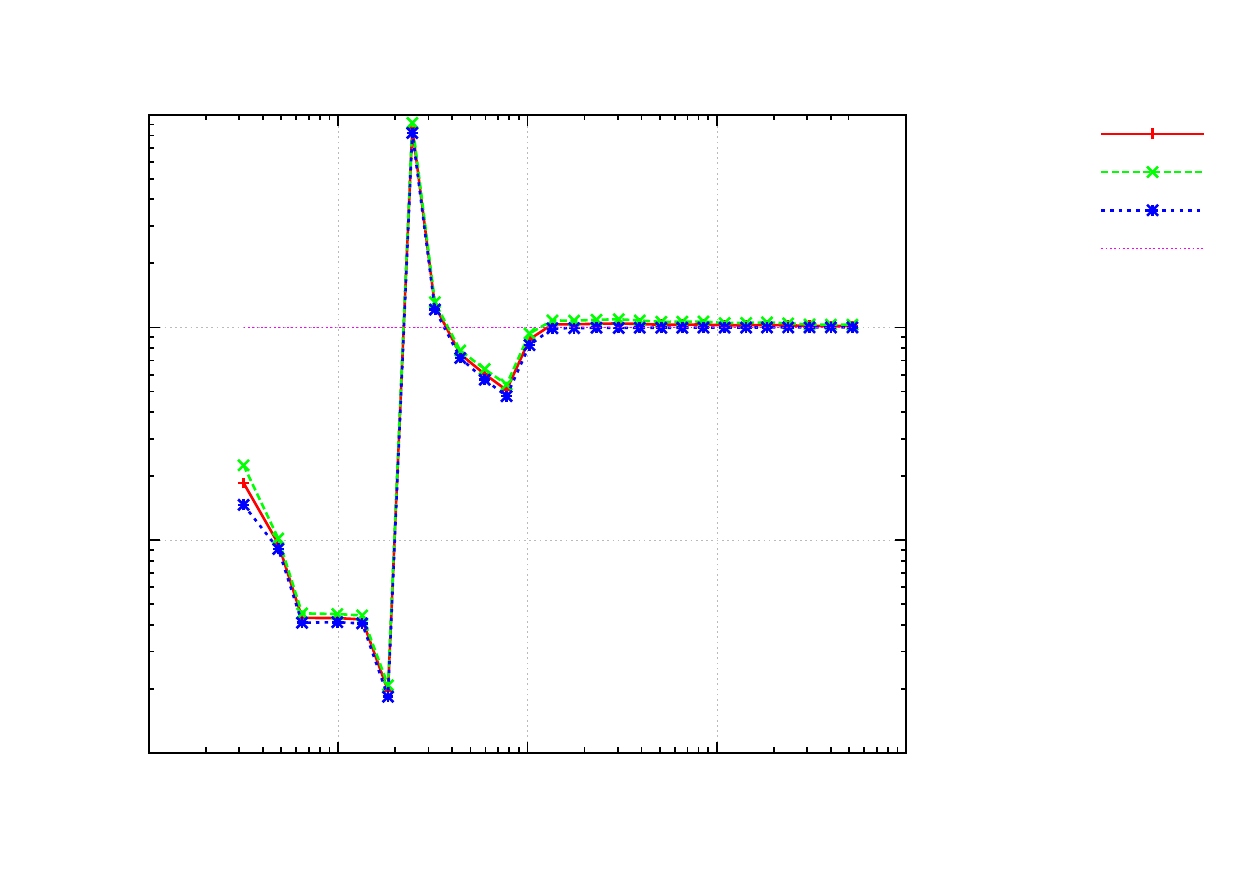}}%
    \gplfronttext
  \end{picture}%
\endgroup

%% file: Figures/YSplitterIeffConfig3.tex
\begingroup
  \makeatletter
  \providecommand\color[2][]{%
    \GenericError{(gnuplot) \space\space\space\@spaces}{%
      Package color not loaded in conjunction with
      terminal option `colourtext'%
    }{See the gnuplot documentation for explanation.%
    }{Either use 'blacktext' in gnuplot or load the package
      color.sty in LaTeX.}%
    \renewcommand\color[2][]{}%
  }%
  \providecommand\includegraphics[2][]{%
    \GenericError{(gnuplot) \space\space\space\@spaces}{%
      Package graphicx or graphics not loaded%
    }{See the gnuplot documentation for explanation.%
    }{The gnuplot epslatex terminal needs graphicx.sty or graphics.sty.}%
    \renewcommand\includegraphics[2][]{}%
  }%
  \providecommand\rotatebox[2]{#2}%
  \@ifundefined{ifGPcolor}{%
    \newif\ifGPcolor
    \GPcolortrue
  }{}%
  \@ifundefined{ifGPblacktext}{%
    \newif\ifGPblacktext
    \GPblacktexttrue
  }{}%
  \let\gplgaddtomacro\g@addto@macro
  \gdef\gplbacktext{}%
  \gdef\gplfronttext{}%
  \makeatother
  \ifGPblacktext
    \def\colorrgb#1{}%
    \def\colorgray#1{}%
  \else
    \ifGPcolor
      \def\colorrgb#1{\color[rgb]{#1}}%
      \def\colorgray#1{\color[gray]{#1}}%
      \expandafter\def\csname LTw\endcsname{\color{white}}%
      \expandafter\def\csname LTb\endcsname{\color{black}}%
      \expandafter\def\csname LTa\endcsname{\color{black}}%
      \expandafter\def\csname LT0\endcsname{\color[rgb]{1,0,0}}%
      \expandafter\def\csname LT1\endcsname{\color[rgb]{0,1,0}}%
      \expandafter\def\csname LT2\endcsname{\color[rgb]{0,0,1}}%
      \expandafter\def\csname LT3\endcsname{\color[rgb]{1,0,1}}%
      \expandafter\def\csname LT4\endcsname{\color[rgb]{0,1,1}}%
      \expandafter\def\csname LT5\endcsname{\color[rgb]{1,1,0}}%
      \expandafter\def\csname LT6\endcsname{\color[rgb]{0,0,0}}%
      \expandafter\def\csname LT7\endcsname{\color[rgb]{1,0.3,0}}%
      \expandafter\def\csname LT8\endcsname{\color[rgb]{0.5,0.5,0.5}}%
    \else
      \def\colorrgb#1{\color{black}}%
      \def\colorgray#1{\color[gray]{#1}}%
      \expandafter\def\csname LTw\endcsname{\color{white}}%
      \expandafter\def\csname LTb\endcsname{\color{black}}%
      \expandafter\def\csname LTa\endcsname{\color{black}}%
      \expandafter\def\csname LT0\endcsname{\color{black}}%
      \expandafter\def\csname LT1\endcsname{\color{black}}%
      \expandafter\def\csname LT2\endcsname{\color{black}}%
      \expandafter\def\csname LT3\endcsname{\color{black}}%
      \expandafter\def\csname LT4\endcsname{\color{black}}%
      \expandafter\def\csname LT5\endcsname{\color{black}}%
      \expandafter\def\csname LT6\endcsname{\color{black}}%
      \expandafter\def\csname LT7\endcsname{\color{black}}%
      \expandafter\def\csname LT8\endcsname{\color{black}}%
    \fi
  \fi
  \setlength{\unitlength}{0.0500bp}%
  \begin{picture}(7200.00,5040.00)%
    \gplgaddtomacro\gplbacktext{%
      \csname LTb\endcsname%
      \put(990,704){\makebox(0,0)[r]{\strut{}0.0001}}%
      \csname LTb\endcsname%
      \put(990,1439){\makebox(0,0)[r]{\strut{}0.001}}%
      \csname LTb\endcsname%
      \put(990,2174){\makebox(0,0)[r]{\strut{}0.01}}%
      \csname LTb\endcsname%
      \put(990,2909){\makebox(0,0)[r]{\strut{}0.1}}%
      \csname LTb\endcsname%
      \put(990,3644){\makebox(0,0)[r]{\strut{}1}}%
      \csname LTb\endcsname%
      \put(990,4379){\makebox(0,0)[r]{\strut{}10}}%
      \csname LTb\endcsname%
      \put(1122,484){\makebox(0,0){\strut{}100}}%
      \csname LTb\endcsname%
      \put(2146,484){\makebox(0,0){\strut{}1000}}%
      \csname LTb\endcsname%
      \put(3171,484){\makebox(0,0){\strut{}10000}}%
      \csname LTb\endcsname%
      \put(4195,484){\makebox(0,0){\strut{}100000}}%
      \csname LTb\endcsname%
      \put(5219,484){\makebox(0,0){\strut{}1e+06}}%
      \put(3170,154){\makebox(0,0){\strut{}\text{DOFs}}}%
      \put(3170,4709){\makebox(0,0){\strut{}Configuration 3}}%
    }%
    \gplgaddtomacro\gplfronttext{%
      \csname LTb\endcsname%
      \put(6212,4269){\makebox(0,0)[r]{\strut{} \footnotesize $I_{eff}$}}%
      \csname LTb\endcsname%
      \put(6212,4049){\makebox(0,0)[r]{\strut{} \footnotesize  $I_{eff,a}$}}%
      \csname LTb\endcsname%
      \put(6212,3829){\makebox(0,0)[r]{\strut{} \footnotesize $I_{eff,p}$}}%
      \csname LTb\endcsname%
      \put(6212,3609){\makebox(0,0)[r]{\strut{}1}}%
    }%
    \gplbacktext
    \put(0,0){\includegraphics{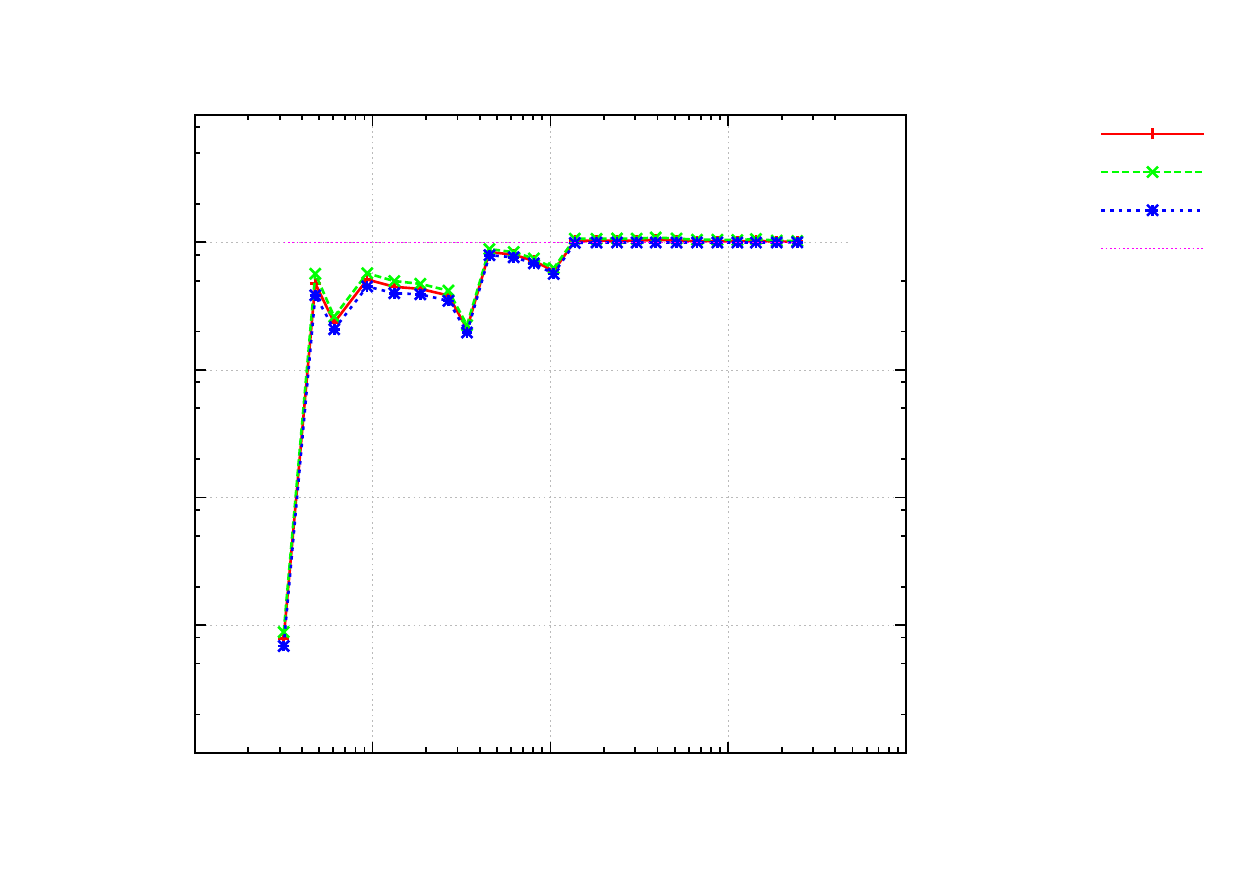}}%
    \gplfronttext
  \end{picture}%
\endgroup

%% file: Figures/YSplitterIeffConfig4.tex
\begingroup
  \makeatletter
  \providecommand\color[2][]{%
    \GenericError{(gnuplot) \space\space\space\@spaces}{%
      Package color not loaded in conjunction with
      terminal option `colourtext'%
    }{See the gnuplot documentation for explanation.%
    }{Either use 'blacktext' in gnuplot or load the package
      color.sty in LaTeX.}%
    \renewcommand\color[2][]{}%
  }%
  \providecommand\includegraphics[2][]{%
    \GenericError{(gnuplot) \space\space\space\@spaces}{%
      Package graphicx or graphics not loaded%
    }{See the gnuplot documentation for explanation.%
    }{The gnuplot epslatex terminal needs graphicx.sty or graphics.sty.}%
    \renewcommand\includegraphics[2][]{}%
  }%
  \providecommand\rotatebox[2]{#2}%
  \@ifundefined{ifGPcolor}{%
    \newif\ifGPcolor
    \GPcolortrue
  }{}%
  \@ifundefined{ifGPblacktext}{%
    \newif\ifGPblacktext
    \GPblacktexttrue
  }{}%
  \let\gplgaddtomacro\g@addto@macro
  \gdef\gplbacktext{}%
  \gdef\gplfronttext{}%
  \makeatother
  \ifGPblacktext
    \def\colorrgb#1{}%
    \def\colorgray#1{}%
  \else
    \ifGPcolor
      \def\colorrgb#1{\color[rgb]{#1}}%
      \def\colorgray#1{\color[gray]{#1}}%
      \expandafter\def\csname LTw\endcsname{\color{white}}%
      \expandafter\def\csname LTb\endcsname{\color{black}}%
      \expandafter\def\csname LTa\endcsname{\color{black}}%
      \expandafter\def\csname LT0\endcsname{\color[rgb]{1,0,0}}%
      \expandafter\def\csname LT1\endcsname{\color[rgb]{0,1,0}}%
      \expandafter\def\csname LT2\endcsname{\color[rgb]{0,0,1}}%
      \expandafter\def\csname LT3\endcsname{\color[rgb]{1,0,1}}%
      \expandafter\def\csname LT4\endcsname{\color[rgb]{0,1,1}}%
      \expandafter\def\csname LT5\endcsname{\color[rgb]{1,1,0}}%
      \expandafter\def\csname LT6\endcsname{\color[rgb]{0,0,0}}%
      \expandafter\def\csname LT7\endcsname{\color[rgb]{1,0.3,0}}%
      \expandafter\def\csname LT8\endcsname{\color[rgb]{0.5,0.5,0.5}}%
    \else
      \def\colorrgb#1{\color{black}}%
      \def\colorgray#1{\color[gray]{#1}}%
      \expandafter\def\csname LTw\endcsname{\color{white}}%
      \expandafter\def\csname LTb\endcsname{\color{black}}%
      \expandafter\def\csname LTa\endcsname{\color{black}}%
      \expandafter\def\csname LT0\endcsname{\color{black}}%
      \expandafter\def\csname LT1\endcsname{\color{black}}%
      \expandafter\def\csname LT2\endcsname{\color{black}}%
      \expandafter\def\csname LT3\endcsname{\color{black}}%
      \expandafter\def\csname LT4\endcsname{\color{black}}%
      \expandafter\def\csname LT5\endcsname{\color{black}}%
      \expandafter\def\csname LT6\endcsname{\color{black}}%
      \expandafter\def\csname LT7\endcsname{\color{black}}%
      \expandafter\def\csname LT8\endcsname{\color{black}}%
    \fi
  \fi
  \setlength{\unitlength}{0.0500bp}%
  \begin{picture}(7200.00,5040.00)%
    \gplgaddtomacro\gplbacktext{%
      \csname LTb\endcsname%
      \put(990,704){\makebox(0,0)[r]{\strut{}0.0001}}%
      \csname LTb\endcsname%
      \put(990,1439){\makebox(0,0)[r]{\strut{}0.001}}%
      \csname LTb\endcsname%
      \put(990,2174){\makebox(0,0)[r]{\strut{}0.01}}%
      \csname LTb\endcsname%
      \put(990,2909){\makebox(0,0)[r]{\strut{}0.1}}%
      \csname LTb\endcsname%
      \put(990,3644){\makebox(0,0)[r]{\strut{}1}}%
      \csname LTb\endcsname%
      \put(990,4379){\makebox(0,0)[r]{\strut{}10}}%
      \csname LTb\endcsname%
      \put(1122,484){\makebox(0,0){\strut{}100}}%
      \csname LTb\endcsname%
      \put(2146,484){\makebox(0,0){\strut{}1000}}%
      \csname LTb\endcsname%
      \put(3171,484){\makebox(0,0){\strut{}10000}}%
      \csname LTb\endcsname%
      \put(4195,484){\makebox(0,0){\strut{}100000}}%
      \csname LTb\endcsname%
      \put(5219,484){\makebox(0,0){\strut{}1e+06}}%
      \put(3170,154){\makebox(0,0){\strut{}\text{DOFs}}}%
      \put(3170,4709){\makebox(0,0){\strut{}Configuration 4}}%
    }%
    \gplgaddtomacro\gplfronttext{%
      \csname LTb\endcsname%
      \put(6212,4269){\makebox(0,0)[r]{\strut{} \footnotesize $I_{eff}$}}%
      \csname LTb\endcsname%
      \put(6212,4049){\makebox(0,0)[r]{\strut{} \footnotesize  $I_{eff,a}$}}%
      \csname LTb\endcsname%
      \put(6212,3829){\makebox(0,0)[r]{\strut{} \footnotesize $I_{eff,p}$}}%
      \csname LTb\endcsname%
      \put(6212,3609){\makebox(0,0)[r]{\strut{}1}}%
    }%
    \gplbacktext
    \put(0,0){\includegraphics{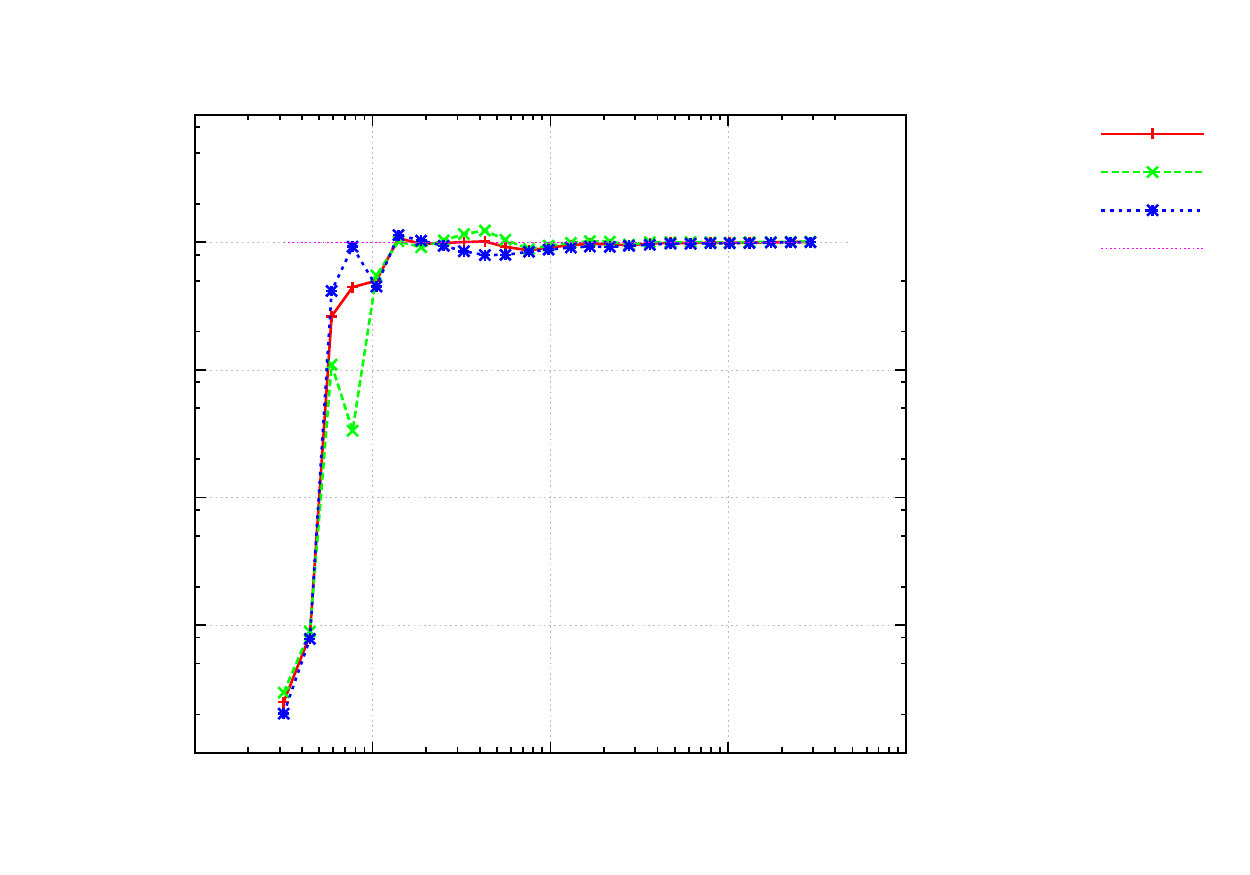}}%
    \gplfronttext
  \end{picture}%
\endgroup

%% file: Figures/YSplitterIeffConfig5.tex
\begingroup
  \makeatletter
  \providecommand\color[2][]{%
    \GenericError{(gnuplot) \space\space\space\@spaces}{%
      Package color not loaded in conjunction with
      terminal option `colourtext'%
    }{See the gnuplot documentation for explanation.%
    }{Either use 'blacktext' in gnuplot or load the package
      color.sty in LaTeX.}%
    \renewcommand\color[2][]{}%
  }%
  \providecommand\includegraphics[2][]{%
    \GenericError{(gnuplot) \space\space\space\@spaces}{%
      Package graphicx or graphics not loaded%
    }{See the gnuplot documentation for explanation.%
    }{The gnuplot epslatex terminal needs graphicx.sty or graphics.sty.}%
    \renewcommand\includegraphics[2][]{}%
  }%
  \providecommand\rotatebox[2]{#2}%
  \@ifundefined{ifGPcolor}{%
    \newif\ifGPcolor
    \GPcolortrue
  }{}%
  \@ifundefined{ifGPblacktext}{%
    \newif\ifGPblacktext
    \GPblacktexttrue
  }{}%
  \let\gplgaddtomacro\g@addto@macro
  \gdef\gplbacktext{}%
  \gdef\gplfronttext{}%
  \makeatother
  \ifGPblacktext
    \def\colorrgb#1{}%
    \def\colorgray#1{}%
  \else
    \ifGPcolor
      \def\colorrgb#1{\color[rgb]{#1}}%
      \def\colorgray#1{\color[gray]{#1}}%
      \expandafter\def\csname LTw\endcsname{\color{white}}%
      \expandafter\def\csname LTb\endcsname{\color{black}}%
      \expandafter\def\csname LTa\endcsname{\color{black}}%
      \expandafter\def\csname LT0\endcsname{\color[rgb]{1,0,0}}%
      \expandafter\def\csname LT1\endcsname{\color[rgb]{0,1,0}}%
      \expandafter\def\csname LT2\endcsname{\color[rgb]{0,0,1}}%
      \expandafter\def\csname LT3\endcsname{\color[rgb]{1,0,1}}%
      \expandafter\def\csname LT4\endcsname{\color[rgb]{0,1,1}}%
      \expandafter\def\csname LT5\endcsname{\color[rgb]{1,1,0}}%
      \expandafter\def\csname LT6\endcsname{\color[rgb]{0,0,0}}%
      \expandafter\def\csname LT7\endcsname{\color[rgb]{1,0.3,0}}%
      \expandafter\def\csname LT8\endcsname{\color[rgb]{0.5,0.5,0.5}}%
    \else
      \def\colorrgb#1{\color{black}}%
      \def\colorgray#1{\color[gray]{#1}}%
      \expandafter\def\csname LTw\endcsname{\color{white}}%
      \expandafter\def\csname LTb\endcsname{\color{black}}%
      \expandafter\def\csname LTa\endcsname{\color{black}}%
      \expandafter\def\csname LT0\endcsname{\color{black}}%
      \expandafter\def\csname LT1\endcsname{\color{black}}%
      \expandafter\def\csname LT2\endcsname{\color{black}}%
      \expandafter\def\csname LT3\endcsname{\color{black}}%
      \expandafter\def\csname LT4\endcsname{\color{black}}%
      \expandafter\def\csname LT5\endcsname{\color{black}}%
      \expandafter\def\csname LT6\endcsname{\color{black}}%
      \expandafter\def\csname LT7\endcsname{\color{black}}%
      \expandafter\def\csname LT8\endcsname{\color{black}}%
    \fi
  \fi
  \setlength{\unitlength}{0.0500bp}%
  \begin{picture}(7200.00,5040.00)%
    \gplgaddtomacro\gplbacktext{%
      \csname LTb\endcsname%
      \put(594,704){\makebox(0,0)[r]{\strut{}0.1}}%
      \csname LTb\endcsname%
      \put(594,2542){\makebox(0,0)[r]{\strut{}1}}%
      \csname LTb\endcsname%
      \put(594,4379){\makebox(0,0)[r]{\strut{}10}}%
      \csname LTb\endcsname%
      \put(726,484){\makebox(0,0){\strut{}100}}%
      \csname LTb\endcsname%
      \put(1849,484){\makebox(0,0){\strut{}1000}}%
      \csname LTb\endcsname%
      \put(2973,484){\makebox(0,0){\strut{}10000}}%
      \csname LTb\endcsname%
      \put(4096,484){\makebox(0,0){\strut{}100000}}%
      \csname LTb\endcsname%
      \put(5219,484){\makebox(0,0){\strut{}1e+06}}%
      \put(2972,154){\makebox(0,0){\strut{}\text{DOFs}}}%
      \put(2972,4709){\makebox(0,0){\strut{}Configuration 5}}%
    }%
    \gplgaddtomacro\gplfronttext{%
      \csname LTb\endcsname%
      \put(6212,4269){\makebox(0,0)[r]{\strut{} \footnotesize $I_{eff}$}}%
      \csname LTb\endcsname%
      \put(6212,4049){\makebox(0,0)[r]{\strut{} \footnotesize  $I_{eff,a}$}}%
      \csname LTb\endcsname%
      \put(6212,3829){\makebox(0,0)[r]{\strut{} \footnotesize $I_{eff,p}$}}%
      \csname LTb\endcsname%
      \put(6212,3609){\makebox(0,0)[r]{\strut{}1}}%
    }%
    \gplbacktext
    \put(0,0){\includegraphics{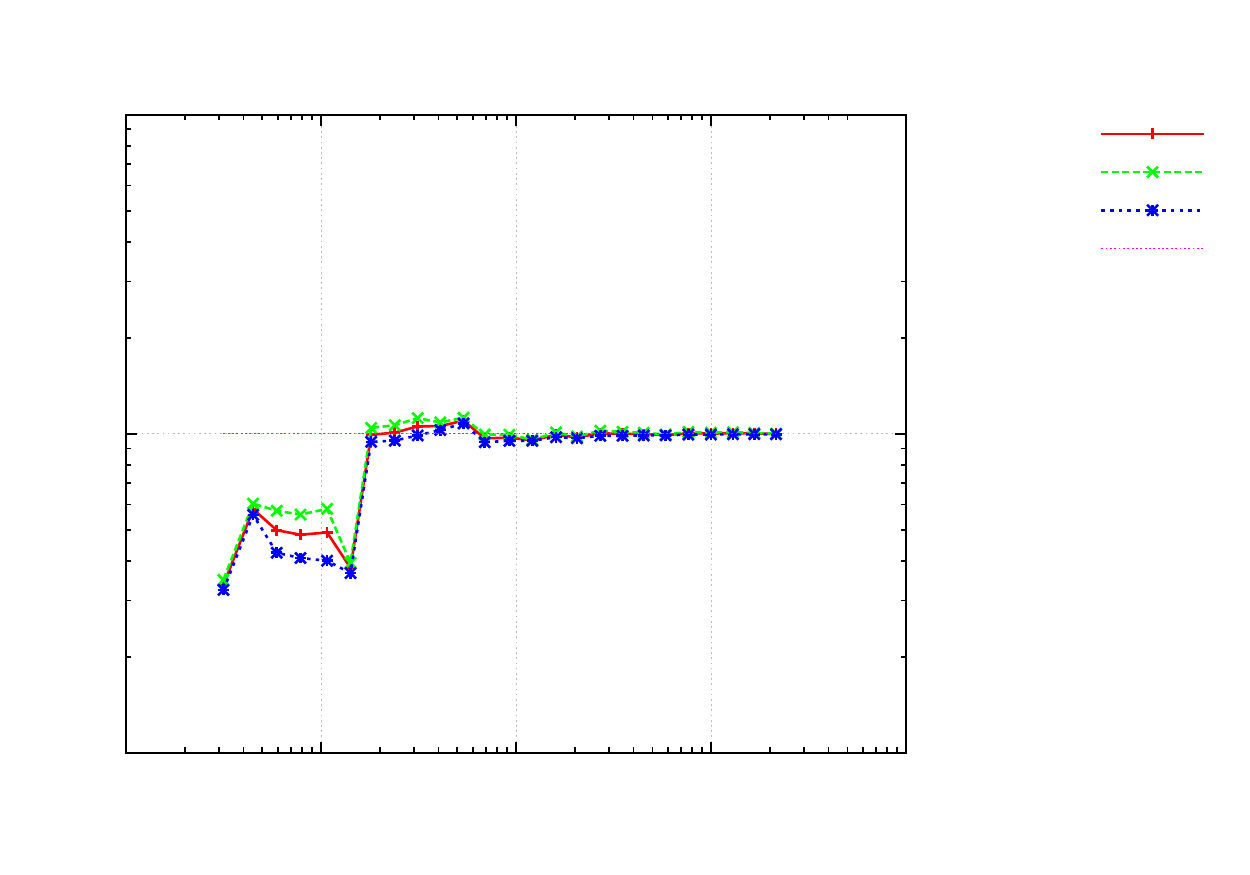}}%
    \gplfronttext
  \end{picture}%
\endgroup

%% file: Figures/YSplitterIeffConfig6.tex
\begingroup
  \makeatletter
  \providecommand\color[2][]{%
    \GenericError{(gnuplot) \space\space\space\@spaces}{%
      Package color not loaded in conjunction with
      terminal option `colourtext'%
    }{See the gnuplot documentation for explanation.%
    }{Either use 'blacktext' in gnuplot or load the package
      color.sty in LaTeX.}%
    \renewcommand\color[2][]{}%
  }%
  \providecommand\includegraphics[2][]{%
    \GenericError{(gnuplot) \space\space\space\@spaces}{%
      Package graphicx or graphics not loaded%
    }{See the gnuplot documentation for explanation.%
    }{The gnuplot epslatex terminal needs graphicx.sty or graphics.sty.}%
    \renewcommand\includegraphics[2][]{}%
  }%
  \providecommand\rotatebox[2]{#2}%
  \@ifundefined{ifGPcolor}{%
    \newif\ifGPcolor
    \GPcolortrue
  }{}%
  \@ifundefined{ifGPblacktext}{%
    \newif\ifGPblacktext
    \GPblacktexttrue
  }{}%
  \let\gplgaddtomacro\g@addto@macro
  \gdef\gplbacktext{}%
  \gdef\gplfronttext{}%
  \makeatother
  \ifGPblacktext
    \def\colorrgb#1{}%
    \def\colorgray#1{}%
  \else
    \ifGPcolor
      \def\colorrgb#1{\color[rgb]{#1}}%
      \def\colorgray#1{\color[gray]{#1}}%
      \expandafter\def\csname LTw\endcsname{\color{white}}%
      \expandafter\def\csname LTb\endcsname{\color{black}}%
      \expandafter\def\csname LTa\endcsname{\color{black}}%
      \expandafter\def\csname LT0\endcsname{\color[rgb]{1,0,0}}%
      \expandafter\def\csname LT1\endcsname{\color[rgb]{0,1,0}}%
      \expandafter\def\csname LT2\endcsname{\color[rgb]{0,0,1}}%
      \expandafter\def\csname LT3\endcsname{\color[rgb]{1,0,1}}%
      \expandafter\def\csname LT4\endcsname{\color[rgb]{0,1,1}}%
      \expandafter\def\csname LT5\endcsname{\color[rgb]{1,1,0}}%
      \expandafter\def\csname LT6\endcsname{\color[rgb]{0,0,0}}%
      \expandafter\def\csname LT7\endcsname{\color[rgb]{1,0.3,0}}%
      \expandafter\def\csname LT8\endcsname{\color[rgb]{0.5,0.5,0.5}}%
    \else
      \def\colorrgb#1{\color{black}}%
      \def\colorgray#1{\color[gray]{#1}}%
      \expandafter\def\csname LTw\endcsname{\color{white}}%
      \expandafter\def\csname LTb\endcsname{\color{black}}%
      \expandafter\def\csname LTa\endcsname{\color{black}}%
      \expandafter\def\csname LT0\endcsname{\color{black}}%
      \expandafter\def\csname LT1\endcsname{\color{black}}%
      \expandafter\def\csname LT2\endcsname{\color{black}}%
      \expandafter\def\csname LT3\endcsname{\color{black}}%
      \expandafter\def\csname LT4\endcsname{\color{black}}%
      \expandafter\def\csname LT5\endcsname{\color{black}}%
      \expandafter\def\csname LT6\endcsname{\color{black}}%
      \expandafter\def\csname LT7\endcsname{\color{black}}%
      \expandafter\def\csname LT8\endcsname{\color{black}}%
    \fi
  \fi
  \setlength{\unitlength}{0.0500bp}%
  \begin{picture}(7200.00,5040.00)%
    \gplgaddtomacro\gplbacktext{%
      \csname LTb\endcsname%
      \put(594,704){\makebox(0,0)[r]{\strut{}0.1}}%
      \csname LTb\endcsname%
      \put(594,2542){\makebox(0,0)[r]{\strut{}1}}%
      \csname LTb\endcsname%
      \put(594,4379){\makebox(0,0)[r]{\strut{}10}}%
      \csname LTb\endcsname%
      \put(726,484){\makebox(0,0){\strut{}100}}%
      \csname LTb\endcsname%
      \put(1849,484){\makebox(0,0){\strut{}1000}}%
      \csname LTb\endcsname%
      \put(2973,484){\makebox(0,0){\strut{}10000}}%
      \csname LTb\endcsname%
      \put(4096,484){\makebox(0,0){\strut{}100000}}%
      \csname LTb\endcsname%
      \put(5219,484){\makebox(0,0){\strut{}1e+06}}%
      \put(2972,154){\makebox(0,0){\strut{}\text{DOFs}}}%
      \put(2972,4709){\makebox(0,0){\strut{}Configuration 6}}%
    }%
    \gplgaddtomacro\gplfronttext{%
      \csname LTb\endcsname%
      \put(6212,4269){\makebox(0,0)[r]{\strut{} \footnotesize $I_{eff}$}}%
      \csname LTb\endcsname%
      \put(6212,4049){\makebox(0,0)[r]{\strut{} \footnotesize  $I_{eff,a}$}}%
      \csname LTb\endcsname%
      \put(6212,3829){\makebox(0,0)[r]{\strut{} \footnotesize $I_{eff,p}$}}%
      \csname LTb\endcsname%
      \put(6212,3609){\makebox(0,0)[r]{\strut{}1}}%
    }%
    \gplbacktext
    \put(0,0){\includegraphics{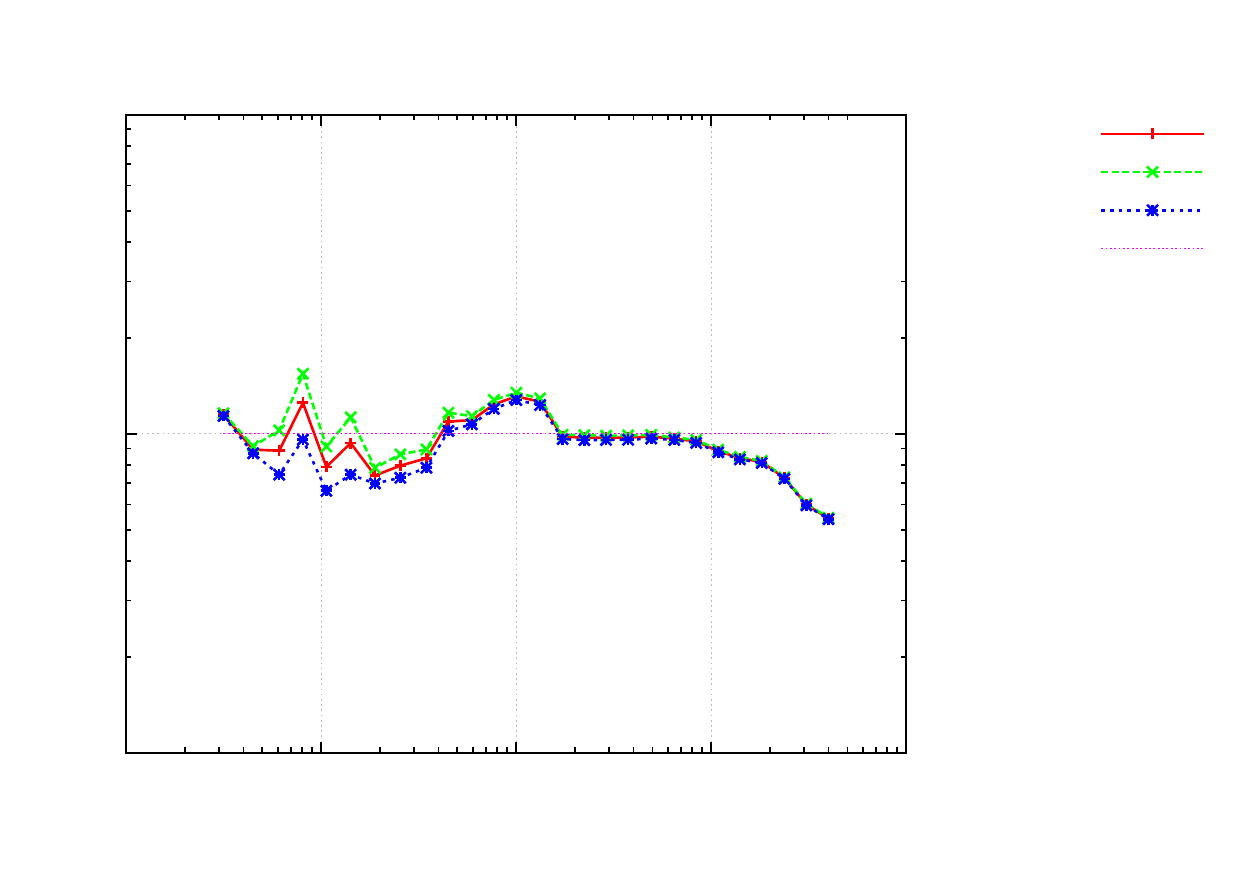}}%
    \gplfronttext
  \end{picture}%
\endgroup